\documentclass[12pt]{amsart}

\usepackage{amssymb, amscd, latexsym, color}
\usepackage{url}
\usepackage[final]{graphicx}

\newtheorem{thm}{Theorem}[section]
\newtheorem{prop}[thm]{Proposition}
\newtheorem{lem}[thm]{Lemma}
\newtheorem{cor}[thm]{Corollary}

\newtheorem{defn}[thm]{Definition}
\newtheorem{rmk}[thm]{Remark}

\newtheorem{notation}[thm]{Notation}

\begin{document}

\title[A basis of the fixed point subgroup of an automorphism ...]{An efficient algorithm for finding a basis\\ of the fixed point subgroup of an automorphism of a free group}

\author{Oleg Bogopolski$^{\ast}$}
\thanks{$^{\ast}$Mathematisches Institut, D\"{u}sseldorf University, D\"{u}sseldorf, Germany\\
{\sl and} Sobolev Institute of Mathematics of SB RAS, Novosibirsk, Russia.\\
e-mail: {\it bogopolski@math.uni-duesseldorf.de}}

\author{Olga Maslakova$^{\dag}$ }
\thanks{$^{\dag}$Sobolev Institute of Mathematics of SB RAS, Novosibirsk, Russia.\\
e-mail: {\it o.s.maslakova@gmail.com}}



\begin{abstract}
We prove that for any automorphism $\alpha$ of a free group $F$ of finite rank, one can efficiently compute
a basis of the fixed point subgroup ${\text{\rm Fix}}(\alpha)$.
\end{abstract}

\maketitle

\tableofcontents
\setcounter{tocdepth}{2}



\section{Introduction}
Let $F_n$ be the free group of finite rank $n$. For any endomorphism $\varphi$ of $F_n$ the {\it fixed point subgroup} of $\varphi$ is
$${\text{\rm Fix}}(\varphi)=\{x\in F_n\,|\, \varphi(x)=x\}.$$
More generally, if $S$ is a set of endomorphisms of $F_n$, the {\it fixed point subgroup} of $S$ is
${\text{\rm Fix}}(S)=\underset{\varphi\in S}{\bigcap}{\text{\rm Fix}}(\varphi)$.

Dyer and Scott showed in~\cite{DS} that
if $\varphi$ is an automorphism of finite order of $F_n$, then ${\text{\rm Fix}}(\varphi)$ is a free factor of $F_n$,
in particular the rank of ${\text{\rm Fix}}(\varphi)$ does not exceed~$n$. Scott later (1978) conjectured that $$\text{\rm rk}\,{\text{\rm Fix}}(\varphi)\leqslant n$$ for any automorphism of $F_n$.

Gersten~\cite{Gersten0, Gersten}, Goldstein and Turner~\cite{GT2}, and Cooper~\cite{Cooper} independently proved a weaker form of this conjecture that the group ${\text{\rm Fix}}(\varphi)$ is finitely generated.\break
In~\cite{Thomas}, Thomas generalized this result for an arbitrary set of automorphisms of~$F_n$.
In~\cite{St1}, Stallings proved that the equalizer subgroup
$${\text{\rm  Eq}}(\varphi,\psi):=\{x\in F_n\,|\, \varphi(x)=\psi(x)\}$$
is finitely generated for any two monomorphisms $\varphi,\psi$ from $F_n$ into a free group~$G$.
Goldstein and Turner~\cite{GT1} showed that the same holds in case where $\varphi$ is a homomorphism and $\psi$ is a monomorphism from $F_n$ to $F_n$.

In the seminal paper~\cite{BH}, Bestvina and Handel proved the Scott conjecture for automorphisms of free groups completely. The relative train track technique developed in this paper enabled
the solution of a sequel of difficult problems about automorphisms of free groups, see~\cite{BFH,BFH1,BFH2,CT,BG,Brinkmann1,Brinkmann0}.
Later Gaboriau, Levitt and Lustig~\cite{GLL} gave a dendrological proof of the Scott conjecture.
A stronger inequality, which also takes into
account infinite fixed words, is given in~\cite{GJLL}.
Another dendrological proof
was given by Sela in~\cite{Sela}.
A short dendrological proof of the Scott conjecture based on~papers~\cite{GLL},~\cite{Sela}, and~\cite{Paulin} is given
in the survey~\cite{Bestvina} of Bestvina on $\mathbb{R}$-trees.


Imrich and Turner~\cite{IT} proved that $\text{\rm rk}\,{\text{\rm Fix}}(\varphi)\leqslant n$ for an arbitrary endomorphism $\varphi$ of $F_n$.
Dicks and Ventura~\cite{DV} proved that,
$\text{\rm rk}\,{\text{\rm Fix}}(S)\leqslant n$ for any
set $S$ of injective endomorphisms of~$F_n$. Bergman, using derivations in group rings, proved that the same holds for an arbitrary set of endomorphisms of~$F_n$. However, the following problem has been open for almost 20~years.



\medskip

{\bf Problem A.} Find an efficient algorithm for computing a basis of ${\text{\rm Fix}}(\varphi)$,
where $\varphi$ is an automorphism of a free group $F$ of finite rank.

\medskip

Variations of this problem are formulated in the book of Dicks and Ventura~\cite[pages 71-72]{DV}; see a commentary there. A weaker form of this problem is formulated in~\cite[Problem (F1) (a)]{www}.








Problem A has been solved in three special cases: for positive automorphisms in the paper~\cite{CL} of Cohen and Lustig,
for special irreducible automorphisms in the paper of Turner~\cite[Proposition B]{Turner}, and for all automorphisms of $F_2$ in the paper of Bogopolski~\cite{B1}.

In 1999, Maslakova, a former PhD student of the first named author, attempted to solve this problem in general case.
However, her proof published in~\cite{Masl}, see also~\cite{MaslThes}, was not complete.
So, we have decided to give a full and correct proof.
The main result of this paper is the following.

\begin{thm}\label{thm 1.1}
{\it
Let $F_n$ be the free group of a finite rank $n$.
There exists an efficient algorithm which, given an automorphism $\varphi$ of $F_n$ finds
a basis of its fixed point subgroup
$${\text{\rm Fix}}(\varphi)=\{x\in F_n\,|\, \varphi(x)=x\}.$$
}
\end{thm}

We fix a basis $X=\{x_1,\dots,x_n\}$ of $F_n$. The length of a word $w\in F$ with respect to $X$ is denoted by $|w|$.
The {\it norm} of the automorphism $\varphi$ with respect to $X$ is the number $||\varphi||:=\max\{|\varphi(x_1)|,\dots ,|\varphi(x_n)|\}$.

The input of this algorithm is the sequence of words $\varphi(x_1),\dots,\varphi(x_n)$ in the alphabet
$X^{\pm}$; the output are words in $X^{\pm}$ which freely generate ${\text{\rm Fix}}(\varphi)$.

Under efficient algorithm we understand an algorithm for which
one can  write down a recursive function estimating the number of steps in terms of input.
We don't consider algorithms
which use two processes, one of which eventually terminates, and we don't use diagonalization methods.

As in~\cite{BH}, we use the relative train track techniques.
A relative train track is a homotopy equivalence $f:\Gamma \rightarrow \Gamma$
of a finite graph $\Gamma$ with certain good properties, see Section~4.
In~\cite[Theorem 5.12]{BH}), Bestvina and Handel proved that for any outer automorphism $\mathcal{O}$ of $F_n$, there
exists a relative train track $f:\Gamma \rightarrow \Gamma$ representing~$\mathcal{O}$.

However, to start our algorithm, we need to represent the automorphism $\varphi$ (and not its outer class)
by a relative train track $f:(\Gamma,v)\rightarrow (\Gamma,v)$, see Definition~\ref{defn 2.1}.
This is done in Theorem~\ref{thm 3.5}.

We use $f$ to define an auxiliary graph $D_f$ (first introduced in~\cite{GT1} in another setting, see also ~\cite{Turner}).
The fundamental group of one of the components of $D_f$, denoted $D_f(\mathbf{1}_{v})$, can be
identified with ${\text{\rm Fix}}(\varphi)$ (see Section~7). Thus, to compute a basis of
${\text{\rm Fix}}(\varphi)$, we need to construct the core $Core(D_f(\mathbf{1}_{v}))$ of this component.

At all but finite number of the vertices of $D_f$
there is a preferable outgoing direction. This determines a flow on almost all of $D_f$. The inverse automorphism
$\varphi^{-1}$ determines its own flow on almost all of $D_f$.
According to~\cite{Turner} (see also~\cite{CL}), there is a procedure for constructing a part of  $Core(D_f)$ which contains $Core(D_f(\mathbf{1}_{v}))$ if the latter is non-contractible:
one should start from a finite number of computable exceptional edges and follow the first flow for sufficiently long.
Theoretically we could arrive at a dead vertex, or get a loop, or arrive at a vertex where two rays of this flow meet, or none of these may occur.
To convert this procedure into an algorithm, we must detect at the beginning, which possibility occurs.
For that, we must solve the Finiteness and the Membership problems for vertices and certain subgraphs of $D_f$ (see Sections~3 and~7).
We solve these problems in this paper.\break
A sketch of the proof is given in Section~3.




\section{Preliminaries}
Let $\Gamma$ be a finite connected graph,
$\Gamma^0$ be the set of its vertices, $\Gamma^1$ be the set of its edges.
The initial vertex of an edge $E$ is denoted by $\alpha(E)$, the terminal by~$\omega(E)$, the inverse edge to $E$ is denoted by $\overline{E}$.

The geometric realization of $\Gamma$ is obtained by identification of each edge of $\Gamma$ with a real segment $[a,b]$ of length 1. This realization is denoted again by $\Gamma$. Using this realization,
we can work with partial edges and compute distances between points inside edges without passing to a subdivision.
{\it Partial edges} in $\Gamma$ are identified with subsegments $[a_1,b_1]\subset [a,b]$. Let $l$ be the corresponding metric on $\Gamma$.


For algorithmic aims, we work only with piecewise linear maps. For brevity, we  skip the wording
{\it piecewise linear}, e.g. we say a path instead of a piecewise linear path.



A {\it nontrivial path} in $\Gamma$ is a continuous map $\tau:[0,1]\rightarrow \Gamma$ with the following property:
there exist numbers $0=s_1<s_2<\dots <s_k<s_{k+1}=1$ and a sequence of (possibly partial) edges $E_1,E_2,\dots ,E_k$, such that $\tau|_{[s_i,s_{i+1}]}$ is a linear map onto $E_i$ for each $i=1,\dots,k$.
We will not usually distinguish between $\tau$ and the concatenation of (partial) edges $E_1E_2\dots E_k$. The {\it length} of $\tau$ is $l(\tau):=\sum_{i=1}^kl(E_i)$. For any two occurrences of points $u,v$ in $\tau$,
let $l_{\tau}(u,v)$ be the $l$-length of the subpath of $\tau$ connecting $u$ and~$v$.
A {\it trivial path} in $\Gamma$ is a map $\tau:[0,1]\rightarrow \Gamma$ whose image consists of a single point.
The trivial path whose image is $\{u\}$ is denoted by~$\mathbf{1}_u$; we set $l(\mathbf{1}_u)=0$.
An {\it edge-path} in $\Gamma$ is either a path of the form $E_1E_2\dots E_k$, where all $E_i$ are full edges, or a trivial path $\mathbf{1}_u$, where $u$ is a vertex.

The initial and the terminal points of a path $\tau$ are denoted by $\alpha(\tau)$ and $\omega(\tau)$, respectively.
The inverse path to $\tau$ is denoted by $\overline{\tau}$.
The concatenation of two paths, a reduced path, and homotopic paths are defined in a usual way.
By $[\tau]$ we denote the reduced path in $\Gamma$ which is homotopic to $\tau$ relative to the endpoints of $\tau$.
Let $[[\tau]]$ be the class of paths homotopic to $\tau$ relative to the endpoints of $\tau$.
For two paths $\tau,\mu$, we write $\tau=\mu$ if these paths are homotopic and $\tau\equiv \mu$ if they coincide.
The concatenation of $\tau$ and $\mu$ (if exists) is denoted by $\tau\mu$ or $\tau\cdot \mu$.

Let $\mathcal{PLHE}$ be the class of all homotopy equivalences $f:\Gamma\rightarrow \Gamma$ such that $\Gamma$ is a finite connected graph, $f(\Gamma^0)\subseteq \Gamma^0$,
and for each edge $E$ the following is satisfied:
$f(E)\equiv E_1E_2\dots E_k$, where each $E_i$ is an edge and
$E$ has a subdivision into segments, $E\equiv e_1e_2\dots e_k$, such that $f|{e_i}:e_i\rightarrow E_i$ is surjective and linear with respect to the metric $l$. The abbreviation $\mathcal{PLHE}$ stands for {\it piecewise linear homotopy equivalence}.


The homotopy equivalence $f:\Gamma\rightarrow \Gamma$ is called {\it tight} (resp. {\it nondegenerate})  if for each edge $E$ in $\Gamma$ the path $f(E)$ is reduced (resp. nontrivial).
The {\it norm} of $f$ is the number $||f||:=\max\{l(f(E))\,|\, E\hspace*{2mm}{\text{\rm is an edge of}}\hspace*{2mm}
\Gamma \}.$

Let $f:\Gamma\rightarrow \Gamma$ be a homotopy equivalence from $\mathcal{PLHE}$.
Then for every path $\tau$ in $\Gamma$ the map $f\circ \tau$ is also a path in~$\Gamma$.
We denote this path by~$f(\tau)$.
If $v$ is a distinguished vertex of $\Gamma$, we write $f_{\ast}:\pi_1(\Gamma,v)\rightarrow \pi_1(\Gamma,f(v))$
for the isomorphism given by $[[\tau]]\mapsto [[f(\tau)]]$, where $[[\tau ]]\in \pi_1(\Gamma,v)$.
We use the following rule for composition of two maps: $x(\varphi\psi)=(x\varphi)\psi$.

\vspace*{-1mm}

\section{Sketch of the proof}

In this section we give a sketch of the proof of the main theorem.
Some definitions we use here are given in the following sections.
At this point, it suffices to know the definition of a PL-relative train track (see Section 4).

\vspace*{1mm}

{\bf A.}  Let $\varphi$ be an automorphism of a free group $F$ of finite rank.

\vspace*{1mm}

\begin{defn}~\label{defn 2.1}
{\rm We say that $\varphi$ is {\it represented}\, by a homotopy equivalence~\!$f\!:~\!\!\Gamma\rightarrow~\!\Gamma$, where $\Gamma$ is a finite connected graph, if there is a vertex $v$ in $\Gamma$ fixed by $f$ and there is an isomorphism $j:F\rightarrow \pi_1(\Gamma,v)$
such that the automorphism $j^{-1}\varphi j$ of the group
$\pi_1(\Gamma,v)$ coincides with the induced automorphism~$f_{\ast}:\pi_1(\Gamma,v)\rightarrow \pi_1(\Gamma,v)$.
We also say that $f$ {\it represents $\varphi$}.}
\end{defn}


\vspace*{-2mm}

By Theorem~\ref{thm 3.5}, we may assume that $\varphi$ is algorithmically represented
by a homotopy equivalence $f:(\Gamma ,v) \rightarrow (\Gamma, v)$ which is a PL-relative train track.
Let $\varnothing =G_0\subset \dots \subset G_N=\Gamma$ be the fixed maximal filtration for $f$;
this gives us exponential, polynomial, or zero strata $H_i:={\text {\rm cl}}(G_i\backslash \,G_{i-1})$,
$i=1,\dots ,N$.


In Section~4 we recall that for each exponential stratum $H_r$, there exist an algebraic real number $\lambda_r>1$ and a pseudo-metric $L_r$ on $G_r$ with the following properties:

\begin{itemize}
\item if $E$ is an edge in $G_r$, then $L_r(E)>0$ if and only if $E$ is an edge in  $H_r$;
\item if $p$ is a path in $G_r$, then $L_r(f(p))=\lambda_r L_r(p)$;
\item if $p$ is a reduced $r$-legal path in $G_r$, then $L_r([f(p)])=~\lambda_r L_r(p)$.
\end{itemize}


{\bf B.} In Section~9, we define $r$-cancelation areas. To avoid technical details, we give here
an equivalent easier definition:
Let $H_r$ be an exponential stratum. Let $\tau$ be a reduced path in $G_r$.
An occurrence of a vertex $y$ in $\tau$ is called
an {\it $r$-cancelation point in} $\tau$ if $\tau$ contains a subpath $\bar{a}b$,
where $a$ and $b$ are nontrivial partial edges such that $\alpha(a)=\alpha(b)=y$
and the full edges containing $a$ and $b$ form an illegal $r$-turn.
A reduced path $\tau$ in $G_r$ is called an {\it $r$-cancelation area} if,
for each $k\in \mathbb{N}\cup \{0\}$, there is exactly one $r$-cancelation point in
$[f^k(\tau)]$ and if each proper subpath of $\tau$ does not have this property. These areas are important for
describing splittings of the so-called $r$-stable paths (see Section~10).


One can prove that the set of $r$-cancelation areas
coincides with the set $P_r$ defined before Lemma~4.2.5 in~\cite{BFH}. This lemma asserts
that the set $P_r$ is finite and that if $\tau\in P_r$, then $[f(\tau)]\in P_r$.
However, the proof of this lemma does not give a method for finding the elements of $P_r$.
In Section~9, we show that all $r$-cancelation areas can be efficiently found (Theorem~\ref{find_areas}).

The $r$-cancelation areas closely related to the indivisible periodic Nielsen paths in $G_r$ intersecting ${\text{\rm int}}(H_r)$. (Recall that a nontrivial reduced path $\sigma$ in $\Gamma$ is called a {\it periodic Nielsen path}
for $f:\Gamma\rightarrow \Gamma$ if $[f^k(\sigma)]=\sigma$ for some $k\geqslant 1$. A periodic Nielsen path $\sigma$ is called {\it indivisible} if it cannot be written as a concatenation of periodic Nielsen paths,
see Definition 5.1.1 in~\cite{BFH}.)
Indeed, if $N(r)$ is the number of $r$-cancelation areas for $f$, then the set
$\{[f^{N(r)}(A)]\,|\, A\hspace*{2mm}{\text{\rm  is an}}\hspace*{2mm} r{\text{\rm -cancelation area}}\}$
coincides with the set of all indivisible periodic Nielsen paths in $G_r$ intersecting ${\text{\rm int}}(H_r)$.

In Section~10 we define $r$-stable paths: A reduced path $\tau\subset G_r$ is called {\it $r$-stable} if the number
of $r$-cancelation points in $[f^k(\tau)]$ is the same for each $k\in \mathbb{N}\cup \{0\}$.
We prove that given a reduced edge path $\tau\subset G_r$,
one can efficiently check, whether $\tau$ is $r$-stable or not.
If not, we show how to efficiently find $i_0\in \mathbb{N}$
such that the path $[f^{i_0}(\tau)]$ is $r$-stable (Theorems~\ref{prop 7.2} and~\ref{prop 7.4}).

The following lemma is important for describing the splittings of paths in $G_r$.
A  variation of this lemma for $r$-stable paths is given in~\cite[Lemma~4.2.6]{BFH}. Another variation
is given in~\cite[Proposition~6.2]{Brinkmann0}.
However, the proof there, which is pretended to be constructive, is not correct. We give a correct proof in Section~11.










\medskip

{\bf Splitting lemma.}
{\it For any PL-relative train track $f:\Gamma\rightarrow \Gamma$, the following is satisfied:

Let $H_r$ be an exponential stratum of $\Gamma$ and
let $\tau$ be a reduced edge path in $G_r$. Then, for all $L>0$,  one can efficiently find an exponent $S>0$ such that at least one of the three possibilities occurs:

\begin{enumerate}
\item[{\rm 1)}] $[f^S(\tau)]$ contains an $r$-legal subpath of $r$-length greater than $L$.

\item[{\rm 2)}] $[f^S(\tau)]$ contains fewer illegal $r$-turns than $\tau$.

\item[{\rm 3)}] $[f^S(\tau)]$ is a trivial path or a concatenation of paths each of which is either an indivisible periodic Nielsen path
intersecting\, ${\text{\rm int}} (H_r)$ or an edge path in~$G_{r-1}$.
\end{enumerate}
}


In Section~12, we subdivide $\Gamma$ at the so called {\it $r$-exceptional points}
to obtain a new PL-relative train track $f':\Gamma'\rightarrow \Gamma'$ representing $\varphi$ and
satisfying the following additional condition:

\medskip

(RTT-iv) There is a computable natural number $P=P(f)$ such that for each exponential strata $H_r$
and each $r$-cancelation area $A$ of $f$, the $r$-cancelation area $[f^P(A)]$ is an edge path.

\medskip

Such relative train tracks are more convenient for solving algorithmic problems.

\medskip

{\bf C.} Let $\Gamma$ be a finite connected graph with a distinguished vertex $v$ and let $f:\Gamma\rightarrow \Gamma$ be a homotopy equivalence such that $f(\Gamma^0)\subseteq \Gamma^0$,
$f$ maps edges of $\Gamma$ to reduced edge paths,
and $f(v)=v$. Thus, $f$ is not necessarily a (subdivided) relative train track. In Section~7 and here we define a graph $D_f$ and describe a procedure which helps to
compute a basis of the group $$\overline{\text {\rm Fix}}(f):=\{[[p]]\in \pi_1(\Gamma,v)\,|\, f(p)=p\}.$$
In Sections~13-20, we will convert this procedure into an algorithm in case where $f$ is a subdivided relative train track.
From this we will deduce an algorithm for computing ${\text{\rm Fix}}(\varphi)$.

\medskip

{\bf A definition of the graph $D_f$.}
An {\it $f$-path} in $\Gamma$ is an edge path $\mu$ in $\Gamma$ (possibly trivial) such that $f$ maps the initial point of $\mu$ to the terminal point of~$\mu$.
Thus, if $\mu$ is an $f$-path in $\Gamma$, the path $\mu f(\mu)$ is well defined.
Moreover, if $\mu$ is an $f$-path and $E$ is an edge in $\Gamma$ satisfying $\alpha(E)=\alpha(\mu)$, then
$[\overline{E}\mu f(E)]$ is also an $f$-path.

\medskip

The vertices of the graph $D_f$ are reduced $f$-paths.
A vertex $\mu$ of $D_f$ is called {\it dead} if $\mu=\mathbf{1}_u$ for some vertex $u$ of $\Gamma$ fixed by $f$; otherwise $\mu$ is called {\it alive}.\break
Two vertices $\mu$ and $\tau$ in $D_f$ are connected by an edge (from $\mu$ to $\tau$)
with label $E$ if $E$ is an edge in $\Gamma$ satisfying $[\overline{E}\mu f(E)]=\tau$.

\medskip

We set $\widehat{f}(\mu):=[\overline{E}\mu f(E)]$ if $E$
is the first edge of the $f$-path $\mu$ (in particular, $\mu$ must be alive). Clearly, $\mu$ and $\widehat{f}(\mu)$, considered as vertices of $D_f$, are connected by an edge with the label $E$. The direction of this edge is called {\it preferable} at the vertex $\mu$ (see Figure~3).
Observe that there is a unique preferable direction at each alive vertex of $D_f$.
Preferable directions at all alive vertices determine a flow in $D_f$. Starting at $\mu$ and moving along this flow, we get the vertices
$\mu=\mu_1,\mu_2,\dots $, where  $\mu_{i+1}=\widehat{f}(\mu_i)$, $i\geqslant 1$.
These vertices together with the directed edges we pass
form a subgraph in $D_f$ which
we call the {\it $\mu$-subgraph}. The $\mu$-subgraph is either a finite segment, or a finite segment with a cycle, or a ray (see~Figure 5). We set $\widehat{f}^{\hspace*{1.5mm}0}(\mu):=\mu$ and $\widehat{f}^{\hspace*{1.5mm}i}(\mu):=\widehat{f}^{\hspace*{1.5mm}}(\widehat{f}^{\hspace*{1.5mm}{i-1}}(\mu))$
if $\widehat{f}^{\hspace*{1.5mm}{i-1}}(\mu)$ is defined and alive.

\medskip

An edge $e$ connecting two vertices $u,w$ in $D_f$ is called {\it repelling} if the direction of $e$
is not preferable at $u$ and the direction of the opposite edge $\overline{e}$ is not preferable at~$w$.
In Figure 1, the repelling edges are red.
Endpoints of repelling edges are called {\it repelling vertices}.
By Proposition~\ref{prop 4.1}, there exist only finitely many repelling edges in $D_f$ and they can be algorithmically found.

\medskip

{\bf D.} Let $v$ be a distinguished vertex of $\Gamma$ and $f(v)=v$.
Since $\mathbf{1}_{v}$ is an $f$-path, we can consider $\mathbf{1}_{v}$ as a vertex of~$D_f$. Let $D_f(\mathbf{1}_{v})$ be the component of $D_f$ containing $\mathbf{1}_{v}$. The fundamental group
$\pi_1(D_f(\mathbf{1}_{v}),\mathbf{1}_v)$ can be identified with $\overline{\text{\rm Fix}}(f)$,
see Lemma~\ref{lem 4.2}.

\vspace*{1mm}

A component of $D_f$ is called {\it repelling} if it contains at least one repelling edge.
Let $C_1,\dots ,C_n$ be all repelling components of $D_f$.
For each $C_i$, let $CoRe(C_i)$ be the minimal connected subgraph of $C_i$
which contains all repelling edges of $C_i$ and carries $\pi_1(C_i)$.
We set $C_f:=\cup_{i=1}^n C_i$ and $CoRe(C_f):=\cup_{i=1}^n CoRe(C_i)$.

By Proposition~\ref{prop 4.5}, to compute a basis of $\pi_1(D_f(\mathbf{1}_v),\mathbf{1}_v)$, it suffices
to construct $CoRe(C_f)$ and decide, whether the vertex $\mathbf{1}_v$ lies
in the $\mu$-subgraph for some repelling vertex $\mu$.










\hspace*{16mm}\hspace*{50mm}\includegraphics[scale=0.4]{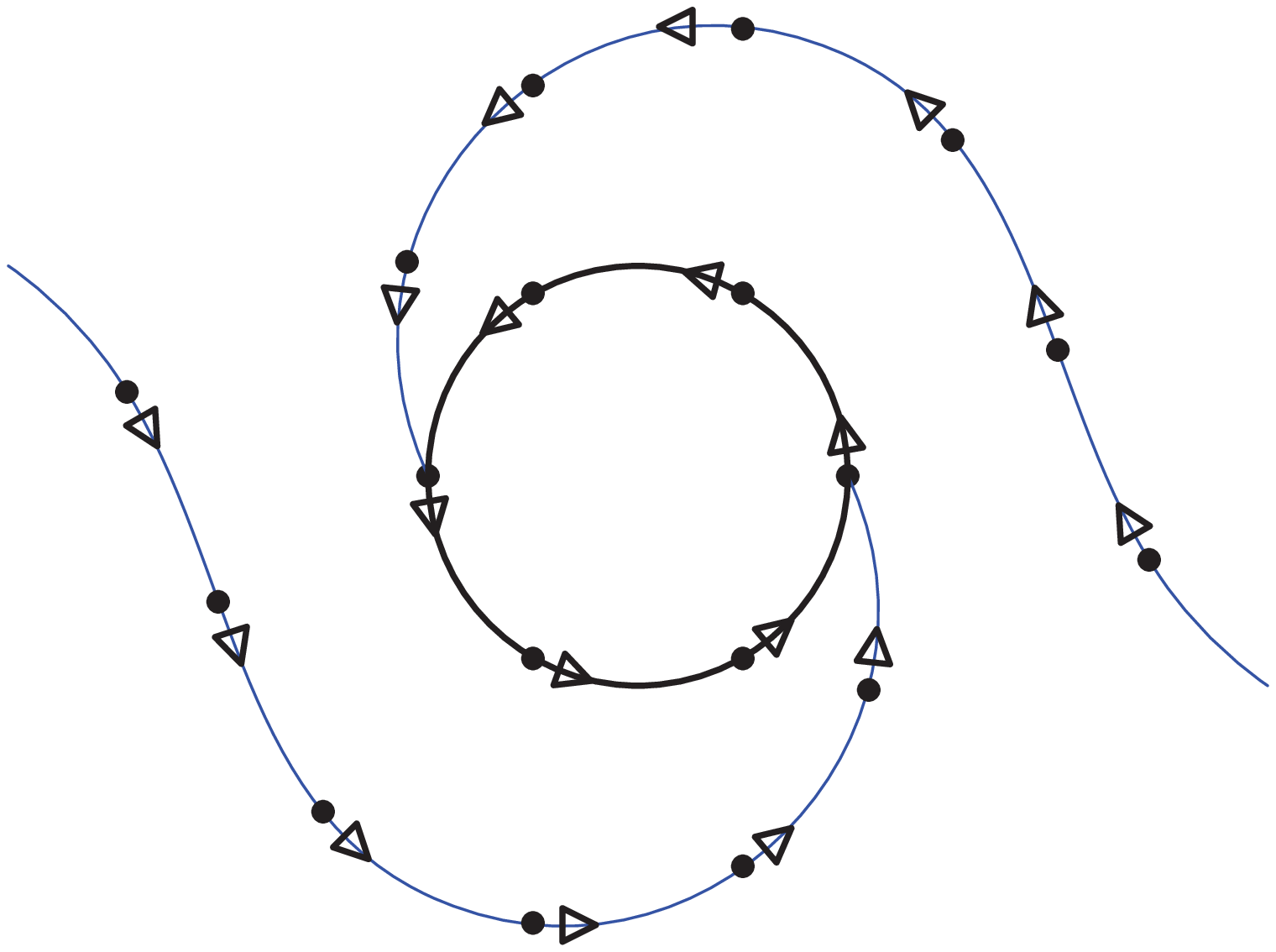}

\vspace*{-65mm}\hspace*{-5mm}\includegraphics[scale=0.5]{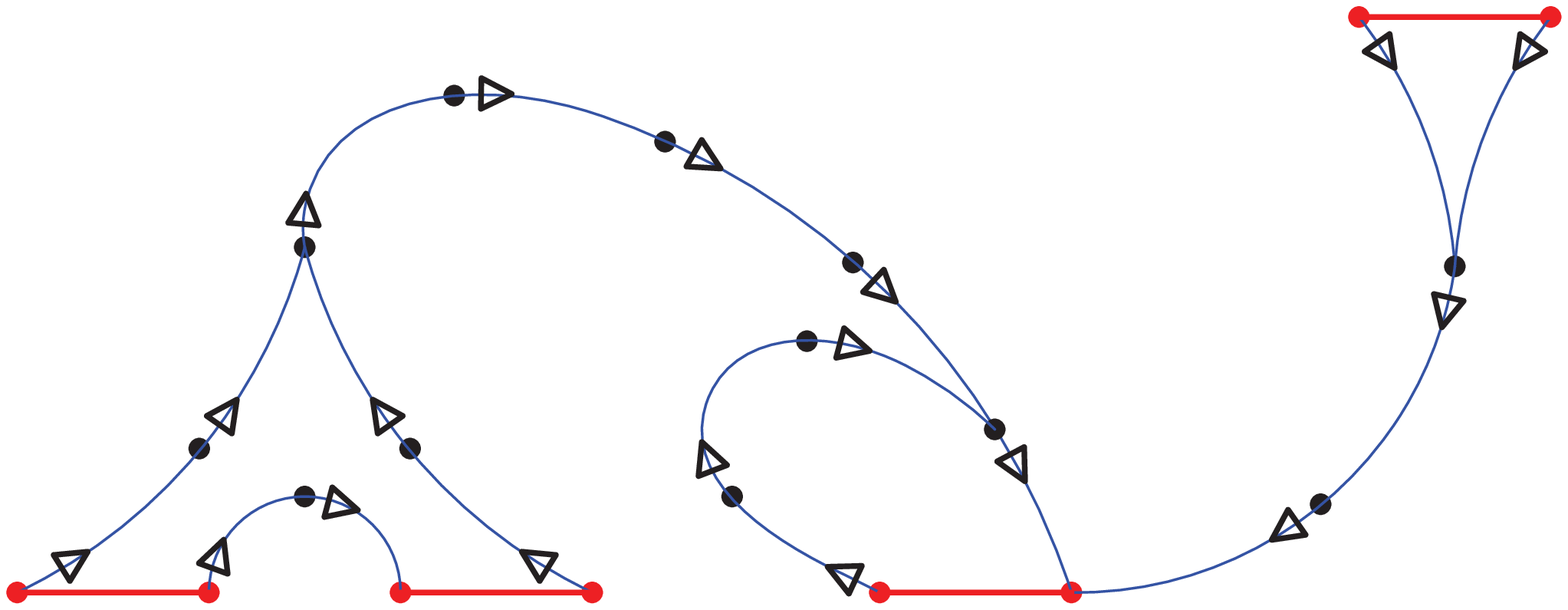}

\vspace*{-120mm}\hspace*{90mm}\includegraphics[scale=0.4]{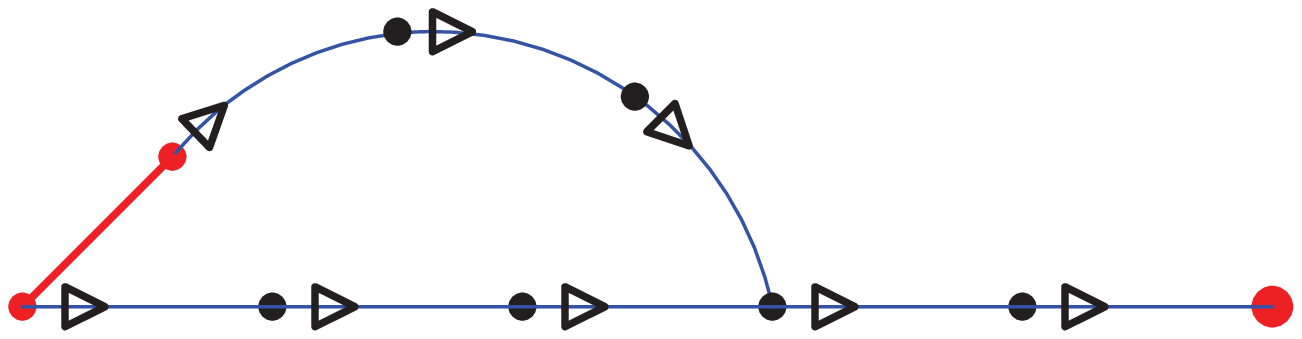}

\vspace*{-200mm}\hspace*{0mm}\includegraphics[scale=0.5]{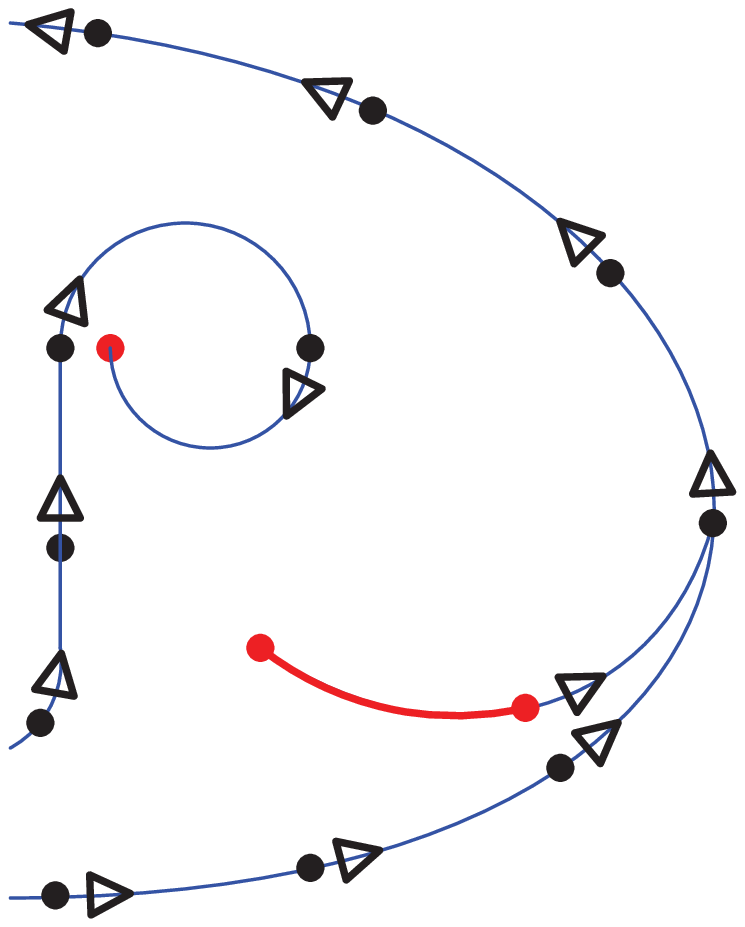}

\vspace*{-42mm}

\begin{center}
Figure 1.

\medskip

\begin{minipage}{9cm}
An example of a graph $D_f$ with three repelling components.
One of them is contractible.
\end{minipage}
\end{center}

\bigskip

{\bf E.}
It turns out that $CoRe(C_f)$ is contained in the union of the repelling edges and the $\mu$-subgraphs, where $\mu$ runs over the set of repelling vertices. So, to construct $CoRe(C_f)$, it suffices to do the following:



\begin{itemize}

\item[(1)] Find all repelling edges in $D_f$.

\item[(2)] For each alive repelling vertex $\mu$ determine, whether the $\mu$-subgraph is finite or~not.

\item[(3)] Compute all elements of all finite $\mu$-subgraphs from (2).

\item[(4)] For each two repelling vertices $\mu$ and $\tau$ with infinite $\mu$-and $\tau$-subgraphs determine,
whether these subgraphs intersect.

\item[(5)] If the $\mu$-subgraph and the $\tau$-subgraph from  (4) intersect,
find their first intersection point and compute
their initial segments up to this point.



\end{itemize}

As it was mentioned above, Step (1) can be done algorithmically.
In Section~7 we show that Steps (2)-(5) can be done algorithmically if the following two problems are solvable:

\medskip

{\bf Finiteness problem.} Given a vertex $\mu$ of the graph $D_f$, determine whether the $\mu$-sub\-graph is finite or not. If the $\mu$-subgraph is finite, construct it.

\medskip

{\bf Membership problem.}
Given two vertices $\mu$ and $\tau$ of the graph $D_f$, verify
whether $\tau$ is contained in the $\mu$-subgraph.

\medskip

In the part F, we sketch solutions of these problems in the case where $f$ is a PL-relative train track satisfying (RTT-iv).
Hence, for such $f$ we can compute a basis of $\overline{\text {\rm Fix}}(f)$.
By the parts A and B, we may assume that $\varphi$ is represented by such $f$.
Using this, we can compute a basis of ${\text{\rm Fix}}(\varphi)$.


\medskip

{\bf F.}
From now on and to the end of this section, we assume that
$f:(\Gamma,v)\rightarrow (\Gamma,v)$ is a PL-relative train track satisfying (RTT-iv).
Let $\varnothing =G_0\subset \dots \subset G_N=\Gamma$ be the maximal filtration associated with $f$.

To solve the Finiteness and the Membership problems, we investigate the $\mu$-sub\-graphs of $D_f$ in details.
Let $\mu\equiv E_1E_2\dots E_s$ be a nontrivial reduced $f$-path, and suppose that $\mu$ lies in $G_m$, but not in $G_{m-1}$.
There may be cancelations in computing $\widehat{f}(\mu)\equiv [E_2\dots E_s \cdot f(E_1)]$.

To control these cancelations,
we introduce the notions {\it $r$-perfect}, {\it $A$-perfect}, and {\it $E$-perfect} vertices in $D_f$, see Definitions~\ref{defn 9.1},~\ref{defn 9.3}, and~\ref{defn 13.1}, respectively.
These definitions imply that if $\mu$ is $r$-perfect or $A$-perfect, then there are no cancelations in passing from $\mu$ to $\widehat{f}(\mu)$, and if $\mu$ is $E$-perfect, then such cancelations are possible only between edges which lie in $G_{m-1}$.

Propositions~\ref{prop 12.3} and~\ref{prop 13.2}.(3) imply that, given a vertex $\mu$ in $D_f$, one can either
efficiently prove that the $\mu$-subgraph is finite, or find a perfect vertex $v_1\neq \mu$ in the $\mu$-subgraph.
In the second case we still have to decide, whether the $\mu$-subgraph is finite or not.

\medskip

{\it Case 1.} Suppose that $v_1$ is $r$-perfect.

Let $\lambda_r>1$ be the Perron-Frobenius eigenvalue corresponding to $H_r$. Then, by Lemma~\ref{lem 14.1}
we have the following:


\begin{enumerate}
\item[(1)] $L_r(\widehat{f}^{\hspace*{1.5mm}i+1}(v_1))\geqslant L_r(\widehat{f}^{\hspace*{1.5mm}i}(v_1))>0$ for all $i\geqslant 0$.
\vspace*{1mm}

\item[(2)] There exist computable natural numbers $m_1<m_2<\dots $, such that\\
$L_r(\widehat{f}^{\hspace*{1.5mm}m_i}(v_1))= \lambda_r^i L_r(v_1)$ for all $i\geqslant 1$.
\end{enumerate}


Hence, the $v_1$-subgraph is infinite, and so the $\mu$-subgraph is infinite, that solves the Finiteness problem in this case. Also the Membership problem for the $v_1$-subgraph, and hence for the $\mu$-subgraph, can be easily solved with the help of~$L_r$.

\medskip

{\it Case 2.}
Suppose that $v_1$ is $A$-perfect.

By Proposition~\ref{prop 9.4}, we can efficiently find a
finite set $\{v_1,v_2,\dots ,v_k\}$ of $A$-perfect vertices in the $v_1$-subgraph
such that the following holds:

\begin{enumerate}

\item[(1)] All $A$-perfect vertices in the $v_1$-subgraph are $[f^{\hspace*{1.5mm}i}(v_j)]$,
$i\geqslant 0$, $1\leqslant j\leqslant k$.

\item[(2)] $[f^{\hspace*{1.5mm}i}(v_j)]=\widehat{f}^{\hspace*{1.5mm}m_{i,j}}(v_1)$ for some computable $m_{i,j}$. Moreover, $m_{i,j}$ are different for different $(i,j)$.

\item[(3)] Given a vertex $u$ in the $v_1$-subgraph, we can find $\ell$
such that $\widehat{f}^{\hspace*{1mm} \ell }(u)$
is an $A$-perfect vertex. The least $\ell$ does not exceed the number of edges in the $f$-path~$u$.

\end{enumerate}





Using this, the Finiteness and the Membership problems for the $v_1$-subgraph can be reduced~to the following problems:

\medskip

{\tt FIN\,$(v_1)$.} Do there exist $p>q\geqslant 0$ such that
$[f^p(v_1)]=[f^q(v_1)]$?

\medskip

{\tt MEM\,$(v_1)$.} Given a reduced $f$-path $\tau$, does there exist
$p\geqslant 0$ such that $[f^p(v_j)]=\tau$ for some $1\leqslant j\leqslant k$?

\medskip
\noindent
These problems are solvable by Corollaries~\ref{cor 3.9} and~\ref{cor 3.8}, which we
deduce from a theorem of Brinkmann, see~\cite[Theorem~0.1]{Brinkmann1}.

\medskip

{\it Case 3.} Suppose that $v_1$ is $E$-perfect.

Then the solution is similar to that given in Case~2, see Sections~19 and~20.






\medskip

\section{Relative train tracks for outer automorphisms of free groups}

First we recall the definition of a relative train track from~\cite{BH}. Since we are interested in algorithmic
problems, we will work only with homotopy equivalences from the class $\mathcal{PLHE}$ defined in Section~2.

Let $\Gamma$ be a finite connected graph and let $f:\Gamma\rightarrow \Gamma$ be a tight and nondegenerate
homotopy equivalence from the class $\mathcal{PLHE}$.


A {\it turn} in $\Gamma$ is an unordered pair of
edges of $\Gamma$ originating at a common vertex. A turn is {\it nondegenerate} if these edges
are distinct, and it is {\it degenerate} otherwise.
The map $f:\Gamma\rightarrow \Gamma$ induces a map $Df:\Gamma^1\rightarrow \Gamma^1$ which sends each edge $E\in\Gamma^1$ to the
first edge of the path $f(E)$. This induces a map $Tf$ on turns in $\Gamma$ by the rule $Tf(E_1,E_2)=(Df(E_1),Df(E_2))$. A turn $(E_1,E_2)$ is {\it legal} if the turns $(Tf)^n(E_1,E_2)$
are nondegenerate for all $n \geqslant  0$; a turn is {\it illegal} if it is not legal.
An edge path $E_1E_2\dots E_m$ in $\Gamma$ is {\it legal} if all its turns $(\overline{E}_i,E_{i+1})$ are legal.
Clearly, a legal edge path is reduced.

From each pair of mutually inverse edges of the graph $\Gamma$ we choose one edge. Let
$\{E_1,\dots,E_k\}$ be the ordered set of chosen edges. The {\it transition matrix} of the
map $f$ (with respect to this ordering) is the matrix $M(f)$ of the size $k\times k$ such that the $i\!j^{\text{\rm th}}$ entry of $M(f)$ is equal
to the total number of occurrences of the edges $E_i$ and $\overline{E}_i$ in the path $f(E_j)$.

\medskip

A {\it filtration} for $f:\Gamma\rightarrow \Gamma$ is an increasing sequence of (not necessarily connected) $f$-invariant subgraphs $\varnothing =G_0\subset \dots \subset G_N=\Gamma$. The subgraph $H_i={\text {\rm cl}}(G_i\backslash \,G_{i-1})$
is called the {\it $i$-th stratum}. Edges in $H_i$ are called {\it $i$-edges}. A turn with both edges in $H_i$ is called an {\it $i$-turn}. A turn with one edge in $H_i$ and another in $G_{i-1}$ is called
{\it mixed} in $(G_i,G_{i-1})$. We assume that the edges of $\Gamma$ are ordered so that the edges from $H_i$
precede the edges from $H_{i+1}$. The edges from $H_i$ define a square submatrix $M_{[i]}$ of $M(f)$.

If the filtration is maximal, then each matrix $M_{[i]}$ is irreducible.
If $M_{[i]}$ is nonzero and irreducible,
then it has the associated Perron-Frobenius eigenvalue $\lambda_i\geqslant 1$.
If $\lambda_i>1$, then the stratum $H_i$ is called {\it exponential}.
If $\lambda_i=1$, then $H_i$ is called {\it polynomial}.
In this case $M_{[i]}$ is a permutation matrix, hence for every edge $E\in H_i^1$ the path $f(E)$ contains exactly one edge of $H_i$, all other edges of $f(E)$ lie in $G_{r-1}$.
A stratum $H_i$ is called a {\it zero stratum} if $M_{[i]}$ is a zero matrix. In this case
$f(E)$ lies in $G_{i-1}$ for every edge $E\in H_i^1$.


\vspace*{-0.2mm}

\begin{defn}
{\rm
Let $\Gamma$ be a finite connected graph and let $f:\Gamma\rightarrow \Gamma$ be a tight and nondegenerate
homotopy equivalence from the class $\mathcal{PLHE}$.
The map $f$ is called a {\it PL-relative train track} if there
exists a maximal filtration $\varnothing =G_0\subset \dots \subset G_N=\Gamma$ for $f$ such that each exponential stratum $H_r$ of this filtration satisfies the following conditions:

\begin{enumerate}
\item[(RTT-i)] $Df$ maps the set of edges of $H_r$ to itself; in particular all mixed turns in
$(G_r,G_{r-1})$ are legal.
\item[(RTT-ii)] If $\rho\subset G_{r-1}$ is a reduced nontrivial edge path with endpoints in  $H_r\cap\,G_{r-1}$,
then $[f(\rho)]$ is a  nontrivial edge path with endpoints in $H_r\cap G_{r-1}$.
\item[(RTT-iii)] For each legal edge path $\rho\subset H_r$, the path $f(\rho)$ does not contain
any illegal turns in $H_r$.

\end{enumerate}
}
\end{defn}






\begin{defn} {\rm We use the above notations. Let $H_r$ be an exponential stratum.\break
A nontrivial reduced path $\rho $ in $G_r$ is called {\it $r$-legal}
if the minimal edge path containing $\rho$ does not contain any illegal turns in~$H_r$.


}
\end{defn}


The following proposition will be often used in the further proof.

\begin{prop}\label{prop 3.1} {\rm (see~\cite[Lemma 5.8]{BH})}
Suppose that $f:\Gamma\rightarrow \Gamma$ is a relative train track and
$H_r$ is an exponential stratum of $\Gamma$.
Let $\rho$ be a reduced $r$-legal path:
$$\rho\equiv b_0\cdot a_1\cdot b_1\cdot \ldots \cdot a_k\cdot b_k,$$
where $k\geqslant 1$, $a_1,\dots,a_k$ are paths in $H_r$, and $b_0,\dots ,b_k$ are paths in $G_{r-1}$,
and all these paths except maybe $b_0$ and $b_k$ are nontrivial. Then
$$[f(\rho)]\equiv [f(b_0)]\cdot f(a_1)\cdot [f(b_1)]\cdot \ldots \cdot f(a_k)\cdot [f(b_k)]$$
and this path is $r$-legal. Moreover, for all $i\geqslant 1$ we have
$$[f^i(\rho)]\equiv [f^i(b_0)]\cdot [f^i(a_1)]\cdot [f^i(b_1)]\cdot \ldots \cdot [f^i(a_k)]\cdot [f^i(b_k)]$$
and these paths are
$r$-legal.
\end{prop}


{\bf The $r$-length function $L_r$.}
Let $f:\Gamma\rightarrow \Gamma$ be a PL-relative train track and
$H_r$ be an exponential stratum. Choose a positive vector $\vec{v}$
satisfying $\vec{v} M_{[r]}=\lambda_r\vec{v}$.
Since $M_{[r]}$ is an integer matrix, we can choose~$\vec{v}$ so that the coordinates of $\vec{v}$ are rational functions of $\lambda_r$ over~$\mathbb{Q}$. If $E_i$ is the $i^{\text {\rm th}}$ edge of $H_r$,
define $L_r(E_i)=v_i$; if $E$ is an edge of $G_{r-1}$, define $L_r(E)=0$.
For an arbitrary edge path $\tau$ in $G_r$, we define its {\it $r$-length} $L_r(\tau)$ as the sum
of $r$-lengths of edges of $\tau$.
Then we have $L_r(f(E))=\lambda_r L_r(E)$.

We extend this definition to all paths (not necessarily edge paths) in $G_r$, as it was done in Lemma~5.10 in \cite{BH}.
For an arbitrary path $\mu$ in $G_r$, let $L^{\bullet}_r (\mu)$ be the sum of $L_r$-lengths of full $r$-edges
which occur in $\mu$ if they exist and zero if not.
For any path $\rho$ in $G_r$, we set $$L_r(\rho):=\underset{k\rightarrow \infty}{\lim} \lambda_r^{-k} L^{\bullet}_r(f^k(\rho)).$$


\begin{lem}\label{zero_length}
The function $L_r$ has the following properties:
\begin{enumerate}
\item[{\rm 1)}] $L_r(f(\rho))=\lambda_r L_r(\rho)$ for any path $\rho$ in $G_r$.

\item[{\rm 2)}] $L_r([f(\rho)])=\lambda_r L_r(\rho)$ for any reduced $r$-legal path $\rho$ in $G_r$.

\item[{\rm 3)}] If $\rho$ is a nontrivial initial or terminal segment of an $r$-edge, then $L_r(\rho)>0$.

\item[{\rm 4)}] If $\rho$ is a nontrivial segment of an $r$-edge, then
there exists $k\in \mathbb{N}$
such that $f^k(\rho)$ does not lie in an $r$-edge.

\item[{\rm 5)}] If $\rho$ is a nontrivial path in $G_r$ with $L_r(\rho)=0$,
then there exists $k\in \mathbb{N}$ such that $f^k(\rho)$ lies in $G_{r-1}$.
\end{enumerate}
\end{lem}

{\it Proof.} 1) follows from the definition of $L_r$, 2) from Proposition~\ref{prop 3.1}, and 3) from (RTT-i)
and 4). We prove 4). For that we use the following claims:

\begin{enumerate}

\item[i)] For each $k\in \mathbb{N}$, the map $f^k$ restricted to each component of $\Gamma \setminus f^{-k}(\Gamma^0)$ is linear with respect to the metric $l$.

\item[ii)] Let $k_0$ be the number of $r$-edges in $H_r$ plus 1. Then for each $r$-edge $E$ the path $f^{k_0}(E)$ contains at least two $r$-edges.
\end{enumerate}

The first claim follows from the assumption that $f$ lies in the class $\mathcal{PLHE}$, the second one from the assumption that
the stratum $H_r$ is irreducible and exponential.

Below we define a number $0<a<1$ satisfying the following property:
if $\rho$ is a nontrivial segment of an $r$-edge and $f^{k_0}(\rho)$ lies in an $r$-edge, then $l(f^{k_0}(\rho))\geqslant l(\rho)/a$.

Let $E$ be an $r$-edge and suppose that $f^{k_0}(E)\equiv E_1b_1E_2\dots b_{s-1}E_s$, where $E_1,\dots ,E_s$ are $r$-edges and $b_1,\dots,b_{s-1}$ are paths in $G_{r-1}$ or trivial.
Write $E\equiv E'_1b'_1E'_2\dots b'_{s-1}E_s'$, where $f^{k_0}(E'_i)\equiv E_i$ and $f^{k_0}(b'_i)\equiv b_i$. Since $s\geqslant 2$,
the number $$a_E:=\max \{l(E'_i)\,|\ i=1,\dots ,s\}$$ is smaller than 1. Let $a$ be the maximum of $a_E$ over all
$r$-edges $E$. Then $a$ has the desired property.

To complete 4), we take the minimal $m\in \mathbb{N}$ with $l(\rho)> a^m$. Then $f^{k_0m}(\rho)$ does not lie in an $r$-edge.

Now we prove 5). Since $L_r(\rho)=0$, the statement 3) implies that $\rho$ lies either in $G_{r-1}$,
or in the interior of an $r$-edge. In the first case we are done. In the second case, by 4), there exists $k\in\mathbb{N}$ such that $f^k(\rho)$ does not lie in an $r$-edge. Again 3) implies that $f^k(\rho)$ lies in $G_{r-1}$.
\hfill$\Box$



\medskip

\noindent
{\bf A representation of an outer automorphism of $F_n$ by a PL-relative train track.}
The rose with $n$ petals $R_n$ is the graph with one vertex $\ast$ and $n$ geometric edges. We assume that the free
group on $n$ letters $F_n$ is identified with $\pi_1(R_n, \ast)$. Obviously, every automorphism $\varphi$ of $F_n$ can be represented by a homotopy equivalence $R_n\rightarrow R_n$.

In~\cite[Theorem 5.12]{BH}), Bestvina and Handel proved that every outer automorphism $\mathcal{O}$ of $F_n$
can be represented by a relative train track $f:\Gamma\rightarrow \Gamma$.
One can show that this proof can be organized in a constructive way. Also, we may assume that $f$ is a PL-relative train track.
Thus, we have the following start-point for our algorithm.

\begin{thm}\label{Bestvina} {\rm (see~\cite[Theorem 5.12]{BH})}
{\it Let $F_n$ be the free group of finite rank~$n$. There is an efficient algorithm which, given
an outer automorphism $\mathcal{O}$ of $F$,
constructs a PL-relative train track $f:\Gamma \rightarrow \Gamma$ and a homotopy equivalence (a marking) $\tau: R_n\rightarrow \Gamma$
such that $f$ represents $\mathcal{O}$ with respect to $\tau$.
}
\end{thm}

The latter means that for any homotopy equivalence $\sigma:\Gamma\rightarrow R_n$ which is a homotopy inverse to $\tau$,
the map $(\tau f \sigma)_{\ast}:\pi_1(R_n,\ast)\rightarrow \pi_1(R_n,\ast)$ represents $\mathcal{O}$.
Below we give a reformulation of this theorem, see Theorem~\ref{thm 3.2a}.

\section{Relative train tracks for automorphisms of free groups}

Let $F$ be a free group of finite rank, $\varphi$ be an automorphism of $F$, and $\mathcal{O}$ be the outer
automorphism class of $\varphi$.
Theorem~\ref{Bestvina} gives a representation of $\mathcal{O}$ by a PL-relative train track. However this is not sufficient for our aims.
The purpose of this section is to show that $\varphi$ itself can be  represented by a PL-relative train track, see
Theorem~\ref{thm 3.5}.

\begin{notation}
{\rm Let $\Gamma$ be a finite connected graph and $f:\Gamma\rightarrow \Gamma$ be a homotopy equivalence.
For each vertex $v\in \Gamma^0$ we define the isomorphism
$$
\begin{array}{rl}
f_{v}: \pi_1(\Gamma,v) & \rightarrow \pi_1(\Gamma,f(v)),\vspace*{2mm}\\
\,\! [[\mu]] & \mapsto [[f(\mu)]], \hspace*{2mm}
{\text{\rm where}}\hspace*{2mm} [[\mu]]\in \pi_1(\Gamma,v).
\end{array}
$$
\noindent
For each path $p$ in $\Gamma$ from $v$ to $f(v)$ we define the automorphism
$$
\begin{array}{rl}
f_{v,p}: \pi_1(\Gamma,v) & \rightarrow \pi_1(\Gamma,v),\vspace*{2mm}\\
\,\! [[\mu]] & \mapsto [[p f(\mu) \bar{p}]], \hspace*{2mm}
{\text{\rm where}}\hspace*{2mm} [[\mu]]\in \pi_1(\Gamma,v).
\end{array}
$$
}
\end{notation}


\begin{rmk}
{\rm By Theorem~\ref{Bestvina}, given an automorphism $\varphi$ of $F$,
one can construct a finite connected graph $\Gamma$,
a PL-relative train track $f:\Gamma\rightarrow \Gamma$, and an isomorphism $i:F\rightarrow \pi_1(\Gamma,v)$,
where $v$ is a vertex of $\Gamma$, such that the automorphism $i^{-1}\varphi\, i:\pi_1(\Gamma,v)\rightarrow \pi_1(\Gamma,v)$ coincides with $f_{v,p}$ for an appropriate path $p\subset \Gamma$ from $v$ to~$f(v)$.

We claim that $p$ can be computed. Indeed, if $q$ is an arbitrary path in $\Gamma$ from $v$ to $f(v)$, then $f_{v,q}$
differs from $f_{v,p}$ by an inner automorphism of $\pi_1(\Gamma,v)$.
Comparing $f_{v,q}$ with $i^{-1}\varphi\, i$, we can compute this inner automorphism and hence $p$.
This gives us the following form of Theorem~\ref{Bestvina}.
}

\end{rmk}

\begin{thm}\label{thm 3.2a}
Let $F$ be a free group of finite rank. There is an efficient algorithm which, given
an automorphism $\varphi$ of $F$, constructs
a PL-relative train track $f:\Gamma\rightarrow \Gamma$ and indicates a vertex $v\in \Gamma^0$,
a path $p\subset \Gamma$ from $v$ to $f(v)$, and an isomorphism $i:F\rightarrow \pi_1(\Gamma,v)$
such that the automorphism $i^{-1}\varphi\, i:\pi_1(\Gamma,v)\rightarrow \pi_1(\Gamma,v)$ coincides with $f_{v,p}$.
\end{thm}

The following theorem says that in Theorem~\ref{thm 3.2a} we can provide $f(v)=v$
and choose $p$ equal to the trivial path at $v$.

\medskip

\begin{thm}~\label{thm 3.5} Let $F$ be a free group of finite rank.
There is an efficient algorithm which, given an automorphism $\varphi$ of $F$, constructs
a PL-relative train track $f_1:\Gamma_1\rightarrow \Gamma_1$ with a vertex $v_1\in \Gamma_1^0$ fixed by $f_1$, and indicates an isomorphism $j:F\rightarrow \pi_1(\Gamma_1,v_1)$ such that $j^{-1}\varphi j=(f_1)_{v_1}$.
\end{thm}

{\it Proof.}
Let $f:\Gamma \rightarrow \Gamma$, $v$, $p$, and $i:F\rightarrow \pi_1(\Gamma,v)$ be the $PL$-relative train track, the vertex, the path,
and the isomorphism from Theorem~\ref{thm 3.2a}, respectively. Then we have $\varphi\, i=i\, f_{v,p}$. Hence, for every $w\in F$, we have
$$
i(\varphi(w))=[[p]]\, [[f (i(w))]]\, [[\bar{p}]].
\eqno{(5.1)} $$

Let $\Gamma_1$ be the graph obtained from $\Gamma$ by adding a new vertex $v_1$ and a new edge $E$
connecting $v_1$ and $f(v)$. We extend the homotopy equivalence $f:\Gamma\rightarrow \Gamma$ to a map
$f_1:\Gamma_1\rightarrow \Gamma_1$  by the rule $f_1(v_1)=v_1$ and $f_1(E):=Ef(p)$. Clearly, $f_1$ is a homotopy
equivalence.
We define a maximal filtration for $f_1$ by extending the maximal filtration for $f$ with the help of the new top polynomial stratum
consisting of the edges $E$ and $\bar{E}$.
Finally, we define the isomorphism $j:F\rightarrow \pi_1(\Gamma_1,v_1)$ by the rule
$$
j(w):=[[E]]\,[[f(i(w))]]\,[[\bar{E}]],\hspace{5mm} w\in F.\eqno{(5.2)}$$ To complete the proof, we verify that the automorphism
$j^{-1}\varphi\, j$ of the group $\pi_1(\Gamma_1,v_1)$ coincides with the induced automorphism $(f_1)_{v_1}:\pi_1(\Gamma_1,v_1)\rightarrow \pi_1(\Gamma_1,v_1)$.
It suffices to check that $(f_1)_{v_1}(j(w))=j(\varphi(w))$ for any $w\in F$:
$$(f_1)_{v_1}(j(w))\,\overset{(5.2)}{=}\,(f_1)_{v_1}\bigl( [[E]]\, [[f(i(w))]]\, [[\bar{E}]]\bigr)=
[[f_1(E)]]\, [[f^2(i(w))]]\, [[f_1(\bar{E})]]=$$
$$[[Ef(p)]]\, [[f^2(i(w))]]\, [[f(\bar{p})\bar{E}]]=
[[E]]\, [[f\bigl(p\, f(i(w))\,\bar{p}\bigr)]]\, [[\bar{E}]]\,\overset{(5.1)}{=}$$ $$\,[[E]]\,[[f(i(\varphi(w)))]]\,[[\bar{E}]]\,\overset{(5.2)}{=}\,j(\varphi(w)). $$
\hfill $\Box$


Thus, for computing a basis of ${\text{\rm Fix}}(\varphi)$, it suffices to compute a basis of the group
$$\overline{{\text{\rm Fix}}}(f_1)=\{[[\mu]]\in \pi_1(\Gamma_1,v_1)\,|\, f_1(\mu)=\mu\},$$
where $\Gamma_1$ is the graph, $v_1$ is the vertex, and $f_1$ is the PL-relative train track from Theorem~\ref{thm 3.5}.
In Section~12, we will show that we may assume that $f_1$ satisfies (RTT-iv).

\section{Auxiliary statements}

\medskip

Let $F$ be a free group of finite rank with a fixed basis $X$. For any element $w\in F$
let $|w|$ be the length of $w$ with respect to $X$.

The following theorem was proven by P.~Brinkmann in~\cite[Theorem 0.1]{Brinkmann1}.

\medskip

\begin{thm}~\label{thm 3.7} There exists an efficient algorithm which, given an automorphism $\varphi$
of a free group $F$ of finite rank and given elements $u,v\in F$, verifies,
whether there exists a natural $N$ such that $\varphi^N(u)=v$.
If such $N$ exists, then the algorithm computes $N$ as well.

\end{thm}

\medskip

\begin{cor}~\label{cor 3.8} There exists an efficient algorithm which, given a finite connected graph $\Gamma$ and
a homotopy equivalence $f:\Gamma\rightarrow \Gamma$ with $f(\Gamma^0)\subseteq \Gamma^0$, and
given two edge paths $\rho, \tau$ in $\Gamma$, decides whether there exists a natural
number $k$ such that $f^k(\rho)=\tau$. If such $k$ exists, then the algorithm computes it.
\end{cor}


{\it Proof.} First we reduce the problem to the case, where $f$ fixes the endpoints of $\rho$.
Let $u_i$ and $v_i$ be the initial and the terminal vertices of $\rho_i:= f^i(\rho)$.
Since $f$~acts on the finite set $\Gamma^0\times \Gamma^0$, there exist natural numbers $r,n$ such that  $(u_i,v_i)=(u_{i+n},v_{i+n})$ for $i\geqslant r$.

First we check, whether $f^k(\rho)=\tau$ for $k<r$. If yes, we are done, if no we investigate the case $k\geqslant r$.
Given such $k$, we can write $k=i+\ell n$ for some $\ell\geqslant 0$ and $r\leqslant i< r+n$.
So, we have $f^k(\rho)=g^{\ell}(\rho_i)$, where $g:=f^n$. Thus we have to investigate $n$ problems: does there exist
$\ell\geqslant 0$ such that $g^{\ell}(\rho_i)=\tau$, $r\leqslant i< r+n$. Note that $g$ fixes the endpoints of $\rho_i$.

So, from the beginning, we may assume that $f$ fixes the endpoints of $\rho$ and
$\alpha(\rho)=\alpha(\tau)$, and $\omega(\rho)=\omega(\tau)$.

Let $\Gamma_1$ be the graph obtained from $\Gamma$ by adding
a new vertex $v$ and two new oriented edges: $E_1$ from $v$ to $\alpha(\rho)$ and $E_2$ from $v$ to $\omega(\rho)$.
Let $f_1:\Gamma_1 \rightarrow \Gamma_1$ be the extension of $f$ mapping $E_1$ to $E_1$ and $E_2$ to $E_2$.
Clearly, $f_1$ is a homotopy equivalence which fixes $v$. Let $(f_1)_v:\pi_1(\Gamma_1,v)\rightarrow \pi_1(\Gamma_1,v)$
be the induced automorphism. We have
$$f^k(\rho)=\tau\hspace*{1mm}\Longleftrightarrow\hspace*{1mm} f_1^k(E_1\rho \bar{E}_2)=E_1\tau \bar{E}_2
\hspace*{1mm}\Longleftrightarrow\hspace*{1mm} (f_1)_v^k([[E_1\rho \bar{E}_2]])=[[E_1\tau \bar{E}_2]].$$
Thus, the problem is solvable by Theorem~\ref{thm 3.7}.\hfill $\Box$

\medskip

\begin{cor}~\label{cor 3.9} There exists an efficient algorithm which, given a finite connected graph $\Gamma$
and a homotopy equivalence $f:\Gamma\rightarrow \Gamma$ with $f(\Gamma^0)\subseteq \Gamma^0$, and given
two edge paths $\rho, \tau$ in $\Gamma$, decides whether there exist natural
numbers $k>s$ such that $f^k(\rho)=f^s(\tau)$. If such $k$ and $s$ exist, then the algorithm  computes them.
\end{cor}


{\it Proof.}
Let $u_i$ and $v_i$ be the initial and the terminal vertices of $f^i(\rho)$, $i\geqslant 0$, and
let $u'_j$ and $v'_j$ be the initial and the terminal vertices of $f^j(\tau)$, $j\geqslant 0$.
First we can decide, whether there exist $i,j$ such that $(u_i,v_i)=(u_j',v_j')$.
If such $i,j$ don't exist, then the desired $k,s$ don't exist. If such $i,j$ exist, we can algorithmically find natural $i,j,n$ with the following properties:

1) $(u_i,v_i)=(u_j',v_j')$;

2) $(u_i,v_i)=(u_{i+n},v_{i+n})$ and $n$ is minimal;

3) $i>j$;

4) $i-j$ is the minimal possible for 1)-3).

\medskip

So, we reduce the problem to the following: does there exist $p\geqslant q\geqslant 0$ such that $f^{i+pn}(\rho)=f^{j+qn}(\tau)$?
We set $\rho_1:=f^i(\rho)$,  $\tau_1:=f^j(\tau)$, $g:=f^n$. Then the endpoints of $\rho_1$ and $\tau_1$ coincide and are fixed by $g$. In this setting we have to decide, whether there exist $p\geqslant q\geqslant 0$ such that $g^p(\rho_1)=g^q(\tau_1)$.

We extend $\Gamma$ to $\Gamma'$ by adding an edge $E$ from $v_i$ to $u_i$ and we extend $g$ to $g':\Gamma'\rightarrow \Gamma'$ by setting $g'|_{\Gamma} =g$ and $g(E)=E$.
Then the problem is equivalent to the following:
does there exist $p\geqslant q\geqslant 0$ such that $g'^p(\rho_1E)=g'^q(\tau_1E)$?

Since $g'$ is a homotopy equivalence and $\rho_1E$ and $\tau_1E$ are loops based at the same point, and this point is fixed by $g'$, we have $$g'^p(\rho_1E)=g'^q(\tau_1E) \Longleftrightarrow g'^{p-q}(\rho_1E)=\tau_1E.$$ Thus, the problem can be reformulated as follows:
does there exist $m\geqslant 0$ such that $g'^m(\rho_1E)=\tau_1E$? This can be decided by Theorem~\ref{thm 3.7}.\hfill $\Box$

\medskip

We need the following bounded cancelation lemma from~\cite{Cooper}, where it is credited to Grayson and Thurston.


\begin{lem}\label{lem 3.12} Let $\Gamma$ be a finite connected graph and $f:\Gamma \rightarrow \Gamma$ be a homotopy equivalence sending edges to edge paths. Let $\tau_1,\tau_2$ be reduced paths in $\Gamma$ such that $\omega(\tau_1)=\alpha(\tau_2)$ and the path $\tau_1\tau_2$ is reduced. Then
$$l([f(\tau_1\tau_2)])\geqslant l([f(\tau_1)])+l([f(\tau_2)])-2C_{\star},$$
where $C_{\star}>0$ is an algorithmically computable constant which depends only on $f$.
\end{lem}


\section{Graphs $D_f$ and $CoRe(C_f)$ for a homotopy equivalence $f:\Gamma\rightarrow \Gamma$}

Let $\Gamma$ be a finite connected graph with a distinguished vertex $v_{\ast}$. Let $f:\Gamma\rightarrow \Gamma$ be a homotopy equivalence
which maps vertices of $\Gamma$ to vertices and edges to reduced edge paths,
and suppose that $f$ fixes $v_{\ast}$. We consider the group
$$\overline{{\text {\rm Fix}}}(f):=\{[[p]]\in \pi_1(\Gamma,v_{\ast})\,|\, f(p)=p\}.$$

In papers~\cite{GT1,Turner}, the authors suggest a procedure  for computation of a basis of $\overline{{\text {\rm Fix}}}(f)$ with the help of
a graph $D_f$.
This procedure is not an algorithm in general case, since one cannot determine from the beginning,
whether it terminates or not. We give a description of this procedure. We also show that the procedure can be converted into an algorithm if the Membership and the Finiteness problems can be algorithmically solved.

First, we recall some constructions and facts from~\cite{GT1,Turner} and ~\cite{CL}.

\medskip

{\bf A. Definition of $f$-paths.} An edge path $\mu$ in $\Gamma$ is called an {\it $f$-path} if the last point of $\mu$
coincides with the first point of $f(\mu)$.

\vspace*{-15mm}
\includegraphics[scale=0.4]{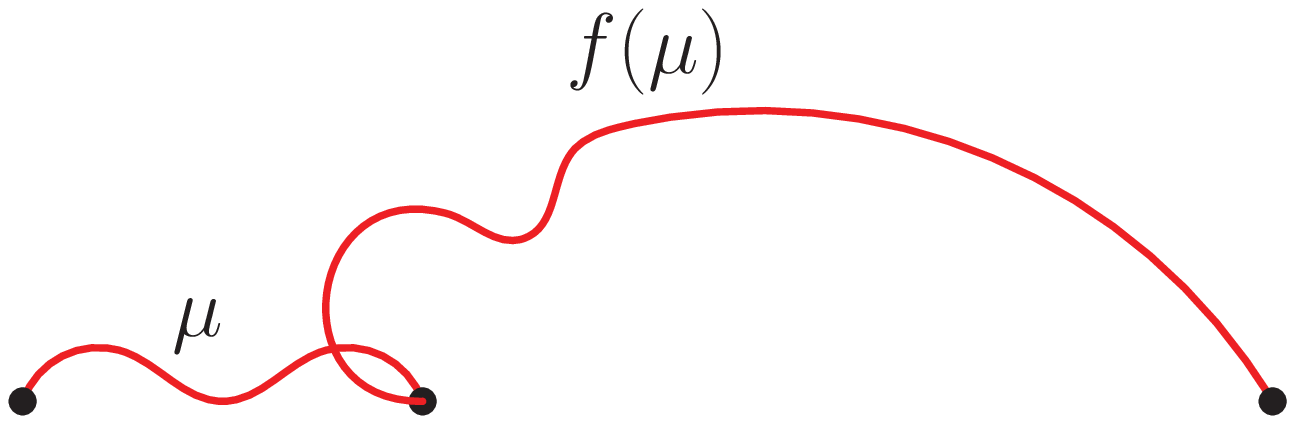}

\vspace*{-110mm}\hspace*{70mm}
\includegraphics[scale=0.4]{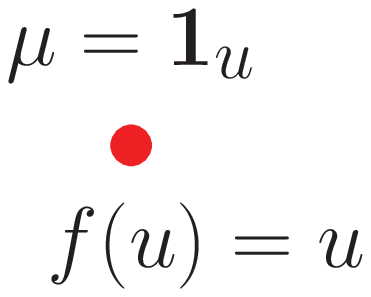}

\vspace*{-85mm}
\begin{center}
Figure 2.

Nontrivial and trivial $f$-paths $\mu$.
\end{center}

\bigskip

We note the following properties of $f$-paths:


\begin{enumerate}
\item[-] the trivial path at a vertex $u$ of $\Gamma$, denoted $\mathbf{1}_u$, is an $f$-path if and only if $u$ is fixed by $f$;

\item[-] if $\mu$ is an $f$-path, then $[\mu]$ is also an $f$-path;

\item[-] if $\mu$ is an $f$-path and $E$ is an edge in $\Gamma$ such that $\alpha(E)=\alpha(\mu)$, then $\overline{E}\mu f(E)$ is also an $f$-path.
\end{enumerate}
\vspace*{-1mm}

{\bf B. Definition of the graph $D_f$.} The vertices of $D_f$ are {\sl reduced} $f$-paths in $\Gamma$.
Let $\mu$ be a reduced $f$-path in $\Gamma$
and let $E_1,\dots ,E_n$ be all edges in $\Gamma$ outgoing from $\alpha(\mu)$.
Then we connect the vertex $\mu$ of $D_f$ to the vertices
$[\overline{E_1}\mu f(E_1)],\dots, [\overline{E_n}\mu f(E_n)]$ by edges with labels $E_1,\dots ,E_n$, respectively, see Figure~3.
The label of a nontrivial edge path in the graph $D_f$ is the product of labels of consecutive edges of this path.
The label of a trivial edge path at a vertex $\mu$ of $D_f$ is $\bold{1}_{\alpha(\mu)}$.

In general, the graph $D_f$ can have infinitely many connected components and some of them can be infinite.
 For a vertex $\mu$ of $D_f$, let $D_f(\mu)$ be the component of $D_f$ containing $\mu$. Lemma~\ref{lem 4.2} says that $\pi_1(D(\mathbf{1}_{v_{\ast}}),\mathbf{1}_{v_{\ast}})\cong \overline{{\text {\rm Fix}}}(f)$.
This was first proved in~\cite{GT1} with the help of preferable directions at vertices of $D_f$.


\medskip

{\bf C. Preferable directions at vertices of $D_f$, dead and alive vertices of~$D_f$.}
For a reduced nontrivial $f$-path $\mu$ in $\Gamma$, we set $\widehat{f}(\mu):=[\overline{E}\mu f(E)]$,
where $E$ is the first edge of $\mu$.
Then $\mu$ and $\widehat{f}(\mu)$ are vertices of the graph $D_f$ connected by the edge with the label $E$. The direction of this edge is called {\it preferable} at the vertex~$\mu$. We will put the symbol $\vartriangleright$ on this edge
near the vertex~$\mu$.

\vspace*{5mm}
\hspace*{20mm} {\bf in $\Gamma$:\hspace*{70mm} in $D_f$:}

\vspace*{-8mm}
\hspace*{-10mm}
\includegraphics[scale=0.4]{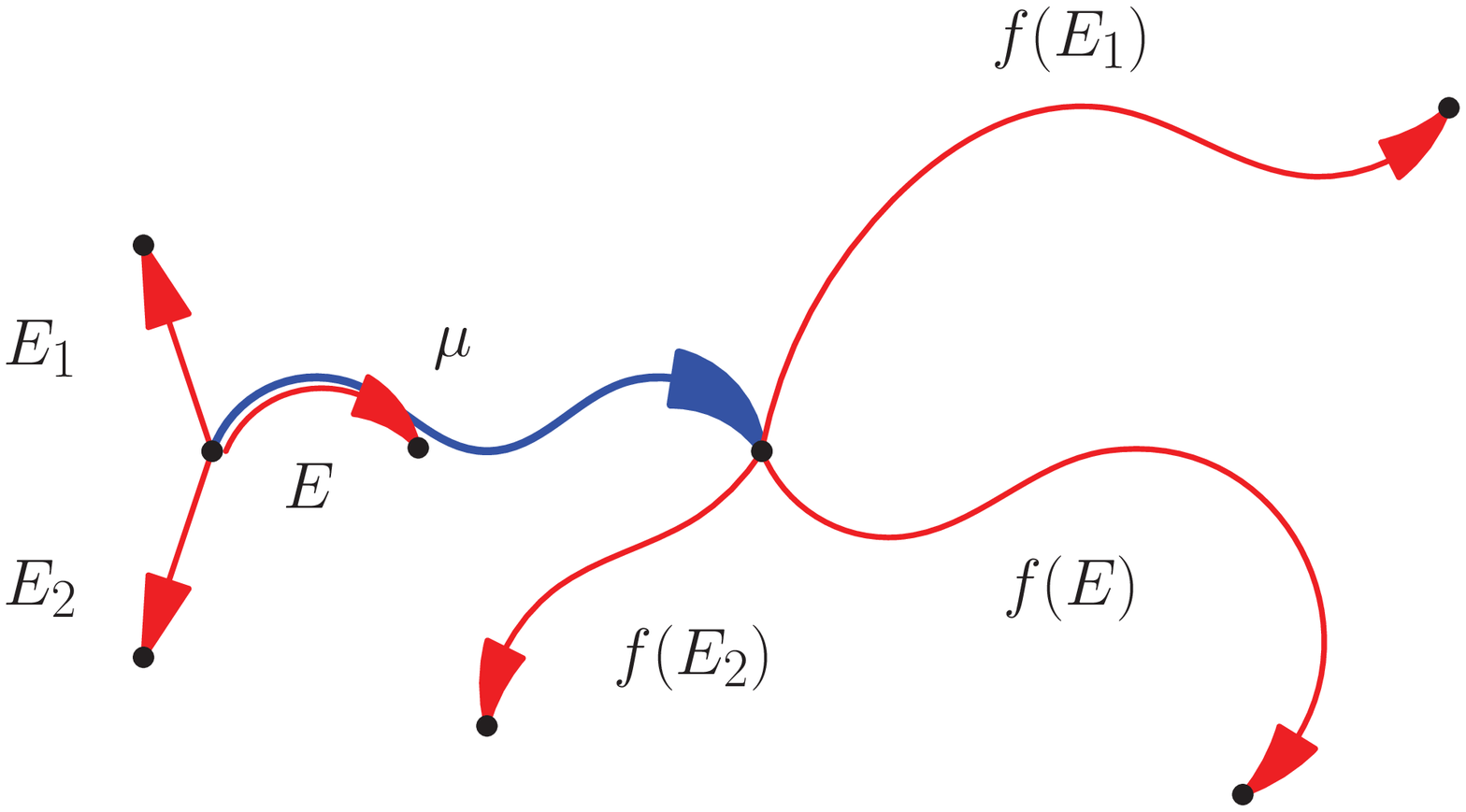}

\vspace*{-125mm}\hspace*{80mm}
\includegraphics[scale=0.466]{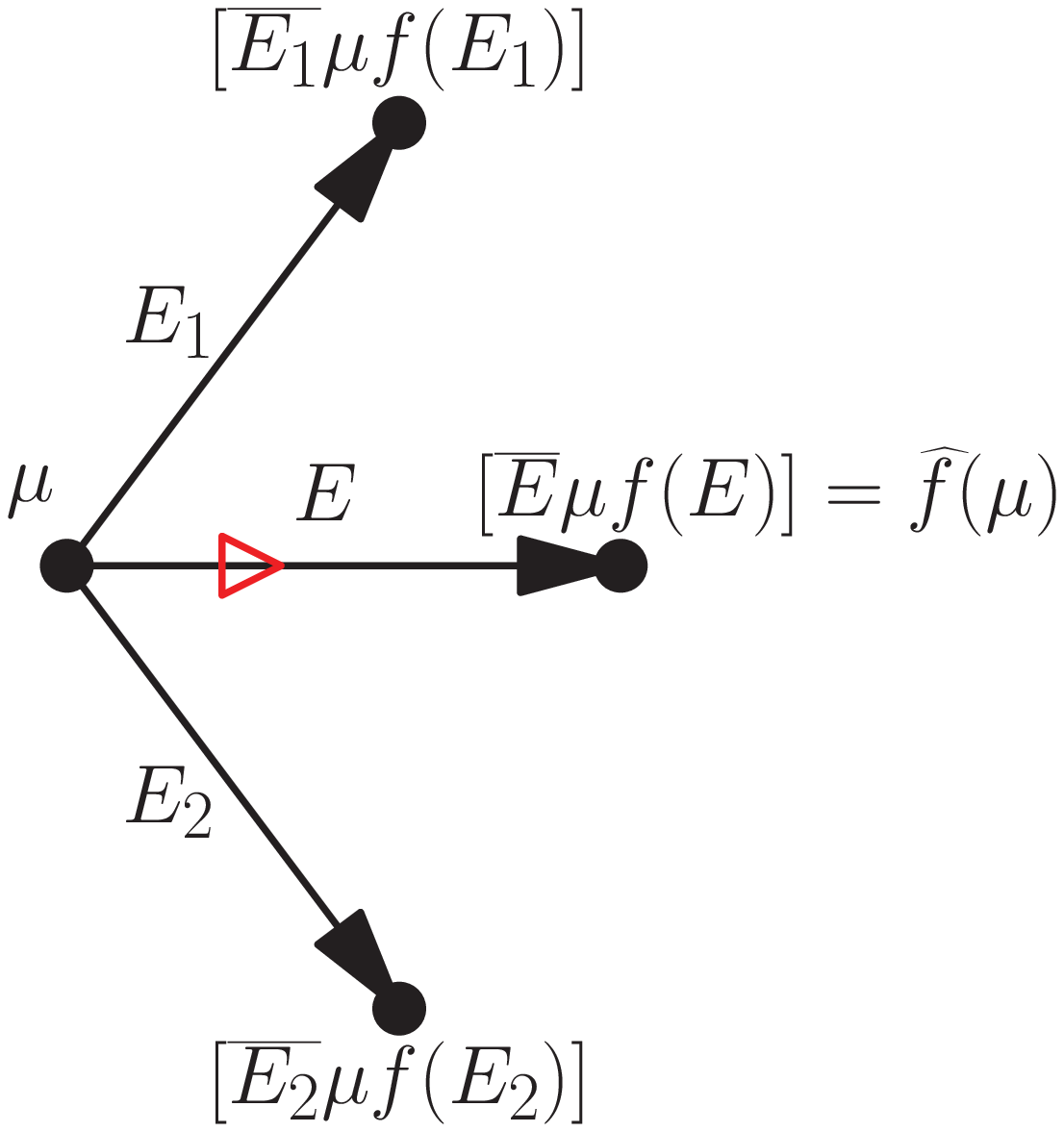}

\vspace*{-58mm}

\begin{center}
Figure 3.


\begin{minipage}{11cm}
From the graph $\Gamma$ to the graph $D_f$. On the left we consider $\mu$ as a path in $\Gamma$,
on the right as a vertex in $D_f$. The red triangle on the right shows the preferable direction at the vertex $\mu$.
\end{minipage}
\end{center}

\bigskip

Note that only the vertices $\mathbf{1}_{w}$, where $w\in \Gamma^0$ and $f(w)=w$, do not
admit a preferable direction. We call such vertices {\it dead} and all other vertices of $D_f$  {\it alive}.
Observe that at each vertex of $D_f$, there is at most one outwardly $\vartriangleright$-directed edge.

\medskip

{\bf D. Ordinary, repelling and attracting edges of $D_f$.}

The following definition is illustrated by Figure~4.

\begin{defn}\label{ora}

{\rm
Let $e$ be an edge of $D_f$, let $p,q$ be the initial and the terminal vertices of $e$, and let $E\in \Gamma^1$ be the label of $e$.

\begin{itemize}
\item[(1)] The edge $e$ is called {\it ordinary} in $D_f$ if one of the following holds:
\begin{itemize}
\item[(a)] $E$ is the first edge of the path $p$ in $\Gamma$ and $\overline{E}$ is not the first edge of the path $q$ in $\Gamma$.

\item[(b)] $E$ is not the first edge of the path $p$ in $\Gamma$ and $\overline{E}$ is the first edge of the path $q$ in $\Gamma$.
\end{itemize}

\medskip

\item[(2)] The edge $e$ is called {\it repelling} in $D_f$ if $E$ is not the first edge of the path $p$ in $\Gamma$ and $\overline{E}$ is not the first edge of the path $q$ in $\Gamma$.

\noindent A vertex of $D_f$ is called {\it repelling} if it is the initial or the terminal vertex of a repelling edge.

\medskip

\item[(3)] The edge $e$ is called {\it attracting} in $D_f$ if $E$ is the first edge of the path $p$ in $\Gamma$ and $\overline{E}$ is the first edge of the path $q$ in $\Gamma$.
\end{itemize}

An edge of $D_f$ is called {\it exceptional} if it is attracting or repelling.
}
\end{defn}

\medskip

\vspace*{-10mm}
\hspace*{10mm}
{\unitlength 1mm
\linethickness{0.4pt}
\ifx\plotpoint\undefined\newsavebox{\plotpoint}\fi 
\begin{picture}(74.75,25)(9.5,64)
\put(65,75){\makebox(0,0)[cc]{\it Ordinary edges:}}
\put(40,65){\vector(1,0){0.5}}\put(27,65){\line(1,0){26}}
\put(27,65){\circle*{2}}
\put(53,65){\circle*{2}}
\put(40,68){\makebox(0,0)[cc]{$E$}}
\put(28,60){\makebox(0,0)[cc]{$p$}}
\put(53,60){\makebox(0,0)[cc]{$q$}}
\put(30.5,65){\makebox(0,0)[cc]{$\vartriangleright$}}

\put(90,65){\vector(1,0){0.5}}\put(77,65){\line(1,0){26}}
\put(77,65){\circle*{2}}
\put(103,65){\circle*{2}}
\put(90,68){\makebox(0,0)[cc]{$E$}}
\put(78,60){\makebox(0,0)[cc]{$p$}}
\put(103,60){\makebox(0,0)[cc]{$q$}}
\put(100,65){\makebox(0,0)[cc]{$\vartriangleleft$}}

%
\put(65,45){\makebox(0,0)[cc]{\it Exceptional edges:}}

\put(40,35){\vector(1,0){0.5}}\put(27,35){\line(1,0){26}}
\put(27,35){\circle*{2}}
\put(53,35){\circle*{2}}
\put(40,38){\makebox(0,0)[cc]{$E$}}
\put(28,30){\makebox(0,0)[cc]{$p$}}
\put(53,30){\makebox(0,0)[cc]{$q$}}
\put(30.5,35){\makebox(0,0)[cc]{$\vartriangleright$}}
\put(50,35){\makebox(0,0)[cc]{$\vartriangleleft$}}
\put(40,22){\makebox(0,0)[cc]{\it Attracting edge}}

\put(90,35){\vector(1,0){0.5}}\put(77,35){\line(1,0){26}}
\put(77,35){\circle*{2}}
\put(103,35){\circle*{2}}
\put(90,38){\makebox(0,0)[cc]{$E$}}
\put(78,30){\makebox(0,0)[cc]{$p$}}
\put(103,30){\makebox(0,0)[cc]{$q$}}
\put(90,22){\makebox(0,0)[cc]{\it Repelling edge}}

\end{picture}
}

\vspace*{50mm}

\begin{center}
Figure 4. Different types of edges in $D_f$.
\end{center}

\medskip

\begin{prop}\label{prop 4.1} {\rm (see~\cite{GT1, Turner} and~\cite{CL})} {\rm (a)} The repelling edges of $D_f$ are in 1-1 correspondence with the occurrences of edges $E$ in $f(E)$, where $E\in \Gamma^1$.
More precisely, there exists a bijection of the type:

$$\hspace*{-10mm}f(E)\equiv \overline{u}\cdot E\cdot v \hspace*{3mm}\Longleftrightarrow
\begin{cases}
\vspace*{-10mm}
\hspace*{0mm}
{\unitlength 1mm
\linethickness{0.4pt}
\ifx\plotpoint\undefined\newsavebox{\plotpoint}\fi 
\begin{picture}(74.75,65)(20,5)
\put(40,65){\vector(1,0){0.5}}\put(20,65){\line(1,0){40}}
\put(27,65){\circle*{2}}
\put(53,65){\circle*{2}}
\put(40,68){\makebox(0,0)[cc]{$E$}}
\put(27,61.5){\makebox(0,0)[cc]{$u$}}
{
\linethickness{0.2pt}
\put(27,65){\line(-1,1){4}}
\put(27,65){\line(-1,-1){4}}
}
\put(53,61.5){\makebox(0,0)[cc]{$v$}}
\put(56.5,65){\makebox(0,0)[cc]{$\vartriangleright$}}
\put(23.5,65){\makebox(0,0)[cc]{$\vartriangleleft$}}
\put(93.5,65){\makebox(0,0)[cc]{\rm if $u$ and $v$ are nonempty,}}
\hspace*{0mm}
\put(40,50){\vector(1,0){0.5}}\put(20,50){\line(1,0){40}}
\put(27,50){\circle*{2}}
\put(53,50){\circle*{2}}
\put(40,53){\makebox(0,0)[cc]{$E$}}
\put(27,46.5){\makebox(0,0)[cc]{$\mathbf{1}_{\alpha(E)}$}}
{
\linethickness{0.2pt}
\put(27,50){\line(-1,1){4}}
\put(27,50){\line(-1,-1){4}}
}
\put(53,46.5){\makebox(0,0)[cc]{$v$}}
\put(56.5,50){\makebox(0,0)[cc]{$\vartriangleright$}}
\put(93.5,50){\makebox(0,0)[cc]{\rm if $u$ is empty and $v$ not,}}
\hspace*{0mm}
\put(40,35){\vector(1,0){0.5}}\put(20,35){\line(1,0){40}}
\put(27,35){\circle*{2}}
\put(53,35){\circle*{2}}
\put(40,38){\makebox(0,0)[cc]{$E$}}
\put(27,31.5){\makebox(0,0)[cc]{$u$}}
{
\linethickness{0.2pt}
\put(27,35){\line(-1,1){4}}
\put(27,35){\line(-1,-1){4}}
}
\put(53,31.5){\makebox(0,0)[cc]{$\mathbf{1}_{\omega(E)}$}}
\put(23.5,35){\makebox(0,0)[cc]{$\vartriangleleft$}}
\put(93.5,35){\makebox(0,0)[cc]{\rm if $v$ is empty and $u$ not,}}
\hspace*{0mm}
\put(40,20){\vector(1,0){0.5}}\put(20,20){\line(1,0){40}}
\put(27,20){\circle*{2}}
\put(53,20){\circle*{2}}
\put(40,23){\makebox(0,0)[cc]{$E$}}
\put(27,16.5){\makebox(0,0)[cc]{$\mathbf{1}_{\alpha(E)}$}}
{
\linethickness{0.2pt}
\put(27,20){\line(-1,1){4}}
\put(27,20){\line(-1,-1){4}}
}
\put(53,16.5){\makebox(0,0)[cc]{$\mathbf{1}_{\omega(E)}$}}
\put(93.5,20){\makebox(0,0)[cc]{\rm if $u$ and $v$ are empty.}}
\end{picture}
} & \\
\end{cases}
$$

\bigskip

{\rm (b)}
The attracting edges of $D_f$ are in 1-1 correspondence with the occurrences of edges $\overline{E}$ in $f(E)$, where $E\in \Gamma^1$.
More precisely, there exists a bijection of the type:

\vspace*{-10mm}
$$\hspace*{-10mm}f(E)\equiv \overline{u}\cdot \overline{E}\cdot v \hspace*{3mm}\Longleftrightarrow
\vspace*{0mm}
\hspace*{10mm}
{\unitlength 1mm
\linethickness{0.4pt}
\ifx\plotpoint\undefined\newsavebox{\plotpoint}\fi 
\begin{picture}(74.75,15)(26.5,64)
\put(40,65){\vector(1,0){0.5}}\put(20,65){\line(1,0){40}}
\put(27,65){\circle*{2}}
\put(53,65){\circle*{2}}
\put(40,68){\makebox(0,0)[cc]{$E$}}
\put(28,60){\makebox(0,0)[cc]{$E\cdot u$}}
{
\linethickness{0.2pt}
\put(27,65){\line(-1,1){4}}
\put(27,65){\line(-1,-1){4}}
}
\put(53,60){\makebox(0,0)[cc]{$\overline{E}\cdot v$}}
\put(50,65){\makebox(0,0)[cc]{$\vartriangleleft$}}
\put(30.5,65){\makebox(0,0)[cc]{$\vartriangleright$}}
\end{picture}
}
$$

\vspace*{10mm}
{\rm (c)} There exist only finitely many exceptional edges in $D_f$ and they can be algorithmically found.

\end{prop}

\medskip

The following fundamental lemma was first proved by Goldstein and Turner in~\cite{GT1}. We reproduce here their nice proof for completeness.

\begin{lem}\label{lem 4.2} {\rm (see~\cite{GT1})} The fundamental group of each component of $D_f$ is finitely generated.
Moreover,
$\pi_1(D_f(\mathbf{1}_{v_{\ast}}),\mathbf{1}_{v_{\ast}})\cong \overline{{\text {\rm Fix}}}(f)$.
\end{lem}

{\it Proof.}
Observe that a connected locally finite graph has finite rank if and only if the edges of this graph can be directed so that at all but a finite number of vertices, there is at most one outwardly directed edge.

Let $D_f'$ be the graph obtained from $D_f$ by removing all exceptional edges. Each component of $D_f'$ has only ordinary edges, and these edges carry preferable directions which satisfy the above observation. Hence each component of $D_f'$ has finite rank. Since there is only finitely many exceptional edges in $D_f$, each component of $D_f$ has finite rank.


Let us prove the second claim in a more general context.
Let $\mu$ be a vertex in $D_f$ and suppose that $c$ is a closed path in $D_f$ based at $\mu$. Let $\ell=E_1E_2\dots E_k$ be the label of $c$; so all $E_i$ are edges in $\Gamma$. Then $$\mu,\hspace*{2mm} [\overline{E}_1\mu f(E_1)],\hspace*{2mm} [\overline{E}_2\overline{E}_1\mu f(E_1)f(E_2)],\hspace*{2mm} \dots, \hspace*{2mm}[\overline{\ell}\mu f(\ell)]$$ are consecutive vertices of $c$ and we have
$[\bar{\ell}\mu f(\ell)]=\mu$.

So, $\ell$ is a closed edge path in $\Gamma$ based at $\alpha(\mu)$ and satisfying $[\mu f(\ell)\overline{\mu}]=[\ell]$.
The correspondence $c\mapsto \ell$ induces the isomorphism
$$\pi_1(D_f(\mu),\mu)\cong \{[[p]]\in \pi_1(\Gamma,\alpha(\mu))\,|\, [[\mu f(p)\bar{\mu}]]=[[p]])\}.$$


Setting $\mu:=\mathbf{1}_{v_{\ast}}$, we obtain the second claim of the lemma.  \hfill $\Box$

\bigskip

{\bf E. Definition of a $\mu$-subgraph of $D_f$.}
Let $\mu$ be a vertex in $D_f$. If $\mu$ is not
a dead vertex, i.e. if $\mu\equiv E_1E_2\dots E_m$ for some edges $E_i\in \Gamma^1$, $m\geqslant 1$, we can pass from $\mu$ to the vertex $\widehat{f}(\mu)\equiv [E_2\dots E_mf(E_1)]$
by using the direction which is preferable at $\mu$.

The vertices of the $\mu$-{\it subgraph} are the vertices $\mu_1,\mu_2,\dots $ of $D_f$
such that $\mu_1=\mu$ and $\mu_{i+1}=\widehat{f}(\mu_i)$ if the vertex $\mu_i$ is not dead, $i\geqslant 1$.
The edges of the $\mu$-subgraph are those which connect $\mu_i$ with $\mu_{i+1}$
and carry the preferable direction at $\mu_i$.




\bigskip

Note that the $\mu$-subgraph is finite if and only if starting from $\mu$ and moving along the preferable directions we will
came to a dead vertex or to a vertex which we have seen earlier. If the $\mu$-subgraph is infinite, we call it a {\it $\mu$-ray}.
Thus, any $\mu$-subgraph is one of the following four types:

\vspace*{25mm}
\hspace*{5mm}
{
\unitlength=1.00mm
\special{em:linewidth 0.4pt}
\linethickness{0.4pt}
\begin{picture}(84.00,41.00)
\put(5.00,55.00){\line(1,0){20.00}}
\put(30.50,55.00){\makebox(0,0)[cc]{$\cdots$}}
\put(35.00,55.00){\line(1,0){20.00}}
\put(5.00,55.00){\circle*{1.00}}
\put(15.00,55.00){\circle*{1.00}}
\put(25.00,55.00){\circle*{1.00}}
\put(35.00,55.00){\circle*{1.00}}
\put(45.00,55.00){\circle*{1.00}}
\put(55.00,55.00){\circle*{2.00}}
\put(62.00,55.00){\makebox(0,0)[lc]{\small a segment ending at a dead vertex}}
\put(5.00,53.00){\makebox(0,0)[ct]{$\mu$}}

\put(7.00,54.95){\makebox(0,0)[cc]{$\triangleright$}}
\put(17.00,54.95){\makebox(0,0)[cc]{$\triangleright$}}
\put(37.00,54.95){\makebox(0,0)[cc]{$\triangleright$}}
\put(47.00,54.95){\makebox(0,0)[cc]{$\triangleright$}}
\put(5.00,40.00){\line(1,0){20.00}}
\put(30.50,40.00){\makebox(0,0)[cc]{$\cdots$}}
\put(35.00,40.00){\line(1,0){20.00}}
\put(5.00,40.00){\circle*{1.00}}
\put(15.00,40.00){\circle*{1.00}}
\put(25.00,40.00){\circle*{1.00}}
\put(35.00,40.00){\circle*{1.00}}
\put(45.00,40.00){\circle*{1.00}}
\put(55.00,40.00){\circle*{1.00}}
\put(62.00,42.00){\makebox(0,0)[lc]{\small a segment with an attracting edge}}
\put(62.00,38.00){\makebox(0,0)[lc]{\small which can be considered as a cycle}}
\put(5.00,38.00){\makebox(0,0)[ct]{$\mu$}}

\put(7.00,39.95){\makebox(0,0)[cc]{$\triangleright$}}
\put(17.00,39.95){\makebox(0,0)[cc]{$\triangleright$}}
\put(37.00,39.95){\makebox(0,0)[cc]{$\triangleright$}}
\put(47.00,39.95){\makebox(0,0)[cc]{$\triangleright$}}
\put(52.8,39.95){\makebox(0,0)[cc]{$\triangleleft$}}
\put(5.00,25.00){\line(1,0){20.00}}
\put(30.50,25.00){\makebox(0,0)[cc]{$\cdots$}}
\put(35.00,25.00){\line(1,0){20.00}}
\put(5.00,25.00){\circle*{1.00}}
\put(15.00,25.00){\circle*{1.00}}
\put(25.00,25.00){\circle*{1.00}}
\put(35.00,25.00){\circle*{1.00}}
\put(45.00,25.00){\circle*{1.00}}
\put(55.00,25.00){\circle*{1.00}}
\put(5.00,23.00){\makebox(0,0)[ct]{$\mu$}}

\put(7.00,24.95){\makebox(0,0)[cc]{$\triangleright$}}
\put(17.00,24.95){\makebox(0,0)[cc]{$\triangleright$}}
\put(37.00,24.95){\makebox(0,0)[cc]{$\triangleright$}}
\put(47.00,24.95){\makebox(0,0)[cc]{$\triangleright$}}
\put(62.00,25.00){\circle{14.00}}
\put(58.00,31.00){\circle*{1.00}}
\put(58.00,19.00){\circle*{1.00}}
\put(70.00,26.00){\makebox(0,0)[lc]{$\vdots$}}
\put(79.00,25.00){\makebox(0,0)[lc]{\small a segment with a cycle}}
\put(65.00,31.00){\circle*{1.00}}
\put(55.50,22.00){\makebox(0,0)[cc]{\rotatebox[origin=c]{-60}{$\triangleleft$}}}
\put(55.75,28.00){\makebox(0,0)[cc]{\rotatebox[origin=c]{60}{$\triangleright$}}}
\put(62.00,31.95){\makebox(0,0)[cc]{$\triangleright$}}
\put(5.00,10.00){\line(1,0){30.00}}
\put(5.00,10.00){\circle*{1.00}}
\put(15.00,10.00){\circle*{1.00}}
\put(25.00,10.00){\circle*{1.00}}
\put(35.00,10.00){\circle*{1.00}}
\put(40.00,10.00){\makebox(0,0)[cc]{$\cdots$}}
\put(5.00,8.00){\makebox(0,0)[ct]{$\mu$}}
\put(49.00,10.00){\makebox(0,0)[lc]{\small a ray}}
\put(7.00,9.95){\makebox(0,0)[cc]{$\triangleright$}}
\put(17.00,9.95){\makebox(0,0)[cc]{$\triangleright$}}
\put(27.00,9.95){\makebox(0,0)[cc]{$\triangleright$}}
\end{picture}
}

\begin{center} Figure 5.

\medskip

The types of $\mu$-subgraphs in $D_f$.

\end{center}

\bigskip

Let $\mu$ and $\tau$ be two vertices of $D_f$. Clearly, if the $\mu$-subgraph and the $\tau$-subgraph intersect,
then they differ only by their finite ``initial subsegments''. The vertex $\mu_i$ from this intersection with
minimal $i$ is called the {\it intersection point} of these subgraphs.

Note the following properties of $\mu$-subgraphs:

\begin{enumerate}
\item[-] if $\mu_0$ is a vertex of a $\mu$-subgraph, then the $\mu_0$-subgraph is contained in the $\mu$-subgraph;

\item[-] if $\mu_0$ is a vertex of a $\mu$-subgraph and $\tau_0$ is a vertex of a $\tau$-subgraph, then the $\mu$-subgraph and the $\tau$-subgraph intersect if and only if  the $\mu_0$-subgraph and the $\tau_0$-subgraph intersect;

\item[-] a $\mu$-ray does not intersect a finite $\tau$-subgraph.
\end{enumerate}

\medskip

From this point, we start to develop the above approach.

\bigskip

{\bf F. Definitions of the graphs $C_f$ and $CoRe(C_f)$.} A component of $D_f$ is called {\it repelling} if it contains at least one repelling edge.
Let $C_1,\dots ,C_n$ be all repelling components of $D_f$.
For each $C_i$, let $CoRe(C_i)$ be the minimal connected subgraph of $C_i$ which contains all repelling edges of $C_i$
and carries $\pi_1(C_i)$.
We set
$C_f :=\underset{i=1}{\overset{n}{\cup}} C_i$ and $CoRe(C_f) :=\underset{i=1}{\overset{n}{\cup}} CoRe(C_i)$.

\bigskip

Below we show how to compute a basis of the group $\pi_1(D_f(\mathbf{1}_{v_{\ast}}), \mathbf{1}_{v_{\ast}})$
if we know how to construct the graph $CoRe(C_f)$.

\bigskip

\begin{lem}\label{lem 4.3} Let $\mathbf{1}_u$ be a dead vertex of $D_f$. If the component $D_f(\mathbf{1}_u)$ is non-contractible,
then it lies in $C_f$.
\end{lem}




\medskip

{\it Proof.}
Suppose that $D_f(\mathbf{1}_u)$ is non-contractible.
Then there exists an edge path $p=E_1E_2\dots E_k$ in $D_f(\mathbf{1}_u)$ such that $\omega(E_k)=\mathbf{1}_u$,
the edges of $p$ are distinct, and $\alpha(E_1)=\omega(E_s)$ for some $1\leqslant s\leqslant k$.
We show that $p$ contains a repelling edge.
Suppose not, then the direction of $E_k$ is preferable at $\alpha(E_k)$. By induction, the direction of $E_i$ is
preferable at the point $\alpha(E_i)$ for every $i=1,\dots,k$ (see Figure 6).
In particular, $\alpha(E_1)\neq \mathbf{1}_u$, and hence $k>s$.
Then there are two preferable directions at $\alpha(E_1)$, namely the direction of $E_1$
and the direction of $E_{s+1}$, a contradiction. Thus, $p$ must contain a repelling edge. \hfill $\Box$



\vspace*{0mm}\hspace*{30mm}\includegraphics[scale=0.5]{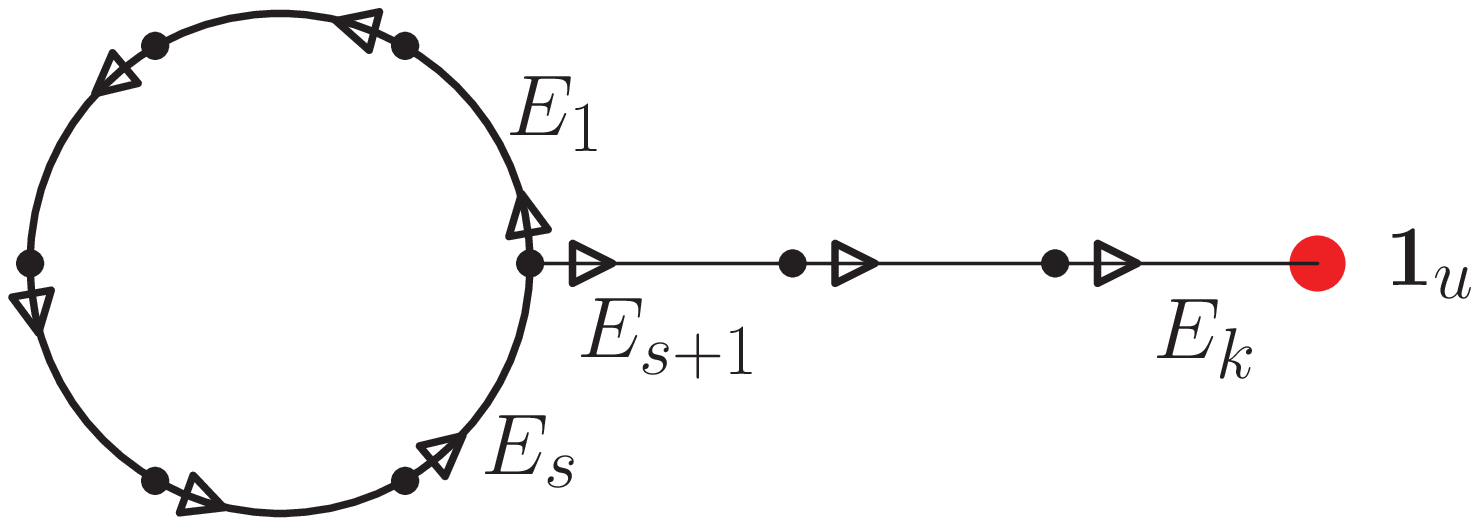}

\vspace*{-110mm}
\begin{center}
Figure 6.
\end{center}




\bigskip


\begin{lem}\label{lem 4.4} Let $\mathbf{1}_u$ be a dead vertex of $D_f$. The vertex $\mathbf{1}_u$ lies in $C_f$ if and only if it lies in the $\mu$-subgraph for some repelling vertex $\mu$.
\end{lem}


{\it Proof.} Suppose that $\mathbf{1}_u$ lies in $C_f$. Then there exists a shortest path $E_1E_2\dots E_k$, where the edge $E_1$ is repelling and $\omega(E_k)=\mathbf{1}_u$. If $k=1$, then $\mathbf{1}_u$ is repelling as a vertex of a repelling edge, and we are done. If $k\geqslant 2$,
the edges $E_2,\dots ,E_k$ are not repelling. In particular, the preferable direction at $\alpha(E_k)$ coincides with the direction of $E_k$. By induction one can prove that the preferable direction at $\alpha(E_i)$ coincides with the direction of $E_i$ for $i\geqslant 2$.
Then $\mathbf{1}_u$ lies in the $\mu$-subgraph for $\mu:=\omega(E_1)$ and the vertex $\mu$ is repelling.
The converse claim is clear.
\hfill $\Box$

\medskip

\begin{prop}~\label{prop 4.5} Suppose we can algorithmically construct $CoRe(C_f)$
and decide whether there exists a repelling vertex $\mu$ such that the $\mu$-subgraph contains
the vertex~$\mathbf{1}_{v_{\ast}}$, and if such $\mu$ exists, we can compute the $\mu$-subgraph.
Then we can compute
a basis of $\pi_1(D_f(\mathbf{1}_{v_{\ast}}),\mathbf{1}_{v_{\ast}})$.
\end{prop}

\medskip

{\it Proof.} Suppose that $\mathbf{1}_{v_{\ast}}$ does not lie in the $\mu$-subgraph for any repelling vertex~$\mu$. Then, by Lemma~\ref{lem 4.4}, $\mathbf{1}_{v_{\ast}}\not\in C_f$ and, by Lemma~\ref{lem 4.3}, $D_f(\mathbf{1}_{v_{\ast}})$ is contractible.

Now suppose that $\mathbf{1}_{v_{\ast}}$ lies in the $\mu$-subgraph for some repelling vertex $\mu$.
Since $CoRe(C_f)$ is supposed to be constructible and each repelling vertex lies in $CoRe(C_f)$,
we can find the component of $CoRe(C_f)$ containing~$\mu$.
Let $\Delta$ be the union of this component and the $\mu$-subgraph; note that the $\mu$-subgraph terminates at $\mathbf{1}_{v_{\ast}}$.
Then~$\Delta$ is a core of $D_f(\mathbf{1}_{v_{\ast}})$ containing $\mathbf{1}_{v_{\ast}}$. In particular, we can compute a basis of $\pi_1(D_f(\mathbf{1}_{v_{\ast}}),\mathbf{1}_{v_{\ast}})$.
\hfill $\Box$





\bigskip

{\bf G. To construct the graph $CoRe(C_f)$, it suffices to do the following:}

\begin{enumerate}
\item[(1)] Find all repelling edges of $D_f$.

\item[(2)] For each alive repelling vertex $\mu$ determine, whether the $\mu$-subgraph is finite or not.

\item[(3)] Compute all elements of all finite $\mu$-subgraphs from (2).

\item[(4)] For each two repelling vertices $\mu$ and $\tau$ with infinite $\mu$-and $\tau$-subgraphs determine,
whether these subgraphs intersect.

\item[(5)] If the $\mu$-subgraph and the $\tau$-subgraph from (4) intersect,
find their first intersection point and compute their initial segments up to this point.
\end{enumerate}

\medskip







To convert this procedure to an algorithm,
we shall construct algorithms for steps (2) and (4). In papers~\cite{CL} and~\cite{Turner} these algorithms are given
only in some special cases (for positive automorphisms and for irreducible automorphisms represented by train tracks for which each fixed point is a vertex).
The main idea in these papers is to use an {\it inverse preferred direction} at vertices in the graph $D_f$.
This direction can be constructed algorithmically (in general case) with the help of a homotopic inverse to $f$.
It determines its own repelling edges and repelling and dead vertices; they can be algorithmically~found.

\medskip

{\bf H. Inverse preferred directions in $D_f$.} We will realize the following plan. First we
define a map $g:\Gamma\rightarrow \Gamma$ which is a homotopy inverse to $f:\Gamma\rightarrow \Gamma$.
Then we show that there is a label preserving graph map $\Phi:D_f\rightarrow D_g$.
Finally we define the inverse preferred directions at vertices in $D_f$ by pulling back the preferred directions in $D_g$ by $\Phi$. This idea is due to Turner~\cite{Turner}, and has sources in the paper of Cohen and Lustig~\cite{CL}. Note that in~\cite{Turner}, the map $\Phi$ is claimed to be locally injective (see Proposition in Section 3 there), and we claim that $\Phi$ is an isomorphism.

\begin{defn}\label{constant_K}
{\rm For the given homotopy equivalence $f:\Gamma\rightarrow \Gamma$, we can efficiently construct
a homotopy equivalence $g:\Gamma\rightarrow \Gamma$  such that $g$ maps vertices of $\Gamma$ to vertices, edges to egde paths, and the maps $h:=g\circ f$ and $f\circ g$ are homotopic to the identity on~$\Gamma$.
From now on, we fix $g$. Let $H:\Gamma\times [0,1]\rightarrow \Gamma$ be a homotopy from the identity $id$ to $h$.
For each point $u$ in $\Gamma$, let $p_u$ be the path from $u$ to $h(u)$ determined by the homotopy $H$: namely $p_u(t)=H(u,t)$, $t\in [0,1]$. We set $K_{\star}(f):=\max\{l(p_u):\, u\in \Gamma^0 \}$.}
\end{defn}

First we define a map $\Phi$ from the set of vertices of $D_f$ to the set of vertices of~$D_g$.
Let $\mu$ be a vertex in $D_f$. We consider $\mu$ as a reduced $f$-path in $\Gamma$ and let $u$ be the initial vertex of $\mu$. Then we set $\Phi(\mu)=[p_u g(\overline{\mu})]$. Clearly, $\Phi(\mu)$ is a reduced $g$-path in $\Gamma$.
Hence $\Phi(\mu)$ can be considered as a vertex in $D_g$.


\begin{lem}\label{lem 4.6}
The map $\Phi$ can be continued to a graph map $\Phi: D_f\rightarrow D_g$ preserving the labels of edges.
\end{lem}

\medskip

{\it Proof.} Let $\mu$ and $\mu_1$ be two vertices in $D_f$ connected by an edge with label $E$, i.e. $\mu_1=[\overline{E}\mu f(E)]$.
We must show that $\Phi(\mu)$ and $\Phi(\mu_1)$ are connected by an edge with the label $E$, i.e.
$\Phi(\mu_1)=[\overline{E}\Phi(\mu) g(E)]$.
Let $u$ and $w$ be the initial and the terminal vertices of $E$. Then $u$ and $w$ are the initial vertices of $\mu$ and $\mu_1$, respectively. We have
$$\Phi(\mu_1)=[p_w g\bigl(f(\overline{E})\overline{\mu}E\bigr)]=[p_w h(\overline{E})g(\overline{\mu})g(E)]=
[\overline{E}p_u g(\overline{\mu})g(E)]=[\overline{E}\Phi(\mu) g(E)].$$
Here we use the fact that $H$ is a homotopy and hence
$$[h(\ell)]=[\overline{p}_{\alpha(\ell)} \ell p_{\omega(\ell)}]\eqno{(7.1)}$$
for any path $\ell$ in $\Gamma$.
\hfill $\Box$


\begin{rmk}\label{rem 4.6a}
{\rm
Let $\mu$ be a vertex in $D_f$.  Then the following holds:
\begin{enumerate}
\item[1)] The $f$-path $\mu$ and the $g$-path $\Phi(\mu)$ have the same initial vertices in $\Gamma$.
\item[2)] Let $E_1,\dots,E_n$ be the edges outgoing from $\alpha(\mu)$ in $\Gamma$.
Then the vertices $\mu$ and $\Phi(\mu)$ of the graphs $D_f$ and $D_g$ have degree $n$
and the labels of edges outgoing from each of these vertices are $E_1,\dots,E_n$.
\item[3)] $\Phi$ maps the star of the vertex $\mu$ to the star of the vertex $\Phi(\mu)$ bijectively
and label preserving.
\end{enumerate}
}
\end{rmk}

\begin{prop}\label{prop 4.7}
The map $\Phi: D_f\rightarrow D_g$ is an isomorphism of graphs.
\end{prop}

\medskip

{\it Proof.} By Remark~\ref{rem 4.6a}.~3), it suffices to show that $\Phi$ is bijective on vertices.

First we show that $\Phi$ is injective on vertices.
Let $\mu_1,\mu_2$ be two different vertices of $D_f$.
If the $f$-paths $\mu_1$ and $\mu_2$ have different initial vertices in $\Gamma$, then, by Remark~\ref{rem 4.6a}.~1), the
$g$-paths $\Phi(\mu_1)$ and $\Phi(\mu_2)$ have different initial vertices in~$\Gamma$ too, hence $\Phi(\mu_1)\neq \Phi(\mu_2)$.

Suppose that the initial vertices of the $f$-paths $\mu_1$ and $\mu_2$ coincide and equal to~$u$.
Then their terminal vertices also coincide and equal to $f(u)$. Since the $f$-paths $\mu_1,\mu_2$ are reduced, $\mu_1\neq \mu_2$, and $g$ is a homotopy equivalence, we have $[g(\mu_1)]\neq [g(\mu_2)]$, hence $\Phi(\mu_1)=[p_ug(\overline{\mu}_1)]\neq [p_ug(\overline{\mu}_2)]=\Phi(\mu_2)$.

Now we show that $\Phi$ is surjective on vertices. Let $\tau$ be a vertex in $D_g$, i.e. $\tau$ is a reduced $g$-path in $\Gamma$. Let $u$ be the initial vertex of the path $\tau$. We will find a reduced $f$-path $\mu$ in $\Gamma$ such that $\Phi(\mu)=\tau$.
Let $\mu_1$ be an arbitrary path in $\Gamma$ from $u$ to $f(u)$. Then the paths $\tau$ and $p_ug(\overline{\mu}_1)$
have the same endpoints, so $\overline{\tau}p_ug(\overline{\mu}_1)$ is a loop based at $g(u)$. Hence, there  exists a loop $\sigma$ in $\Gamma$ based at $u$ such that $g(\sigma)=\overline{\tau}p_ug(\overline{\mu}_1)$. We set $\mu:=[\sigma\mu_1]$. Then $\mu$ is an $f$-path and $\Phi(\mu)=[p_ug(\overline{\mu})]=[p_ug(\overline{\mu}_1)g(\overline{\sigma})]=\tau$.
\\ \hspace*{144mm} $\Box$

\vspace*{-3mm}





\begin{defn}\label{defn 4.8}
{\rm The {\it inverse preferred direction} at a vertex $\mu$ in $D_f$
is the preimage of the preferred direction at the vertex $\Phi(\mu)$ in $D_g$ under $\Phi$.
}
\end{defn}

We formulate this more detailed.
Recall that $\Phi(\mu)=[p_u g(\overline{\mu})]$, where $u$ is the initial vertex of the $f$-path $\mu$.
First suppose that the $g$-path $\Phi(\mu)$ is nontrivial and
let $E$ be the first edge of this path. Then the inverse preferred direction at the vertex $\mu$ of $D_f$
is the direction of the edge of $D_f$ which starts at $\mu$ and has the label $E$.

If the $g$-path $\Phi(\mu)$ is trivial in $\Gamma$, the inverse preferred direction at $\mu$ in $D_f$ is not defined.

\begin{prop}\label{prop 4.8a} The inverse preferred direction is defined at almost all vertices of $D_f$.
\end{prop}

{\it Proof.} If the inverse preferred direction at a vertex $\mu$ in $D_f$ is not defined,
then $\Phi(\mu)$ lies in the finite set $\{\bold{1}_u \,|\, u\in \Gamma^0\}$. Since $\Phi$ is injective,
the number of such $\mu$ is finite. \hfill $\Box$

\begin{defn}\label{defn 4.9}
{\rm Preimages, with respect to $\Phi$, of repelling edges, repelling vertices and dead vertices of $D_g$ are
called {\it inv-repelling}\, edges, {\it inv-repelling} vertices and {\it inv-dead} vertices of $D_f$, respectively.}
\end{defn}

By Proposition~\ref{prop 4.1} applied to $g$,  there are only finitely many inv-repelling edges and inv-repelling and inv-dead vertices in $D_f$, and they can be algorithmically found.


\bigskip

{\bf I. Normal vertices}

\begin{defn}\label{defn 4.10}
{\rm A vertex of $D_f$ is called {\it normal} if the preferred and the inverse preferred directions at this vertex exist and do not coincide.}
\end{defn}

The main purpose of this subsection are Propositions~\ref{prop 4.11} and~\ref{prop 4.12}; they will help us to decide, whether two rays in $D_f$ (given by their initial vertices) meet.

The following lemma easily follows from Lemma~\ref{lem 3.12}.

\begin{lem}\label{Cooper_1} Let $\Gamma$ be a finite connected graph and $f:\Gamma \rightarrow \Gamma$ be a homotopy equivalence sending edges to edge paths.
Let $p$ be an initial subpath of a reduced path $q$. Write $[f(p)]\equiv ab $, where $a$ is the maximal common initial subpath of $[f(p)]$ and $[f(q)]$. Then $l(b)\leqslant C_{\star}(f)$.
\end{lem}

The source of the following lemma is Proposition (4.3) in~\cite{CL}.

\begin{lem}\label{lem_(4.3)}
Let $R$ be a $\mu$-subgraph with consecutive vertices $\mu=\mu_0,\mu_1,\dots$,
and with labels of consecutive edges $E_1,E_2,\dots $.
For each $j\geqslant 0$ with alive vertex~$\mu_j$, let $k(j)$ be the maximal natural number
such that $\mu_j\equiv E_{j+1}\dots E_{j+k(j)}\cdot Z_j$ for some~$Z_j$.
If $j>l(\mu_0)$ and $R$ has at least $j+k(j)+2$ vertices,
then $l(Z_j)\leqslant C_{\star}(f)$.




\end{lem}

\medskip

{\it Proof.} With notation $X_j=E_1E_2\dots E_j$, we have $\mu_j\equiv [\overline{X}_j\mu_0f(X_j)]$.
Hence, $f(X_j)=\overline{\mu}_0X_j\mu_j$. Therefore $[f(X_j)]\equiv [\overline{\mu}_0X_j]\cdot E_{j+1}\dots E_{j+k(j)}\cdot Z_j$. Indeed, the condition $j>l(\mu_0)$ guarantees that the last edge of $[\overline{\mu}_0X_j]$
is $E_j$ which is not inverse to $E_{j+1}$.
Applying the same arguments to $\mu_{j+k(j)}$, we have
$$
\begin{array}{ll}
[f(X_{j+k(j)})] & \equiv [\overline{\mu}_0X_{j+k(j)}]\cdot E_{j+k(j)+1}\dots E_{j+k(j)+k(j+k(j))}\cdot Z_{j+k(j)}\vspace*{2mm}\\
& \equiv [\overline{\mu}_0X_j]\cdot E_{j+1}\dots E_{j+k(j)}\cdot E_{j+k(j)+1}\dots E_{j+k(j)+k(j+k(j))}\cdot Z_{j+k(j)}.
\end{array}
$$
From Lemma~\ref{Cooper_1} applied to $X_j$ and $X_{j+k(j)}$, we deduce that $l(Z_j)\leqslant C_{\star}(f)$.\hfill $\Box$

\medskip

The source of the following lemma is Proposition (4.10) from~\cite{Turner}. The map $g$ and the constant  $K_{\star}(f)$ were defined
in~Definition~\ref{constant_K}.

\begin{lem}\label{lem_(4.10)}
Let $R$ be a $\mu$-subgraph with consecutive vertices $\mu=\mu_0,\mu_1,\dots$,
and with labels of consecutive edges $E_1,E_2,\dots $.
Suppose that $j>l(\mu_0)$ and $l(\mu_j)> C_{\star}(f)\cdot (||g||+1)+K_{\star}(f)$. If
$R$ has at least $j+k(j)+2$ vertices, then $\mu_{j+k(j)}$ is normal.
(Here $k(j)$ is as in Lemma~\ref{lem_(4.3)}.)


\end{lem}

\medskip

{\it Proof.}
It suffices to show that the first edge of the $g$-path $\Phi(\mu_{j+k(j)})$ is $\overline{E}_{j+k(j)}$.
Then, by Definition~\ref{defn 4.8}, the inv-preferred direction at $\mu_{j+k(j)}$ in $D_f$ will coincide with the direction of the edge outgoing from $\mu_{j+k(j)}$ and having the label $\overline{E}_{j+k(j)}$.
On the other hand, the (direct) preferred direction at $\mu_{j+k(j)}$ in $D_f$ coincides with the direction of the edge outgoing from $\mu_{j+k(j)}$ and having the label $E_{j+k(j)+1}$.
Since these labels do not coincide, the vertex $\mu_{j+k(j)}$ is normal.

By Lemma~\ref{lem_(4.3)}, we have
$$\mu_j\equiv E_{j+1}\dots E_{j+k(j)}\cdot Z_j\hspace*{3mm} {\text {\rm with}}\hspace*{3mm} l(Z_j)\leqslant C_{\star}(f).\eqno{(7.2)}$$
This implies
$$\mu_{j+k(j)}=\overline{E_{j+1}\dots E_{j+k(j)}}\,\,\mu_j\, f(E_{j+1}\dots E_{j+k(j)})=Z_jf(E_{j+1}\dots E_{j+k(j)}).$$

Recall that $\Phi(\mu)=[p_{\alpha(\mu)} g(\overline{\mu})]$.
Then, using (7.1), where $h=g\circ f$, we have
$$
\begin{array}{ll}
\Phi(\mu_{j+k(j)}) & \equiv [p_{\omega(E_{j+k(j)})} (g\circ f)(\overline{E_{j+1}\dots E_{j+k(j)}})\,g(\overline{Z}_j)]\vspace*{2mm}\\
& \equiv [\overline{E_{j+1}\dots E_{j+k(j)}}\, p_{\alpha(E_{j+1})}g(\overline{Z}_j)].
\end{array}
$$


From (7.2) and the assumption in this lemma, we have
$$\begin{array}{ll}
l(\overline{E_{j+1}\dots E_{j+k(j)}})=k(j) & \geqslant l(\mu_j)-C_{\star}(f)\vspace*{2mm}\\
& > K_{\star}(f)+ C_{\star}(f)\cdot||g||\vspace*{2mm}\\
& \geqslant l(p_{\alpha(E_{j+1})})+l(g(\overline{Z}_j)).
\end{array}
$$
Therefore the first edge of $\Phi(\mu_{j+k(j)})$ is
$\overline{E}_{j+k(j)}$.\hfill $\Box$

\begin{prop}\label{prop 4.11}
There exists an efficient algorithm which, given an $f$-path $\mu$,
either proves that the $\mu$-subgraph $R$ is finite or finds a normal vertex in $R$.
\end{prop}

{\it Proof.} Computing consecutive vertices of $R$, $\mu=\mu_0,\mu_1,\dots$, we either prove that $R$ is finite, or find the first $j$ with $j>l(\mu_0)$ and $l(\mu_j)>C_{\star}(f)\cdot (||g||+1)+K_{\star}(f)$. If we find such $j$, we compute $k(j)$
(note that $k(j)\leqslant l(\mu_j)$) and check, whether $\mu_0,\mu_1,\dots ,\mu_{j+k(j)+1}$ exist and different.
If the result is negative, then $R$ is finite; if positive, then the vertex $\mu_{j+k(j)}$ is normal
by Lemma~\ref{lem_(4.10)}. \hfill  $\Box$


The following proposition is contained in Claim b) in the proof of Theorem~A in ~\cite{Turner}.
This claim was inspired by Lemma (4.8) and Proposition~(4.10) from \cite{CL}.
The proof of this proposition is valid in general situation, i.e. for any homotopy equivalence $f:\Gamma\rightarrow \Gamma$ sending edges to edge paths. We give it for completeness.



\begin{prop}\label{prop 4.12} Let $R_1$ and $R_2$ be a $\mu_1$-ray and a $\mu_2$-ray in $D_f$, respectively. Suppose that they do not contain inv-repelling vertices and that their initial vertices
$\mu_1$ and $\mu_2$ are normal. Then $R_1$ and $R_2$ are either disjoint or one of them is contained in the other.
\end{prop}


{\it Proof.}
Figure~7 illustrates the proof. Suppose that the rays $R_1$ and $R_2$ intersect
and none of them is contained in the other.
We indicate the preferred directions by red triangles and the inv-preferred directions by blue triangles.
Since $\mu_1$ and $\mu_2$ are normal, the blue and the red triangles at $\mu_1$ and at $\mu_2$
look in different directions, see Figure~7~(a). Since  $R_1$ and $R_2$ do not contain inv-repelling vertices, we can
inductively reconstruct the inv-preferred directions at the vertices of $R_1$ and $R_2$ until the first
intersection point of these rays, see Figure~7~(b). We obtain two inv-preferred directions
at this point, a contradiction.\hfill $\Box$

\vspace*{-11mm}\hspace*{-13.0mm}\includegraphics[scale=0.4]{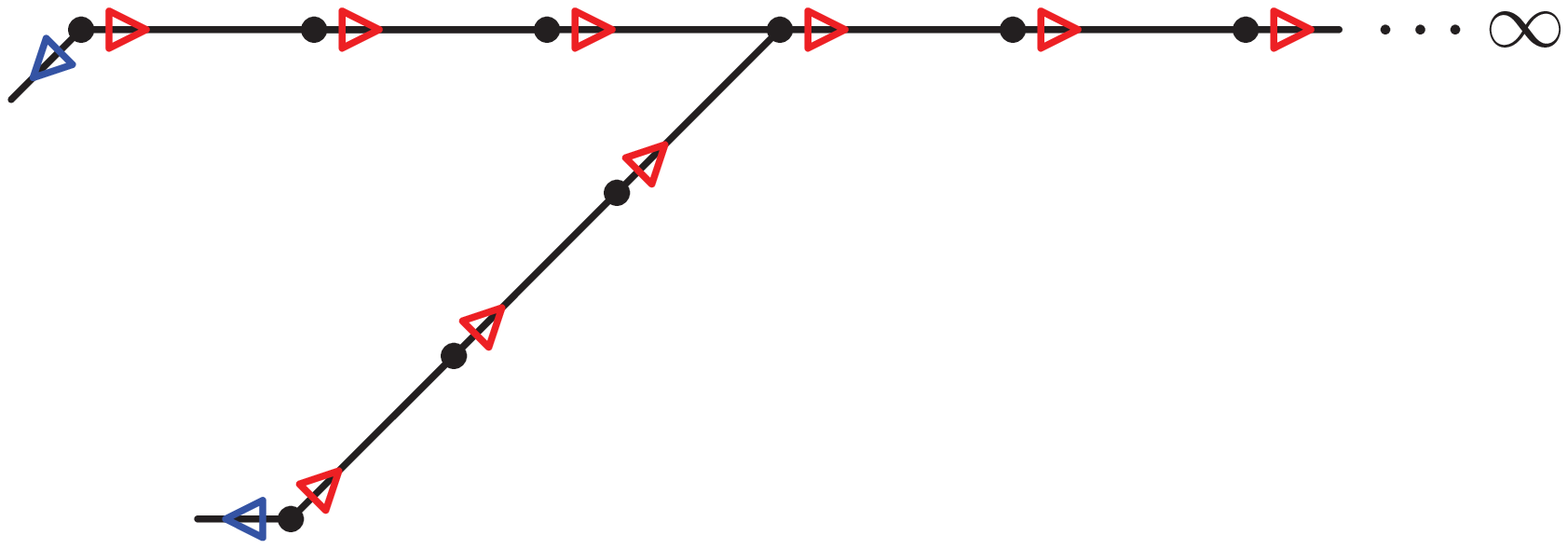}

\vspace*{-100mm}


\vspace*{-19mm}\hspace*{65.mm}\includegraphics[scale=0.4]{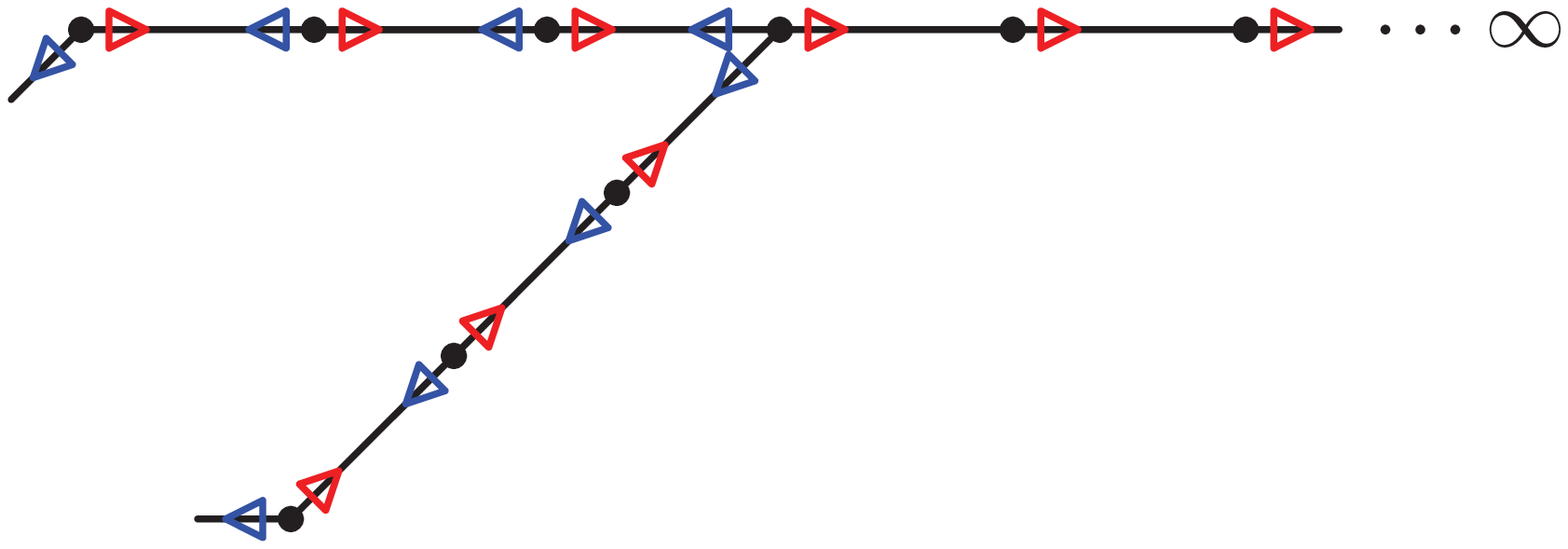}

\vspace*{-70mm}

\hspace*{29mm}(a) \hspace*{72mm}(b)

\begin{center} Figure 7.
\end{center}

{\bf J. How to convert the procedure (1)-(5) for construction of the graph $CoRe(C_f)$ into an algorithm}

As it was observed, it suffices to find algorithms for steps (2) and (4).

Using Propositions~\ref{prop 4.11} and~\ref{prop 4.12}, Step (4) can be replaced by the following three steps.


\begin{enumerate}
\item[(4.1)] For each repelling vertex $\mu$ whose $\mu$-subgraph is a ray, find in this
$\mu$-ray a vertex $\mu'$ such that the $\mu'$-ray does not contain inv-repelling
vertices.

\item[(4.2)] Find a normal vertex $\mu''$ in the $\mu'$-ray.

\item[(4.3)] For every two repelling vertices $\mu$ and $\tau$ whose $\mu$- and $\tau$-subgraphs are rays,
verify whether $\tau''$ is contained in the $\mu''$-ray
or $\mu''$ is contained in the $\tau''$-ray.
\end{enumerate}


Step (4.2) can be done algorithmically by Proposition~\ref{prop 4.11}.
Steps (4.1) and (4.3) can be done if we find an algorithm for the following problem.

\medskip


{\bf Membership problem.}
Given two vertices $\mu$ and $\tau$ of the graph $D_f$, verify
whether $\tau$ is contained in the $\mu$-subgraph.

\medskip

Indeed, for Step (4.1) we first find all inv-repelling vertices in $D_f$.
Then we detect those of them which lie in the $\mu$-ray. Let
$I$ be the minimal initial segment of the $\mu$-ray which contains all these vertices.
We can take $\mu'$ equal to the first vertex in the $\mu$-ray which lies outside $I$.
Step (4.3) is a partial case of the above problem.


Step (2) can be done if we find an algorithm for the following problem:

{\bf Finiteness problem.} Given a vertex $\mu$ of the graph $D_f$, determine whether the $\mu$-subgraph is finite or not. If the $\mu$-subgraph is finite, construct it.

\medskip

Thus, to construct $CoRe(C_f)$ algorithmically, it suffices to find algorithms for these problems.
Moreover, using the solvability of the Membership problem, we can decide whether the vertex $\mathbf{1}_{v_{\ast}}$
lies in the $\mu$-subgraph for some repelling vertex $\mu$. Then, by Proposition~\ref{prop 4.5}, we can
compute a basis of $\pi_1(D_f(\mathbf{1}_{v_{\ast}}),\mathbf{1}_{v_{\ast}})$. Lemma~\ref{lem 4.2}
identifies this group with $\overline{\text{\rm Fix}}(f)$.

In Section~20 we will present algorithms for the above problems in case where $f:\Gamma\rightarrow \Gamma$
is a PL-relative train track satisfying (RTT-iv).
As explained in Section~5, this will provide an algorithm for computing a basis of ${\text{\rm Fix}}(\varphi)$.

\section{$r$-cancelation points in paths in $G_r$}

Let $f:\Gamma\rightarrow \Gamma$ be a PL-relative train track with the maximal filtration
$\varnothing=G_0\subset \dots \subset G_N=\Gamma$.

\begin{defn}\label{defn 6.1}
{\rm
1) Let $p,q$ be reduced paths in $\Gamma$ with the same initial point. By $I(p,q)$ we denote the largest common initial subpath of $p$ and $q$. Then $p\equiv I(p,q)\cdot p'$ and $q\equiv I(p,q)\cdot q'$ for some reduced paths $p',q'$. We denote by $\Lambda(p,q)$ the ordered pair of paths $(p',q')$. These notations are motivated by the figure
$$\overset{\hspace*{0mm}\big|}{\vspace*{-1mm}\diagup\hspace*{-0.4mm}\diagdown}$$
\medskip
2) Let $\tau\equiv\bar{p}\cdot q$ be a reduced path in $\Gamma$.
For $k\geqslant 1$, we set
$$(p_k,q_k)\equiv \Lambda\bigl([f^k(p)],[f^k(q)]\bigr)\hspace*{10mm}{\text{\rm and}}\hspace*{10mm}
I_k\equiv I\bigl( [f^k(p)],[f^k(q)]\bigr)$$
Then $[f^k(\tau)]\equiv \bar{p}_k\cdot q_k$. The occurrence $y^k:=\alpha (p_k)=\alpha(q_k)$ in $[f^k(\tau)]$ is called the {\it $k$-successor} of $y:=\alpha(q)$.
}
\end{defn}

\begin{defn}\label{defn 6.3}
{\rm Let $H_r$ be an exponential stratum. Let $\tau\equiv \bar{p}\cdot q$ be a reduced path in $G_r$, where $p$ and $q$ are $r$-legal paths.

For $k\in \mathbb{N}$, let $c_k$ be the maximal initial subpath of $p$ such that
$[f^k(c_k)]$ is a subpath of $I_k$ and the terminal (possibly partial) edge of $c_k$ lies in $H_r$ if it exists.
Clearly, $c_1\subseteq c_2\subseteq \dots \subseteq p$.
Let $(p)_{\tiny {\text {\rm max}}}$ be the minimal path containing all $c_k$. Note that the terminal point of
$(p)_{\tiny {\text {\rm max}}}$ is not necessarily a vertex.
We define $(q)_{\tiny {\text {\rm max}}}$ analogously.

If $\sigma$ is a subpath of some $c_k$, then we say that $\sigma$ can be {\it virtually moved} into $I_k$ in $k$ steps. The same notion is defined for subpaths of $(q)_{\tiny {\text {\rm max}}}$.
}
\end{defn}


\begin{rmk}\label{Jahr_2014}
{\rm
(a) In Definition~\ref{defn 6.3}, we use the word {\it virtually}, since only for nontrivial subpaths
of $[f^k(\sigma)]$ which lie in $H_r$, we can guarantee that they lie in $I_k$.

(b) If $\sigma$ is a subpath of $(p)_{\tiny {\text {\rm max}}}$
with $\omega(\sigma)\neq \omega((p)_{\tiny {\text {\rm max}}})$, then there exists $k\geqslant 1$ such
that $\sigma$ can be virtually moved in $I_k$ in $k$ steps.
}
\end{rmk}


\begin{defn}\label{defn 5.1}
{\rm Let $H_r$ be an exponential stratum. Let $\tau$ be a reduced path in $G_r$.
An occurrence of a vertex $y$ in $\tau$ is called
an {\it $r$-cancelation point in} $\tau$ if $\tau$ contains a subpath $\bar{a}b$,
where $a$ and $b$ are nontrivial partial edges such that $\alpha(a)=\alpha(b)=y$
and the full edges containing $a$ and $b$ form an illegal $r$-turn.
}
\end{defn}




\begin{lem}\label{lem 6.4} Let $H_r$ be an exponential stratum.
Suppose that $\tau\equiv \bar{p}q$ is a reduced path in $G_r$
such that the paths $p$ and $q$ are $r$-legal
and the common initial point of $p$ and $q$ is an $r$-cancelation point in $\tau$.
Then the following statements hold:

\medskip


{\rm 1)} The initial and the terminal (possibly partial) edges of the paths $(p)_{\tiny {\text {\rm max}}}$
and $(q)_{\tiny {\text {\rm max}}}$ lie in $H_r$.

{\rm 2)} The number of $r$-edges in $(p)_{\tiny {\text {\rm max}}}$
and in $(q)_{\tiny {\text {\rm max}}}$, including their terminal (possibly partial) $r$-edges, is bounded from above
by a computable natural number $n_{\text{\rm critical}}$ depending only on $f$.
In particular, $L_r\bigl((p)_{\tiny {\text {\rm max}}}\bigr)<L_{\text{\rm critical}}$
and $L_r\bigl((q)_{\tiny {\text {\rm max}}}\bigr) <L_{\text{\rm critical}}$, where
$$L_{\text{\rm critical}}:=\max\{L_r(E)\,|\, E\in H_r^1\}\cdot n_{\text{\rm critical}}.$$
\end{lem}

{\it Proof.}
1) The initial (partial) edges of the paths $(p)_{\tiny {\text {\rm max}}}$
and $(q)_{\tiny {\text {\rm max}}}$ lie in $H_r$, since the common initial vertex of these paths is an $r$-cancelation point in $\tau$. The terminal (partial) edges of $(p)_{\tiny {\text {\rm max}}}$
and $(q)_{\tiny {\text {\rm max}}}$  lie in $H_r$ by definition of these paths.


2) We claim that this statement holds for $n_{\text{\rm critical}}:=2\lceil K\rceil+2$, where $K$ is the protection constant defined in the proof of~\cite[Lemma 4.2.2]{BFH}. Suppose the contrary.
Then there exists a $K$-protected $r$-edge $E$ in the interior of $(p)_{\tiny {\text {\rm max}}}$,
and by this lemma, $[f^n(E)]$ is a subpath of $[f^n(\tau)]$ for any $n\in \mathbb{N}$.
This contradicts Remark~\ref{Jahr_2014}\!~(b).
Note that $K$ is computable and depends only on $f$.
\hfill $\Box$

\begin{rmk}\label{vier_cases}
{\rm Let $H_r$ be an exponential stratum.
Let $\mu, \tau\subset G_r$ be reduced paths with $\alpha(\mu)=\alpha(\tau)$.
Let $y_1,\dots,y_k$ and $z_1,\dots ,z_s$
be consequent $r$-cancelation points in $\mu$ and in $\tau$, respectively.  Let $t$ be the terminal vertex of $I(\mu,\tau)$.
Then we have one of the following four cases for the path $\bar{\mu}\tau$ (see~Figure~8).

\vspace*{5mm}
\hspace*{15mm}
{\unitlength=1.00mm
\special{em:linewidth 0.4pt}
\linethickness{0.4pt}
\begin{picture}(93.31,75.00)
\put(4.81,55.00){\line(1,0){18.00}}
\put(22.81,55.00){\line(0,1){20.00}}
\put(24.81,75.00){\line(0,-1){20.00}}
\put(24.81,55.00){\line(1,0){18.00}}
\put(22.81,63.00){\circle*{1.00}}
\put(24.81,63.00){\circle*{1.00}}
\put(24.81,71.00){\circle*{1.00}}
\put(22.81,71.00){\circle*{1.00}}
\put(4.81,55.00){\circle*{1.00}}
\put(9.81,55.00){\circle*{1.00}}
\put(18.81,55.00){\circle*{1.00}}
\put(28.81,55.00){\circle*{1.00}}
\put(37.81,55.00){\circle*{1.00}}
\put(42.81,55.00){\circle*{1.00}}
\put(20.81,71.00){\makebox(0,0)[rc]{$y_1$}}
\put(20.81,63.00){\makebox(0,0)[rc]{$y_{\ell}$}}
\put(20.81,68.00){\makebox(0,0)[rc]{$\vdots$}}
\put(26.81,68.00){\makebox(0,0)[lc]{$\vdots$}}
\put(26.81,71.00){\makebox(0,0)[lc]{$z_1$}}
\put(26.81,63.00){\makebox(0,0)[lc]{$z_{\ell}$}}
\put(8.81,53.00){\makebox(0,0)[ct]{$y_k$}}
\put(13.81,53.00){\makebox(0,0)[ct]{$\cdots$}}
\put(19.81,53.00){\makebox(0,0)[ct]{$y_{{\ell}+1}$}}
\put(28.81,53.00){\makebox(0,0)[ct]{$z_{{\ell}+1}$}}
\put(34.81,53.00){\makebox(0,0)[ct]{$\cdots$}}
\put(38.81,53.00){\makebox(0,0)[ct]{$z_s$}}
\put(23.81,55.00){\circle{4.00}}
\put(54.81,55.00){\line(1,0){18.00}}
\put(72.81,55.00){\line(0,1){20.00}}
\put(74.81,75.00){\line(0,-1){20.00}}
\put(74.81,55.00){\line(1,0){18.00}}
\put(72.81,63.00){\circle*{1.00}}
\put(74.81,63.00){\circle*{1.00}}
\put(74.81,71.00){\circle*{1.00}}
\put(72.81,71.00){\circle*{1.00}}
\put(54.81,55.00){\circle*{1.00}}
\put(59.81,55.00){\circle*{1.00}}
\put(68.81,55.00){\circle*{1.00}}
\put(78.81,55.00){\circle*{1.00}}
\put(87.81,55.00){\circle*{1.00}}
\put(92.81,55.00){\circle*{1.00}}
\put(70.81,71.00){\makebox(0,0)[rc]{$y_1$}}
\put(70.81,63.00){\makebox(0,0)[rc]{$y_{{\ell}-1}$}}
\put(70.81,68.00){\makebox(0,0)[rc]{$\vdots$}}
\put(76.81,68.00){\makebox(0,0)[lc]{$\vdots$}}
\put(76.81,71.00){\makebox(0,0)[lc]{$z_1$}}
\put(76.81,63.00){\makebox(0,0)[lc]{$z_{{\ell}-1}$}}
\put(58.81,53.00){\makebox(0,0)[ct]{$y_k$}}
\put(63.81,53.00){\makebox(0,0)[ct]{$\cdots$}}
\put(69.81,53.00){\makebox(0,0)[ct]{$y_{{\ell}+1}$}}
\put(78.81,53.00){\makebox(0,0)[ct]{$z_{{\ell}+1}$}}
\put(84.81,53.00){\makebox(0,0)[ct]{$\cdots$}}
\put(88.81,53.00){\makebox(0,0)[ct]{$z_s$}}
\put(73.81,55.00){\circle{4.00}}
\put(4.89,15.00){\line(1,0){18.00}}
\put(22.89,15.00){\line(0,1){20.00}}
\put(24.89,35.00){\line(0,-1){20.00}}
\put(24.89,15.00){\line(1,0){18.00}}
\put(22.89,23.00){\circle*{1.00}}
\put(24.89,23.00){\circle*{1.00}}
\put(24.89,31.00){\circle*{1.00}}
\put(22.89,31.00){\circle*{1.00}}
\put(4.89,15.00){\circle*{1.00}}
\put(9.89,15.00){\circle*{1.00}}
\put(18.89,15.00){\circle*{1.00}}
\put(28.89,15.00){\circle*{1.00}}
\put(37.89,15.00){\circle*{1.00}}
\put(42.89,15.00){\circle*{1.00}}
\put(20.89,31.00){\makebox(0,0)[rc]{$y_1$}}
\put(20.89,23.00){\makebox(0,0)[rc]{$y_{{\ell}-1}$}}
\put(20.89,28.00){\makebox(0,0)[rc]{$\vdots$}}
\put(26.89,28.00){\makebox(0,0)[lc]{$\vdots$}}
\put(26.89,31.00){\makebox(0,0)[lc]{$z_1$}}
\put(26.89,23.00){\makebox(0,0)[lc]{$z_{{\ell}-1}$}}
\put(8.89,13.00){\makebox(0,0)[ct]{$y_k$}}
\put(13.89,13.00){\makebox(0,0)[ct]{$\cdots$}}
\put(19.89,13.00){\makebox(0,0)[ct]{$y_{{\ell}+1}$}}
\put(28.89,13.00){\makebox(0,0)[ct]{$z_{\ell}$}}
\put(34.89,13.00){\makebox(0,0)[ct]{$\cdots$}}
\put(38.89,13.00){\makebox(0,0)[ct]{$z_s$}}
\put(23.89,15.00){\circle{4.00}}
\put(72.81,56.00){\makebox(0,0)[rb]{$y_{\ell}$}}
\put(74.81,56.00){\makebox(0,0)[lb]{$z_{\ell}$}}
\put(22.89,16.00){\makebox(0,0)[rb]{$y_{\ell}$}}
\put(23.81,41.00){\makebox(0,0)[cc]{(1)}}
\put(73.81,41.00){\makebox(0,0)[cc]{(2)}}
\put(23.89,1.00){\makebox(0,0)[cc]{(3)}}
\put(72.81,55.00){\circle*{1.00}}
\put(74.81,55.00){\circle*{1.00}}
\put(22.89,15.00){\circle*{1.00}}
\
\put(13.81,48.00){\makebox(0,0)[cc]{$\underbrace{\hspace{1.9cm}}$}}
\put(13.81,45.00){\makebox(0,0)[cc]{$\mu$}}
\put(33.81,48.00){\makebox(0,0)[cc]{$\underbrace{\hspace{1.9cm}}$}}
\put(33.81,45.00){\makebox(0,0)[cc]{$\tau$}}
\
\put(63.81,48.00){\makebox(0,0)[cc]{$\underbrace{\hspace{1.9cm}}$}}
\put(63.81,45.00){\makebox(0,0)[cc]{$\mu$}}
\put(83.81,48.00){\makebox(0,0)[cc]{$\underbrace{\hspace{1.9cm}}$}}
\put(83.81,45.00){\makebox(0,0)[cc]{$\tau$}}
\
\put(13.89,8.00){\makebox(0,0)[cc]{$\underbrace{\hspace{1.9cm}}$}}
\put(13.89,5.00){\makebox(0,0)[cc]{$\mu$}}
\put(33.89,8.00){\makebox(0,0)[cc]{$\underbrace{\hspace{1.9cm}}$}}
\put(33.89,5.00){\makebox(0,0)[cc]{$\tau$}}
\
\put(54.75,15.00){\line(1,0){18.00}}
\put(72.75,15.00){\line(0,1){20.00}}
\put(74.75,35.00){\line(0,-1){20.00}}
\put(74.75,15.00){\line(1,0){18.00}}
\put(72.75,23.00){\circle*{1.00}}
\put(74.75,23.00){\circle*{1.00}}
\put(74.75,31.00){\circle*{1.00}}
\put(72.75,31.00){\circle*{1.00}}
\put(54.75,15.00){\circle*{1.00}}
\put(59.75,15.00){\circle*{1.00}}
\put(68.75,15.00){\circle*{1.00}}
\put(78.75,15.00){\circle*{1.00}}
\put(87.75,15.00){\circle*{1.00}}
\put(92.75,15.00){\circle*{1.00}}
\put(70.75,31.00){\makebox(0,0)[rc]{$y_1$}}
\put(70.75,23.00){\makebox(0,0)[rc]{$y_{{\ell}-1}$}}
\put(70.75,28.00){\makebox(0,0)[rc]{$\vdots$}}
\put(76.75,28.00){\makebox(0,0)[lc]{$\vdots$}}
\put(76.75,31.00){\makebox(0,0)[lc]{$z_1$}}
\put(76.75,23.00){\makebox(0,0)[lc]{$z_{{\ell}-1}$}}
\put(58.75,13.00){\makebox(0,0)[ct]{$y_k$}}
\put(63.75,13.00){\makebox(0,0)[ct]{$\cdots$}}
\put(69.75,13.00){\makebox(0,0)[ct]{$y_{\ell}$}}
\put(78.75,13.00){\makebox(0,0)[ct]{$z_{{\ell}+1}$}}
\put(84.75,13.00){\makebox(0,0)[ct]{$\cdots$}}
\put(88.75,13.00){\makebox(0,0)[ct]{$z_s$}}
\put(73.75,15.00){\circle{4.00}}
\put(73.75,1.00){\makebox(0,0)[cc]{(4)}}
\put(63.75,8.00){\makebox(0,0)[cc]{$\underbrace{\hspace{1.9cm}}$}}
\put(63.75,5.00){\makebox(0,0)[cc]{$\mu$}}
\put(83.75,8.00){\makebox(0,0)[cc]{$\underbrace{\hspace{1.9cm}}$}}
\put(83.75,5.00){\makebox(0,0)[cc]{$\tau$}}
\put(74.60,16.00){\makebox(0,0)[lb]{$z_{\ell}$}}
\put(74.60,15.00){\circle*{1.00}}
\end{picture}
}

\vspace*{3mm}

\begin{center}
Figure 8.
\end{center}

\vspace*{2mm}

(1) Suppose that $t\notin \{y_1,\dots ,y_k,z_1,\dots ,z_s\}$ and that
$y_1=z_1,\dots , y_{\ell}=z_{\ell}$, are all $r$-cancelation points in $I(\mu,\tau)$.
Then all $r$-cancelation points in $[\bar{\mu}\tau]$ are $y_k,\dots,y_{\ell+1},z_{\ell+1},\dots,z_s$ and, possibly,
$t$ (see Figure~8 (1)).


(2) Suppose that $t=y_{\ell}=z_{\ell}$ for some ${\ell}$. Then $t$ is an $r$-cancelation point in $[\bar{\mu}\tau]$. This follows
from the fact that the turns at $t$ in $\mu$ and in $\tau$ are illegal and have a common (partial) edge.
So, all $r$-cancelation points in $[\bar{\mu}\tau]$ are $y_k,\dots,y_{\ell+1},t,z_{\ell+1},\dots,z_s$ (see Figure~8 (2)).


(3) Suppose that $t=y_{\ell}$ for some ${\ell}$ and that $t\notin \{z_1,\dots, z_s\}$. Then $y_1=z_1,\dots,y_{\ell-1}=z_{\ell-1}$
are all $r$-cancelation points in the path $I(\mu,\tau)$. Moreover, $t$ is not an $r$-cancelation point
in $[\bar{\mu}\tau]$. This follows from the fact that the turn at $t$ in $\mu$ is $r$-illegal,
while the turn at $t$ in $\tau$ is not and these turns have a common (partial) edge.
So, all $r$-cancelation points in $[\bar{\mu}\tau]$ are $y_k,\dots,y_{\ell+1},z_{\ell},\dots,z_s$ (see Figure~8~\!(3)).


(4) Suppose that $t=z_{\ell}$ for some ${\ell}$ and that $t\notin \{y_1,\dots, y_k\}$ (see Figure~8 (4)).
This case is analogous to Case (3). In particular, $t$ is not an $r$-cancelation point in~$[\bar{\mu} \tau]$.

\medskip

{\sl Thus, the maximal possible number of $r$-cancelation points in $[\bar{\mu}\tau]$ is $k+s+1$ (it can be achieved only in Case (1) for $\ell=0$). If at least one of the $r$-cancelation points in $\mu$ or in
$\tau$ lies in $I(\mu,\tau)$, then the number of $r$-cancelation points in $[\bar{\mu}\tau]$ is strictly less
than $k+s$.}
}
\end{rmk}

\section{$r$-cancelation areas}

Let $f:\Gamma\rightarrow \Gamma$ be a PL-relative train track with the maximal filtration
$\varnothing=G_0\subset \dots \subset G_N=\Gamma$.

\begin{defn}\label{defn 6.2} {\rm Let $H_r$ be an exponential stratum. Suppose that $\tau\equiv \bar{p}q$ is a reduced path in $G_r$
such that the paths $p$ and $q$ are $r$-legal
and the common initial point $y$ of $p$ and $q$ is an $r$-cancelation point in $\tau$.
We say that $y$ is {\it non-deletable} in $\tau$
if for every $k\geqslant 1$ the $k$-successor $y^k$ is an $r$-cancelation point in $[f^k(\tau)]$.
We say that $y$ is {\it deletable} in $\tau$ if this does not hold.

If $y$ is non-deletable in $\tau$, we call the path $A:=\overline{(p)}_{\tiny {\text {\rm max}}}  (q)_{\tiny {\text {\rm max}}}$
the {\it $r$-cancelation area} (in $\tau$).
The number $a:=L_r((p)_{\tiny {\text {\rm max}}})=L_r((q)_{\tiny {\text {\rm max}}})$ is called the {\it $r$-cancelation radius} of $A$.
}
\end{defn}

\begin{rmk}\label{rmk 6.2a} {\rm Let $\tau$ be a reduced path in $G_r$ with a unique $r$-cancelation point $y$.
Then the point $y$ is a non-deletable $r$-cancelation point in~$\tau$ if and only if $y^k$ is a non-deletable $r$-cancelation point in $[f^k(\tau)]$ for some natural $k$.
}
\end{rmk}

\begin{lem}\label{properties_areas}
{\rm Each $r$-cancelation area $A$ satisfies the following properties:

\begin{enumerate}
\item[1)]  Each $[f^k(A)]$ is an $r$-cancelation area. In particular, each $[f^k(A)]$ contains exactly one
$r$-cancelation point.

\item[2)]  The initial and the terminal (possibly partial) edges of each $[f^k(A)]$ are contained in $H_r$.

\item[3)] The number of $r$-edges in $[f^k(A)]$ is bounded independently of $k$.
\end{enumerate}
}
\end{lem}

{\it Proof.} 1) follows from the above definition, 2) and 3) from Lemma~\ref{lem 6.4}. \hfill $\Box$

\begin{rmk}
{\rm
1) The set of $r$-cancelation areas coincides with the set $P_r$ defined before Lemma~4.2.5 in~\cite{BFH}.

2) Let $E$ be an edge of $\Gamma$. We write $f(E)\equiv E_1E_2\dots E_k$, where $E_1,\dots, E_k$ are edges of $\Gamma$.
Then there is a subdivision $E\equiv e_1e_2\dots e_k$ such that $f(e_i)\equiv E_i$. By assumption, the numbers $l(e_i)$ for each $i$ are given and the restriction of $f$ to each $e_i$ is linear. Therefore, for each given $m\in \mathbb{N}$, we can compute the $l$-lengths of the closures of connected components of $E\setminus f^{-m}(\Gamma^0)$.
}
\end{rmk}

\begin{prop}\label{number M} Let $H_r$ be an exponential stratum.

{\rm 1)} Given two $r$-legal paths $\beta,\gamma$ in $G_r$ with $L_r(\beta)>0$ and $L_r(\gamma)>0$,
there exists at most one $r$-cancelation area $A$ such that
$\bar{\beta}$ and $\gamma$ are some initial and terminal subpaths of $A$.

{\rm 2)}
The number of $r$-cancelation areas is at most $M_r:=m_r^2n_r^2$, where $m_r$ is the number of edges in $H_r$ and $n_r$ is the number of sequences $(p_1,p_2,\dots,p_s)$, where all $p_i$ are $r$-legal edge paths in $H_r$
with $\sum_{i=1}^s L_r(p_i)\leqslant L_{\text{\rm critical}}$, $s\in \{0\}\cup \mathbb{N}$.
\end{prop}

\medskip

{\it Proof.}  1) Without loss of generality, we may assume that $L_r(\beta)=L_r(\gamma)$.
Suppose that $A\equiv \bar{p}q$ is an $r$-cancelation area, where $p$ and $q$
are $r$-legal and $\beta$ and $\gamma$ are terminal subpaths of $p$ and $q$, respectively.
Let $k$ be the minimal natural number such that $$\lambda_r^k\cdot L_r(\beta)>L_{\text{\rm critical}}.$$
Then $L_r([f^k(\beta)]=L_r([f^k(\gamma)])>L_{\text{\rm critical}}$. This implies that
$[f^k(\bar{p}q)]=[\bar{b}c]$, where $b$ is obtained from $[f^k(\beta)]$ by deleting the maximal initial subpath
lying in $G_{r-1}$ and $c$ is obtained analogously from $[f^k(\gamma)]$.
Hence $[f^k(\bar{p}q)]$ and so $\bar{p}q$ are completely determined by $\beta$ and $\gamma$.

2) First we introduce notations. For any reduced path $\tau$ in $G_r$, we can write $\tau\equiv c_0\tau_1c_1\dots \tau_sc_s$,
where the paths $c_1,c_2,\dots,c_{s-1}$ lie in $G_{r-1}$ and are nontrivial, the paths $\tau_1,\tau_2,\dots,\tau_s$ lie in $H_r$ and are nontrivial, and $c_0,c_s$ lie in~$G_{r-1}$ or are trivial. We denote $\tau\cap H_r:=(\tau_1,\tau_2,\dots,\tau_s)$.

Let $\tau_s'$ is obtained from $\tau_s$ by deleting the terminal partial edge of $\tau_s$ if it exists.
We set $\lfloor\tau\cap H_r\rfloor:=(\tau_1,\dots ,\tau_{s-1},\tau_s')$ if $\tau_s'$ is not empty and $\lfloor\tau\cap H_r\rfloor:=(\tau_1,\dots,\tau_{s-1})$ if $\tau_s'$ is empty.

The following claim is proven in the proof of Lemma~4.2.5 in~\cite{BFH}:

For any two sequences $\mu:=(\mu_1,\mu_2,\dots ,\mu_s)$, $\sigma:=(\sigma_1,\sigma_2,\dots ,\sigma_t)$ where $\mu_1,\dots,\mu_s$, $\sigma_1,\dots,\sigma_t$
are $r$-legal edge paths in $H_r$, and for any two edges $E_1,E_2$ in $H_r$, there exists at most one $r$-cancelation area $A\equiv\bar{p}q$ such that the paths $p$ and $q$ are $r$-legal, $\lfloor p\cap H_r \rfloor =\mu$, $\lfloor q\cap H_r \rfloor =\sigma$, and the terminal (possibly partial) edge of $p$ is a part of $E_1$,
and the terminal (possibly partial) edge of $q$ is a part of $E_2$.

Clearly, this claim and Lemma~\ref{lem 6.4}.~\!2) imply the statement 2).
\hfill $\Box$

\begin{defn} {\rm Let $H_r$ be an exponential stratum. Let $x$ be a point in an $r$-edge $E$. The {\it $(l,L_r)_E$-coordinates} of $x$ is the pair
$(l(p),L_r(p))$, where $p$ is the initial segment of $E$ with $\omega(p)=x$.}
\end{defn}

\begin{lem}\label{fm_periodic} Let $H_r$ be an exponential stratum. For each $r$-edge $E$ and each $m\in \mathbb{N}$,  the set $\{x\in E\,|\, f^m(x)=x\}$ is finite. Given such $E$ and $m$, we can efficiently compute the set of $(l,L_r)_E$-coordinates of all points of this set.
\end{lem}

{\it Proof.} We show how to compute the $(l,L_r)_E$-coordinates of points of the set $$V:=\{x\in E\,|\, f^m(x)=x\}\setminus \{\alpha(E), \omega(E)\}.$$
Let $x\in V$ and let $p$ be the initial segment of $E$ with $\omega(p)=x$.
Write $f^m(E)\equiv E_1E_2\dots E_k$, where all $E_i$ are edges and write $E\equiv e_1e_2\dots e_k$, where $f^m(e_i)=E_i$ for $i=1,\dots,k$.
Note that all $l(e_i)$ can be computed.
Suppose that $x$ is contained in $e_i$. Then $E_i=E$ or $E_i=\bar{E}$.
First we consider the case $E_i=E$.
We have
$$
l(p)=l_{e_i}(\alpha(e_i),x)+\overset{i-1}{\underset{j=1}{\sum}} l(e_j).\eqno{(9.1)}$$
Using $E=E_i$ and the linearity of $f|_{e_i}:e_i\rightarrow E_i$, we obtain
$$l(p)=l_{E_i}(\alpha(E_i),x)=\frac{l(E_i)}{l(e_i)}\cdot l_{e_i}(\alpha(e_i),x).\eqno{(9.2)}
$$
Then, using $l(E_i)=1$, we deduce from (9.1) and (9.2)
$$l(p)=\frac{1}{1-l(e_i)}\cdot \overset {i-1}{\underset{j=1}{\sum}} l(e_j).\eqno{(9.3)}$$
It is easy to see that $$\lambda_r^m\cdot L_r(p)=L_r(p)+\overset {i-1}{\underset{j=1}{\sum}} L_r(E_j),$$
hence $$L_r(p)=\frac{1}{\lambda_r^m-1}\cdot \overset {i-1}{\underset{j=1}{\sum}} L_r(E_j).\eqno{(9.4)}$$
In case $E_i=\bar{E}$ we have the following analogous of (9.3) and (9.4):
$$l(p)=\frac{1}{1+l(e_i)}\cdot \overset i{\underset{j=1}{\sum}}\, l(e_j),\eqno{(9.3')}$$
$$L_r(p)=\frac{1}{\lambda_r^m+1}\cdot \overset i{\underset{j=1}{\sum}} L_r(E_j).\eqno{(9.4')}$$
\hfill $\Box$

\begin{thm}\label{find_areas}
There is an efficient algorithm finding all $r$-cancelation areas of $f$.
\end{thm}

\medskip

{\it Proof.} Let $\mathcal{A}$ be the set of all $r$-cancelation areas.
Let $U$ be the set of all endpoints of all $r$-cancelation areas. The set $U$ is $f$-invariant and lies in $H_r$ by Lemma~\ref{properties_areas}, and $|U|\leqslant M_r$ by Proposition~\ref{number M}. We consider the subset $U':=\{f^{M_r}(u)\,|\, u\in U\}$ of $U$. Then each point of $U'$ is fixed by $f^m$ for some $0<m\leqslant M_r$.
Therefore $U'$ is contained in the set
$$\overline{U'}:=\bigcup_{E\in H_r^1}\bigcup_{m=1}^{M_r}\,\{x\in E\,|\, f^m(x)=x\}$$
that is finite and computable by Lemma~\ref{fm_periodic}.
Then the set $$\overline{U}:=\{u\in H_r\,|\, f^{M_r}(u)\in \overline{U'}\}$$ is finite and computable, and contains $U$.

Let $\mathcal{P}$ be the set of nontrivial initial segments $\rho$ of $r$-edges with $\omega(\rho)\in \overline{U}$.
Suppose that $A$ is an $r$-cancelation area. We write $A\equiv \bar{p}q$, where $p$ and $q$ are $r$-legal paths. Then $p$ and $q$ have terminal segments $\beta$ and $\gamma$, respectively, which lie in $\mathcal{P}$. By Proposition~\ref{number M}.1),
$A$ is completely determined by $\beta$ and $\gamma$. The proof of this proposition gives us the following algorithm constructing all elements of $\mathcal{A}$:





\begin{enumerate}

\item[(a)]  Compute $\mathcal{L}=\min \{L_r(\rho)|\,\rho\in \mathcal{P}\}$ and the minimal $k\in \mathbb{N}$ such that $\lambda_r^k\cdot \mathcal{L}>L_{\text{\rm critical}}$. Denote $\mathcal{A}_k:=\{[f^k(A)]\,|\, A\in \mathcal{A}\}$. Clearly, $\mathcal{A}_k\subseteq \mathcal{A}$.

\item[(b)] Compute the set $\Psi_k$ of all paths of the form $\bar{b}c$, where $\alpha(b)=\alpha(c)$ is a vertex,
$b$ and $c$ are nontrivial terminal subpaths of $[f^k(\beta)]$ and of $[f^k(\gamma)]$ for some $\beta,\gamma \in \mathcal{P}$, $L_r(b)=L_r(c)\leqslant L_{\text{\rm critical}}$, and
the first (possibly partial) edges of $b$ and $c$ form a nondegenerate illegal $r$-turn.
Then $\mathcal{A}_k\subseteq \Psi_k$.

\item[(c)] Compute the set $\Psi$ of reduced paths $d\subset G_r$ such that $[f^k(d)]\in \Psi_k$ and $d$ contains
exactly one $r$-cancelation point. Then $\mathcal{A}\subseteq \Psi$.

\item[(d)] Compute the set $\widetilde{\Psi}=\{\tau\in \Psi\,|\, [f^i(\tau)]\in \Psi,\, i=1,\dots ,|\Psi|\}$. This is possible since $\Psi$ is finite. Then $\mathcal{A}=\widetilde{\Psi}$.
\hfill $\Box$

\end{enumerate}

\section{Stable paths and their $A$-decompositions}

Let $f:\Gamma\rightarrow \Gamma$ be a PL-relative train track with the maximal filtration
$\varnothing=G_0\subset \dots \subset G_N=\Gamma$.
We use notations from Definition~\ref{defn 6.1}.
The following lemma is an immediate consequence of Definition~\ref{defn 6.3}.

\begin{lem}\label{lem 6.8}
Let $H_r$ be an exponential stratum in $\Gamma$.
Let $\tau\equiv \bar{p}_0q_0$ be a reduced path in $G_r$ such that the paths $p_0$ and $q_0$ are $r$-legal,
and the common initial point of $p_0$ and $q_0$ is an  $r$-cancelation point in $\tau$.
Suppose that, for some $k\geqslant 1$, the path $[f^{k}(\tau)]\equiv \bar{p}_kq_k$ has an $r$-cancelation point and this point is the common initial point of $p_k$ and $q_k$. Then
$$
\bigl((p_{k})_{\tiny {\text {\rm max}}},
(q_{k})_{\tiny {\text {\rm max}}}\bigr) = \Lambda
\bigl([f^k((p_0)_{\tiny {\text {\rm max}}})],
[f^k((q_{0})_{\tiny {\text {\rm max}}})]\bigr).
$$
\end{lem}

\begin{prop}\label{prop 6.5} Let $H_r$ be an exponential stratum in $\Gamma$. There exists a computable integer number $T\geqslant 0$ (depending only on $f$) with the following property:

Suppose that $\tau\equiv \bar{p}_0q_0$ is a reduced edge path in $G_r$
such that the paths $p_0$ and~$q_0$ are $r$-legal,
and the common initial point of $p_0$ and $q_0$ is an  $r$-cancelation point in~$\tau$.
Then either $[f^{T}(\tau)]$ is $r$-legal, or $y$ is a non-deletable $r$-cancelation point in~$\tau$.


\end{prop}

{\it Proof.} We define the following constants:

\bigskip

\hspace*{-5mm}\framebox[145mm][l]
{\noindent
\begin{minipage}[t]{150mm}
$k_0$ is the minimal natural number such that
$\lambda^{k_0}\min \{L_r(E)\,|\, E\in H_r^1\}> L_{\text{\rm critical}}$.

\noindent
$N_0$ is the number of edge paths in $G_r$ of $l$-length at most $||f||^{k_{0}}$.

\noindent
$k_1:=k_0+N_0^2$.

\noindent
$N_1$ is the number of edge paths in $G_r$ of $l$-length at most $||f||^{k_1}n_{\text{\rm critical}}$.

\noindent
$k_2:=k_1+N_0\cdot N_1$.

\noindent
$N_2$ is the number of edge paths in $G_r$ of $l$-length at most $||f||^{k_2}n_{\text{\rm critical}}$.
\end{minipage}
}

\bigskip

We prove that the proposition is valid for $T:=k_2+N_2^2$. Suppose that
$[f^T(\tau)]$ is not $r$-legal. Then

{\sl $[f^i(\tau)]$ contains an $r$-cancelation point for each $0\leqslant i\leqslant T$.}\hfill (10.1)

\medskip

We use the following definition: Let $\mu$ be a path in $G_r$. Then any maximal nontrivial subpath of $\mu$
which lies in $G_{r-1}$ is called a {\it $G_{r-1}$-piece} of $\mu$.

By Lemma~\ref{lem 6.4}.~\!1), we can write $(p_0)_{\tiny {\text {\rm max}}}\equiv a_1b_1\dots a_sb_sa_{s+1}$ and $(q_0)_{\tiny {\text {\rm max}}}\equiv a_1'b_1'\dots a_t'b_t'a_{t+1}'$,
where $s\geqslant 0$, $t\geqslant 0$, all paths $a_i$ and $a_i'$ are nontrivial and lie in $H_r$, and all paths $b_i$ and~$b_i'$ are
nontrivial and lie in $G_{r-1}$.


Recall  that $[f^k(p_0)]\equiv I_kp_k$  and $[f^k(q_0)]\equiv I_kq_k$, where $I_{k}=I([f^{k}(p_0)],[f^{k}(q_0)])$.
The initial (but not necessarily the terminal) point of $(p_k)_{\tiny {\text {\rm max}}}$ is a vertex.
Let $(p_k)_{\tiny {\text {\rm e-max}}}$ be the minimal edge path containing $(p_k)_{\tiny {\text {\rm max}}}$.
We define $(q_k)_{\tiny {\text {\rm e-max}}}$ analogously.

\medskip

{\tt Claim 1.} For each $k\geqslant k_0$ such that $[f^k(\tau)]$ contains an $r$-cancelation point,
and for each $b_i$, $i=1,\dots ,s$, one of the following is satisfied:

\begin{enumerate}

\item[1)] $[f^k(b_i)]$ is a subpath of $I_k$.

\item[2)] $[f^k(b_i)]$ is a $G_{r-1}$-piece of $p_k$. Moreover, $p_k\equiv p_{k,i,1}[f^k(b_i)]p_{k,i,2}$
for some paths $p_{k,i,1}$, $p_{k,i,2}$ such that $0<l(p_{k,i,1})\leqslant ||f||^{k_0}$.
    \\ (We stress that the last number does not depend on $k$.)

\end{enumerate}

The same alternative holds for each $b'_{j}$, $j=1,\dots,t$.

\medskip

{\it Proof.}
We fix $k$ and $i$ in the above intervals.
Let $e$ be the last edge of $a_i$, so $e$ is an $r$-edge and $\omega(e)=\alpha(b_i)$.
Let $E$ be the last edge of $[f^{k-k_0}(e)]$. Then $E$ satisfies the following properties:

$\bullet$ $E$ is an $r$-edge and $\omega([f^{k_0}(E)])=\alpha([f^k(b_i)])$.

$\bullet$ $[f^{k_0}(E)]$ is a subpath of $[f^{k}(e)]$, which is a subpath of $[f^{k}(p_0)]\equiv I_kp_k$.

$\bullet$ $[f^{k_0}(E)]$ does not lie in $p_k$.\hfill (10.2)

\medskip

The first two properties are evident. We prove the last one. Suppose that $[f^{k_0}(E)]$ lies in $p_k$.
Since $\omega(e)\neq \omega((p_0)_{\tiny {\text {\rm max}}})$,
$e$ can be virtually moved into some $I_{\ell}$ in $\ell$ steps.
Then $[f^{k_0}(E)]$ can be virtually moved into $I_{k+\ell}$ in $\ell$ steps.
Therefore $[f^{k_0}(E)]$ lies in~$(p_k)_{\tiny {\text {\rm max}}}$.
From this we have $L_r([f^{k_0}(E)])\leqslant L_r((p_k)_{\tiny {\text {\rm max}}})<L_{\text{\rm critical}}$ by Lemma~\ref{lem 6.4},
that contradicts the definition of $k_0$.
Now we are ready to finish the proof of Claim~1.

 \medskip

 {\it Case 1.} Suppose that $[f^{k_0}(E)]$ lies in $I_k$. Then $\alpha([f^k(b_i)])$ lies in $I_k$.
 We prove that the whole path $[f^k(b_i)]$ lies in $I_k$.
 Suppose the contrary, then $[f^k(b_i)]$ covers the first edge of $p_k$, hence this edge lies in $G_{r-1}$,
 a contradiction to the assumption that $[f^k(\tau)]$ contains an $r$-cancelation point.
 Thus, in this case the statement 1)~holds.

 {\it Case 2.} Suppose that $[f^{k_0}(E)]$ does not lie in $I_{k}$. By (10.2), $[f^{k_0}(E)]$ does not lie in $p_k$. Hence the first edge of $[f^{k_0}(E)]$ lies in $I_k$ and the last edge lies in $p_k$.
Since $\omega([f^{k_0}(E)])=\alpha([f^{k}(b_i)])$, we can write
$p_k\equiv p_{k,i,1}[f^k(b_i)]p_{k,i,2}$ for some paths $p_{k,i,1}$,~$p_{k,i,2}$. Moreover,
$$0<l(p_{k,i,1})\leqslant l([f^{k_0}(E)])\leqslant ||f||^{k_0}.$$
Thus, in this case the statement 2) is valid. \hfill $\Box$

\bigskip

{\tt Claim 2.} If both paths $b_{s}$ and  $b_t'$ exist, then  at least one of them
can be virtually moved into $I_{k_1}$ in $k_1$ steps.

\medskip

{\it Proof.} 
Suppose the contrary.
Then, for each $k=k_0,k_0+1,\dots ,k_1$, and for $i=s$ (for $j=t$) the second statement of Claim 1 is valid:

\medskip

\noindent
$p_k\equiv X_k[f^k(b_s)]Y_k$,\hspace*{2mm}  where $0<l(X_k)\leqslant ||f||^{k_0}$ and $[f^k(b_s)]$ is a $G_{r-1}$-piece of $p_k$;

\medskip
\noindent
$q_k\equiv X_k'[f^k(b'_t)]Y_k'$,\hspace*{2mm} where $0<l(X_k')\leqslant ||f||^{k_0}$ and $[f^k(b'_t)]$ is a $G_{r-1}$-piece of $q_k$.

\medskip

\noindent
Then the pairs of paths $(X_k,X_k')$ repeat, hence, for any $\ell\geqslant 1$,
the paths $b_s$ and $b_t'$ cannot be  virtually moved into $I_{\ell}$ in $\ell$ steps.
This contradicts Remark~\ref{Jahr_2014}\!~(b).
\hfill$\Box$



\medskip

{\tt Claim 3.} Let $(p_{0})_{\text{\rm max}}\equiv X_0b_jY_0$, where $b_j$ is a $G_{r-1}$-piece of $(p_{0})_{\text{\rm max}}$.
Suppose that $b_j$ cannot be virtually moved into $I_{k_1}$ in $k_1$ steps. Then the following statements are satisfied.

\begin{enumerate}
\item[1)] $(p_{k_1})_{\rm max}\equiv X_{k_1}[f^{k_1}(b_j)]Y_{k_1}$ for some paths $X_{k_1}$, $Y_{k_1}$ such that $$0<l(X_{k_1})\leqslant ||f||^{k_0}.$$
          Moreover, $[f^{k_1}(b_j)]$ is a $G_{r-1}$-piece of $(p_{k_1})_{\rm max}$.
          \vspace*{1mm}

\item[2)] $l((q_{k_1})_{\text{\rm max}})\leqslant ||f||^{k_1}n_{\text{\rm critical}}.$
\end{enumerate}

\medskip

{\it Proof.}  The first statement follows from Claim~1. We prove the second one.
First we consider the case where $b_t'$ does not exist, i.e. $(q_0)_{\rm max}\equiv a_1'$.
By Lemma~\ref{lem 6.8},  $(q_{k_1})_{\text{\rm max}}$ is a subpath of $[f^{k_1}((q_0)_{\rm max})]$, hence
(using Lemma~\ref{lem 6.4})
$$l((q_{k_1})_{\text{\rm max}})\leqslant l([f^{k_1}(a_1')])
\leqslant ||f||^{k_1}\lceil l(a_1')\rceil\leqslant ||f||^{k_1}n_{\text{\rm critical}},$$
and the statement 2) in this case is valid.

Now we consider the case where $b_t'$ exists. The assumption of Claim~3 imply that $b_s$ cannot be virtually moved into $I_{k_1}$ in $k_1$ steps. By Claim~2, the path $b_t'$  can be virtually moved into $I_{k_1}$ in $k_1$ steps.
Then the initial vertex of $[f^{k_1}(a'_{t+1})]$  lies in~$I_{k_1}$.
Observe that $[f^{k_1}(a'_{t+1})]$ is a terminal subpath of $[f^{k_1}((q_0)_{\text{\rm max}})]$.
Then, by Lemma~\ref{lem 6.8},  $(q_{k_1})_{\text{\rm max}}$ is a subpath of $[f^{k_1}(a'_{t+1})]$, and we complete as above.
\hfill $\Box$

\medskip

{\tt Claim 4.} Each of the paths $b_s$, $b_t'$ can be virtually moved  into $I_{k_2}$ in $k_2$ steps.
\medskip

{\it Proof.} We assume that $b_s$ exists and prove the claim for $b_s$.
If $b_s$ can be  virtually moved into $I_{k_1+i}$ in $k_1+i$ steps for some $0\leqslant i\leqslant N_0N_1$, we are done.
So, suppose that this is not valid. Then we apply Claim~3 consequently to $(p_{i})_{\text{\rm max}}\equiv X_{i}[f^{i}(b_s)]Y_{i}$ for $i=0,1,\dots , N_0N_1$. The inequalities in Claim~3 and  the choice of $N_0,N_1$ imply that there must be a repetition in the sequence of pairs $(X_{k_1+i},(q_{k_1+i})_{\text{\rm e-max}})$, $i=0,1,\dots, N_0N_1$.
Then $b_s$ cannot be  virtually moved into $I_j$ in $j$ steps for $j=1,2,\dots$.
This contradicts Remark~\ref{Jahr_2014}\!~(b).
\hfill $\Box$

\medskip

{\tt Claim 5.} 
The paths $(p_{k_2})_{\tiny {\text {\rm max}}}$ and
$(q_{k_2})_{\tiny {\text {\rm max}}}$ coincide with some terminal subpaths of $[f^{k_2}(a_{s+1})]$ and $[f^{k_2}(a'_{t+1})]$, respectively.
\medskip

{\it Proof.} 
By Claim 4, the initial vertices of the paths $[f^{k_2}(a_{s+1})]$ and $[f^{k_2}(a'_{t+1})]$ lie in $I_{k_2}$, and their terminal points coincide with the terminal
points of the paths $[f^{k_2}((p_0)_{\text{\rm max}})]$ and $[f^{k_2}((q_0)_{\text{\rm max}})]$, respectively.
Then $\Lambda([f^{k_2}(a_{s+1})],[f^{k_2}(a'_{t+1})])\equiv\Lambda([f^{k_2}((p_0)_{\text{\rm max}})], [f^{k_2}((q_0)_{\text{\rm max}})])$,
and the claim follows from Lemma~\ref{lem 6.8}.
\hfill $\Box$


\medskip

{\tt Claim 6.} The chosen constant $T$ satisfies the proposition.

\medskip

{\it Proof.} 
Since $a_{s+1},a_{t+1}'\subset H_r$ are subpaths of $(p_0)_{\tiny {\text {\rm max}}}$ and
$(q_{0})_{\tiny {\text {\rm max}}}$, respectively, we have $l(a_{s+1})\leqslant n_{\text{\rm critical}}$ and
$l(a_{t+1}')\leqslant n_{\text{\rm critical}}$ by Lemma~\ref{lem 6.4}.
From this and from Claim~5, we deduce
$$l((p_{k_2})_{\tiny {\text {\rm max}}})\leqslant ||f||^{k_2}n_{\text{\rm critical}},\hspace*{10mm}
  l((q_{k_2})_{\tiny {\text {\rm max}}})\leqslant ||f||^{k_2}n_{\text{\rm critical}}.$$

\noindent
For any $0\leqslant i\leqslant N_2^2$, if we start from
$[f^i(\tau)]\equiv \bar{p}_iq_i$ instead of $\tau\equiv\bar{p}_0q_0$, we obtain analogously
$$l((p_{i+k_2})_{\tiny {\text {\rm max}}})\leqslant ||f||^{k_2}n_{\text{\rm critical}},\hspace*{10mm}
  l((q_{i+k_2})_{\tiny {\text {\rm max}}})\leqslant ||f||^{k_2}n_{\text{\rm critical}}.$$

This and the definition of $N_2$ imply that there exist $0\leqslant i<j\leqslant N_2^2$ such
that $(p_{i+k_2})_{\text{\rm e-max}}\equiv (p_{j+k_2})_{\text{\rm e-max}}$ and
$(q_{i+k_2})_{\text{\rm e-max}}\equiv (q_{j+k_2})_{\text{\rm e-max}}$.
Then $y$ is a non-deletable $r$-cancelation point in~$\tau$.
\hfill $\Box$


\hfill $\Box$


Now we analyze cancelations in $f$-images of paths in $G_r$ with several $r$-cancelation points.

\begin{defn}\label{defn 7.1}
{\rm Let $H_r$ be an exponential stratum in $\Gamma$. For any reduced path $\tau\subset G_r$, let $P_r(\tau)$ be the number of $r$-cancelation points in $\tau$.
Clearly, $P_r([f^{i}(\tau)])\geqslant P_r([f^{i+1}(\tau)])$ for any $i\geqslant 0$.
We say that $\tau$ is $r$-{\it stable} if $P_r([f^i(\tau)])= P_r(\tau)$ for all $i\geqslant 1$.
If $\tau$ is $r$-stable, we call the $r$-cancelation points of $\tau$ {\it non-deletable}.
Clearly, there exists a nonnegative integer $i_0$ such that the path $[f^{i_0}(\tau)]$ is $r$-stable.
We denote by $N_r(\tau)$ the number of $r$-cancelation points in $[f^{i_0}(\tau)]$.}
\end{defn}


\begin{thm}\label{prop 7.2} Let $H_r$ be an exponential stratum in $\Gamma$. There exists an efficient algorithm which,
given a reduced edge path $\tau\subset G_r$, computes $i_0\geqslant 0$ such
that the path $[f^{i_0}(\tau)]$ is $r$-stable. In particular, one can check, whether $\tau$ is $r$-stable or not.
\end{thm}

{\it Proof.} Let $y_1,\dots,y_k$ be all $r$-cance\-lation points in $\tau$.
Let $y_0$ and $y_{k+1}$ be the initial and the terminal points of $\tau$, respectively.
Let $\tau_i$ be the subpath of $\tau$ from $y_i$ to $y_{i+1}$, $i=0,\dots ,k$.

First we note that $y_i$ is the unique $r$-cancelation point in the edge path $\tau_{i-1}\tau_i$.
For each $1\leqslant i\leqslant k$, we can check in $T$ steps, by Proposition~\ref{prop 6.5}, whether the point $y_i$ in $\tau_{i-1}\tau_i$ is deletable or not. If it is deletable for some $i$, then $P_r([f^T(\tau)])<P_r(\tau)$ and
we can proceed by induction restarting from $[f^T(\tau)]$.
Suppose that $y_i$ is a non-deletable $r$-cancelation point in $\tau_{i-1}\tau_i$ for each $1\leqslant i\leqslant k$. We compute the $r$-cancelation area $A_i\equiv \bar{p}_i q_i$ in $\tau_{i-1}\tau_i$ and the number  $a_i:=L_r(p_i)=L_r(q_i)$.

If $L_r(\tau_i)\geqslant a_i+a_{i+1}$ for $i=1,\dots,k-1$,
then the path $\tau$ is $r$-stable and the points $y_1,\dots,y_k$ are non-deletable.

Suppose that for some $i$ we have $L_r(\tau_i)< a_i+a_{i+1}$. Then there exists a point $x\in \tau_i$ in the middle
between $\alpha(A_{i+1})$ and $\omega(A_i)$ with respect to $L_r$, i.e.
$L_r(\alpha(A_{i+1}),x)=L_r(x, \omega(A_i))=a$, where $a:=\frac{a_i+a_{i+1}-L_r(\tau_i)}{2}$. The points $\alpha(A_{i+1})$, $x$ and $\omega(A_i)$ divide $\tau_i$ into four subpaths: $\tau_i\equiv u_1u_2u_3u_4$.
Let $m$ be the minimal natural number such that
$\lambda_r^m\cdot a>L_{\text{\rm critical}}$.

Since $q_i\equiv u_1u_2u_3$ and $L_r(u_3)=a$,
the path
$[f^m(\tau_{i-1}u_1u_2)]$ is $r$-legal. Analogously, $[f^m(u_3u_4\tau_{i+1})]$ is $r$-legal.
Then  $[f^m(\tau_{i-1}\tau_i\tau_{i+1})]$ is the product of two $r$-legal paths, hence $P_r([f^m(\tau)])<P_r(\tau)$.
Induction by $k$ completes the proof.
\hfill $\Box$

\medskip


For convenience, we reformulate a part of Definition~\ref{defn 7.1}.

\begin{defn}\label{defn 7.3}
{\rm Let $H_r$ be an exponential stratum in $\Gamma$. Let $\tau$ be a reduced path in $G_r$  and $y_1,\dots,y_k$ be all $r$-cance\-lation points in $\tau$.
We say that these $r$-cancelation points in $\tau$ are {\it non-deletable} if the number of
$r$-cancelation points in $[f^i(\tau)]$ is equal to $k$ for every $i\geqslant 0$.}
\end{defn}

\begin{thm}\label{prop 7.4} (Criterium of $r$-stability) Let $H_r$ be an exponential stratum in $\Gamma$.
Let $\tau$ be a reduced edge path in $G_r$ and $y_1,\dots,y_k$ be all $r$-cancelation points in $\tau$.
Let $y_0:=\alpha(\tau)$, $y_{k+1}:=\omega(\tau)$, and
let $\tau_i$ be the subpath of $\tau$ from $y_i$ to $y_{i+1}$. Then the $r$-cancelation points
$y_1,\dots,y_k$ are non-deletable in $\tau$ if and only if each of $y_i$ is non-deletable in $\tau_{i-1}\tau_i$ for $i=1,\dots ,k$
and $L_r(\tau_i)\geqslant a_i+a_{i+1}$ for $i=1,\dots,k-1$,
where $a_i$ is the $r$-cancelation radius for the $r$-cancelation area in $\tau_{i-1}\tau_i$.
\end{thm}

The proof follows from the proof of Theorem~\ref{prop 7.2}.


\begin{defn}\label{defn 7.6}
{\rm Let $H_r$ be an exponential stratum in $\Gamma$.
Let $\tau$ be a reduced $r$-stable path in $G_r$.
Then $\tau$ can be written as $\tau\equiv b_0\cdot A_1\cdot b_1\cdot \ldots \cdot A_k\cdot b_k$, where
$b_0,\dots,b_k$ are $r$-legal or trivial paths in $G_r$ and $A_1,\dots ,A_k$ are all $r$-cancelation areas in $\tau$.
We call such decomposition the {\it $A$-decomposition} of~$\tau$.
}
\end{defn}

\begin{rmk}\label{rmk 7.7} Suppose that $\tau$ has the $A$-decomposition
$$\tau\equiv b_0\cdot A_1\cdot b_1\cdot \ldots \cdot A_k\cdot b_k.$$
Then, for every $i\geqslant 1$, the path $[f^i(\tau)]$ has the $A$-decomposition
$$[f^i(\tau)]\equiv b_0^i\cdot A_1^{\hspace*{1mm}i}\cdot b_1^i\cdot \ldots \cdot A_k^{\hspace*{1mm}i}\cdot b_k^i,$$
where $b_j^i\equiv [f^i(b_j)]$ and $A_j^{\hspace*{1mm}i}\equiv [f^i(A_j)]$ for all possible $j$.
\end{rmk}

\section{Splitting lemma}



\begin{notation}\label{endpoints_notations} {\rm Let $f:\Gamma\rightarrow \Gamma$ be a PL-relative train track.
For any exponential stratum $H_r$, let $\mathcal{E}_r$ be the set of endpoints of all $r$-cancelation areas.
}
\end{notation}




\begin{lem}\label{split_2} Let $f:\Gamma\rightarrow \Gamma$ be a PL-relative train track and $H_r$ be an exponential stratum.

\begin{enumerate}
\item[\rm 1)] The set $\mathcal{E}_r$ lies in $H_r$, is $f$-invariant, computable, and $|\mathcal{E}_r|\leqslant M_r$.
\item[\rm 2)] If $b$ is a nontrivial path in $G_r$ with endpoints in $\mathcal{E}_r$ and with $L_r(b)=0$, then $f^{M_r}(b)$ is an edge path in $G_{r-1}$.
\end{enumerate}
\end{lem}

{\it Proof.}
1) follows from Lemma~\ref{properties_areas}, Theorem~\ref{find_areas} and Proposition~\ref{number M}.2).
We prove~2).
By Lemma~\ref{zero_length}.5),
there exists $k$ such that $f^k(b)$ lies in $G_{r-1}$. The endpoints of $b$ and hence of $f^k(b)$ lie in $\mathcal{E}_r$, and $\mathcal{E}_r$ lies in $H_r$. Therefore the endpoints of $f^k(b)$ are vertices.
Since $\mathcal{E}_r$ is $f$-invariant and $|\mathcal{E}_r|\leqslant M_r$, and since the $f$-images of vertices are vertices, we have that the endpoints of $f^{M_r}(b)$ are vertices. Since $L_r(b)=0$, the path $f^{M_r}(b)$ lies in $G_{r-1}$.
\hfill $\Box$

\begin{lem}\label{Splitting_lemma} {\rm (Splitting lemma)}
For any PL-relative train track $f:\Gamma\rightarrow \Gamma$,
the following is satisfied:

Let $H_r$ be an exponential stratum of $\Gamma$ and
let $\tau$ be a reduced edge path in $G_r$. Then, for all $L>0$,  one can efficiently find an exponent $S>0$ such that at least one of the three possibilities occurs:

\begin{enumerate}
\item[{\rm 1)}] $[f^S(\tau)]$ contains an $r$-legal subpath of $r$-length greater than $L$.

\item[{\rm 2)}] $[f^S(\tau)]$ contains fewer illegal $r$-turns than $\tau$.

\item[{\rm 3)}] $[f^S(\tau)]$ is a trivial path or a concatenation of paths each of which is
either an indivisible periodic Nielsen path intersecting
${\text{\rm int}}(H_r)$ or an edge path in~$G_{r-1}$.
\end{enumerate}
\end{lem}

{\it Proof.}
By Theorem~\ref{prop 7.2}, we can
efficiently find the minimal nonnegative integer $i_0$ such that $[f^{i_0}(\tau)]$ is $r$-stable.
If $i_0>0$, then $\tau$ is not $r$-stable and we have~2) with $S:=i_0$.

Suppose that $i_0=0$. We may assume that $\tau$ is nontrivial. Then $\tau$ is $r$-stable and there exists
the $A$-decomposition
$\tau\equiv b_0\cdot A_1\cdot b_1\cdot \ldots \cdot A_k\cdot b_k$ as in Definition~\ref{defn 7.6}.
(Note that we can recognize all non-deletable $r$-cancelation points in $\tau$ and compute the $r$-cancelation radii
for all $A_j$; hence, we can compute all $L_r(b_j)$.)
Then, for all $i\geqslant 0$, we have the $A$-decomposition
$$[f^i(\tau)]\equiv b_0^i\cdot A_1^{\hspace*{1mm}i}\cdot b_1^i\cdot \ldots \cdot A_k^{\hspace*{1mm}i}\cdot b_k^i$$
with notations from Remark~\ref{rmk 7.7}.
Since all $b_j$ are $r$-legal or trivial, we have $$L_r(b_j^i)=\lambda_r^i\cdot L_r(b_j).$$

If $L_r(b_j)>0$ for some $j$, we compute the minimal natural $S$ with $\lambda_r^S\cdot L_r(b_j)>L$. Then we have 1) for this $S$.

Suppose that $L_r(b_j)=0$ for all $j$. We assume that $k\geqslant 1$, otherwise $\tau\equiv b_0$ is an edge path which lies in $G_{r-1}$ and we are done.

First we prove that $b_k$ is trivial or is an edge path in $G_{r-1}$.
Since $\omega(b_k)=\omega(\tau)$ is a vertex, it suffices to prove that $\alpha(b_k)$ is a vertex too.
Note that $\alpha(b_k)$ lies in an $r$-edge $E$; this follows from $\alpha(b_k)=\omega(A_k)$ and from the fact that the endpoints of $r$-cancelation areas lie in~$H_r$. Suppose that  $\alpha(b_k)$ is not a vertex. Then the first partial edge of $b_k$ is a non\-trivial terminal segment of $E$, hence $L_r(b_k)>0$, a contradiction.

Thus, $b_k$ is a (possibly trivial) edge path with $L_r(b_k)=0$. Then $b_k$ is trivial or is an edge path in $G_{r-1}$. Analogously $b_0$ is trivial or is an edge path in $G_{r-1}$.

Now consider $b_j$ with $j\in \{1,\dots ,k-1\}$. Suppose that $b_j$ is nontrivial. Since $L_r(b_j)=0$,
Lemma~\ref{split_2}.2) implies that $[f^{M_r}(b_j)]$ is an edge path in $G_{r-1}$. Clearly, each $r$-cancelation area $[f^{M_r}(A_j)]$ is an indivisible periodic Nielsen path.
Thus, we have the statement 3) with $S=M_r$ in this case.
\hfill $\Box$


\section{Subdivided relative train track}

It is technically convenient to start the proof of the main theorem in the situation where the following
condition is satisfied:

\medskip

(RTT-iv) There is a computable natural number $P=P(f)$ such that for each exponential strata $H_r$ and each $r$-cancelation area $A$ of $f$, the $r$-cancelation area $[f^P(A)]$ is an edge path.

\medskip

We will show that this condition is satisfied for a map $f':\Gamma'\rightarrow \Gamma'$
obtained from $f:\Gamma\rightarrow \Gamma$ by using subdivisions of edges of $\Gamma$
at exceptional points (see Definition~\ref{exceptional_points}). Moreover, $f'$
will satisfy the properties (RTT-i) -- (RTT-iii) of a relative train track with
respect to a natural filtration on $\Gamma'$.

\begin{defn}\label{exceptional_points} {\rm Let $f:\Gamma\rightarrow \Gamma$ be a PL-relative train track and $H_r$ be an expo\-nential stratum. Recall that $\mathcal{E}_r$ is the set of endpoints of all $r$-cancelation areas.\break
A point $v\in \mathcal{E}_r$ is called an {\it $r$-exceptional} point if $f^n(v)$
is not a vertex for all $n\in \mathbb{N}$.
}
\end{defn}

\begin{lem}\label{prop 8.2}  Let $f:\Gamma\rightarrow \Gamma$ be a PL-relative train track and $H_r$ be an exponential stratum.
Then the following statements are valid:

\begin{itemize}

\item[{\rm 1)}] The set of $r$-exceptional points is $f$-invariant. The $r$-exceptional points lie in the interiors of $r$-edges.

\item[{\rm 2)}] An endpoint $v$ of an $r$-cancelation area is not an $r$-exceptional point
                if and only if $f^i(v)$ is a vertex of $\Gamma$ for some $0\leqslant i<|\mathcal{E}_r|$.

\item[{\rm 3)}] One can efficiently find the set of $r$-exceptional points.

\end{itemize}
\end{lem}

{\it Proof.} 2) is obvious, 1) and 3) follow from Lemma~\ref{split_2}.1). \hfill $\Box$

\begin{defn}\label{defn 8.3} {\rm Let $f:\Gamma\rightarrow \Gamma$ be a PL-relative train track.
Let $\Gamma'$ be the graph obtained from $\Gamma$ by subdivision at all $r$-exceptional points
for all exponential strata $H_r$.
If $E\equiv e_1e_2\dots e_k$, $k\geqslant 2$, is the result of subdivision of an $r$-edge $E$, we call each $e_i$
an {\it $r$-exceptional partial edge}.
By  the statement 3) of Lemma~\ref{prop 8.2},
we can construct the graph $\Gamma'$ efficiently. By the statement 1), the map  $f:\Gamma\rightarrow \Gamma$
induces a natural map $f':\Gamma'\rightarrow \Gamma'$.

We also define a natural filtration $\varnothing =G'_0\subset \dots \subset G'_{N}=\Gamma'$:
for each exponential stratum $H_r$ the corresponding stratum $H_r'$ consists
of all $r$-edges which were not subdivided and of all $r$-exceptional partial edges;  all other strata remain unchanged.
}
\end{defn}

\begin{lem}
{\rm 1)} If $e$ is an $r$-exceptional partial edge in $\Gamma$, then $L_r(e)>0$.

{\rm 2)} The filtration $\varnothing =G'_0\subset \dots \subset G'_{N}=\Gamma'$ is maximal.

{\rm 3)} The homotopy equivalence $f':\Gamma'\rightarrow \Gamma'$ is a PL-relative train track.
\end{lem}

{\it Proof.}
1) If we had $L_r(e)=0$, then by Lemma~\ref{split_2}.\!~2), both endpoints of $f^{M_r}(e)$ were vertices that contradicts the fact that at least one endpoint of $e$ is $r$-exceptional.

2) Let $E'_1,E_2'$ be two edges in $H_r'$, where $H_r$ is an exponential stratum.
It suffices to show that $(f')^k(E_1')$ contains $E_2'$ for some $k\in \mathbb{N}$.

We consider $E'_1$ and $E_2'$ as (partial) $r$-edges in $\Gamma$.
By 1), we have $L_r(E'_1)>0$. By definition of $L_r$,
there exists $n\in \mathbb{N}$ such that $f^n(E'_1)$ contains a full $r$-edge $E$ in $\Gamma$.
Since the filtration for $f$ is maximal,
there exists $m\in \mathbb{N}$ such that $f^m(E)$ contains $E_2'$, and we are done.

3) If $H_r$ is a nonexponential stratum in $\Gamma$, then the transition matrices for $H_r$ and $H_r'$ coincide.
If $H_r$ is an exponential stratum in $\Gamma$, then the Perron-Frobenius eigenvalues for $H_r$ and $H_r'$ coincide.
Thus, $H_r$ is exponential, polynomial, or zero stratum if and only if $H_r'$ is exponential, polynomial, or zero stratum, respectively.

Propertiy (RTT-i) for $f'$ follows from the same property for $f$ with the help of Lemma~\ref{prop 8.2}.1). Properties (RTT-ii) and (RTT-iii) for $f'$ obviously follow from the corresponding properties for $f$.
 \hfill $\Box$

\begin{defn} {\rm The PL-relative train track $f':\Gamma\rightarrow \Gamma'$ constructed from\break $f:\Gamma\rightarrow \Gamma$ as above is called the {\it subdivided PL-relative train track} associated with~$f$.}
\end{defn}








\begin{rmk}\label{PPP}
{\rm After appropriate identification of $\Gamma$ and $\Gamma'$, the length functions $L_r$ and $L_r'$ coincide and the $r$-cancelation areas of $f$ and $f'$ coincide.
By Lemma~\ref{prop 8.2}.2), the property (RTT-iv)  is satisfied for $f'$ with
$$P:=\max\{|\mathcal{E}_r|\,|\, H_r\hspace*{2mm} {\text{\rm is an exponential stratum of}}\hspace*{2mm} \Gamma\}.$$
This number is computable by Lemma~\ref{split_2}.1).
}
\end{rmk}

\noindent
{\bf Agreement.} From now on we will work with $f'$ and never with the old $f$. So, we will simplify notation and skip dashes by writing $f, \Gamma, G_i, H_i$ instead of $f', \Gamma', G'_i, H_i'$, respectively.

\medskip

Below we collect the properties of the new $f$ which we will use later.


\begin{prop}\label{prop 8.5} The subdivided PL-relative train track $f:\Gamma\rightarrow \Gamma$ satisfies the following properties:

{\rm 1)} This $f$ represents the same automorphism of $F$ as the original one.

{\rm 2)} For each exponential stratum $H_r$,
there exists only finitely many $r$-cancelation areas. These areas can be computed.

{\rm 3)} Each $r$-cancelation area  has an initial and a terminal subpaths which lie in~$H_r$ and have nonzero $r$-lengths.

{\rm 4)} If $A$ is an $r$-cancelation area in a reduced $r$-stable path $\tau\subset G_r$, then $[f(A)]$ is an $r$-cancelation area in the reduced $r$-stable path $[f(\tau)]$.

{\rm 5)} One can compute a natural number $R_{\star}=R_{\star}(f)$ such that for each $r$-cancelation area $A=\bar{p}q$
where $p$ and $q$ are $r$-legal, we have $l(p)\leqslant R_{\star}$ and $l(q)\leqslant R_{\star}$.

{\rm 6)} One can compute a natural number $P=P(f)$ such that for every exponential stratum $H_r$
and every $r$-cancelation area $A$, the $r$-cancelation area $[f^P(A)]$ is an edge path.


\end{prop}

{\it Proof.} The statement 1) is obvious. The statements~2), 3), and 4) follow from Theorem~\ref{find_areas}, Lemma~\ref{properties_areas}, and Lemma~\ref{zero_length}.3). The statement~5) follows from 2).
The statement 6) is contained in Remark~\ref{PPP}. \hfill $\Box$

\section{$r$-superstable, $r$-perfect, and $A$-perfect paths}



in Sections~16 and~18 we need a stronger version of $r$-stability of paths, which we call $r$-superstability.

\medskip

\begin{defn}\label{defn 12.1} {\rm Let $H_r$ be an exponential stratum.
A reduced $f$-path $\tau\subset G_r$ is called {\it $r$-superstable}
if all $r$-cancelation points in $\tau$ and in $[\tau f(\tau)]$ are non-deletable
and all $r$-cancelation areas in these paths are edge paths.}
\end{defn}

\medskip

Note that if $\tau$ is $r$-superstable, then $[f^{i}(\tau)]$ is $r$-superstable for all $i\geqslant 0$.

\medskip

\begin{lem}\label{lem 12.2} Let $H_r$ be an exponential stratum. For any reduced $f$-path
$\tau\subset G_r$, one can efficiently compute a natural number $S=S(\tau)$  such that the path $[f^S(\tau)]$ is $r$-superstable.
\end{lem}

\medskip

{\it Proof.} By Theorem~\ref{prop 7.2} and Proposition~\ref{prop 8.5}.6),
we can compute a number $S_1$ such that
all $r$-cancelation points in $[f^{S_1}(\tau)]$ are non-deletable and all $r$-cancelation areas in $[f^{S_1}(\tau)]$ are edge paths.
Also, we can compute a number $S_2$ such that
all $r$-cancelation points in $[f^{S_2}([\tau f(\tau)])]$ are non-deletable and all $r$-cancelation areas in $[f^{S_2}([\tau f(\tau)])]$ are edge paths.
We set $S=\max\{S_1,S_2\}$.
Then the path $[f^S(\tau)]$ is $r$-superstable. \hfill $\Box$



\medskip

\begin{defn}\label{defn 9.1}
{\rm Let $H_r$ be an exponential stratum. An edge path $\tau\subset G_r$ is called {\it $r$-perfect} if
the following conditions are satisfied:
\begin{enumerate}
\item[(i)] $\tau$ is a reduced $f$-path and its first edge belongs to $H_r$,

\item[(ii)] $\tau$ is $r$-legal,

\item[(iii)] $[\tau f(\tau)]\equiv \tau\cdot [f(\tau)]$ and the turn of this path at the point between $\tau$
 and~$[f(\tau)]$ is legal.
\end{enumerate}

\noindent
A vertex in $D_f$ is called {\it $r$-perfect} if the corresponding $f$-path in $\Gamma$ is
$r$-perfect.
}
\end{defn}

Note that these conditions imply that $[\tau f(\tau)]$ is $r$-legal.
In the following proposition we formulate some important properties of $r$-perfect paths; they
can be proved directly from the above definition.

\medskip

\begin{prop}\label{prop 9.2} Let $H_r$ be an exponential stratum and let
$\tau$ be an $r$-perfect path in $G_r$. Then the following statements are satisfied:

\begin{enumerate}
\item[(1)] For every $i\geqslant 0$, the path $[f^i(\tau)]$ is $r$-perfect.

\medskip

\item[(2)] For every $i\geqslant 0$, the vertex $[f^i(\tau)]$ of $D_f$ lies in the $\tau$-subgraph.
Moreover, $[f^i(\tau)]\equiv \widehat{f}^{\hspace*{1.5mm}m_i}(\tau)$ for some computable $m_i$ satisfying $m_0=0$, $m_i<m_{i+1}$.

\medskip

\item[(3)] For every $i\geqslant 0$, the path $\widehat{f}^{\hspace*{1.5mm}i}(\tau)\subset G_r$ is $r$-legal, contains edges from $H_r$,
and $L_r(\widehat{f}^{\hspace*{1.5mm}i+1}(\tau))\geqslant L_r(\widehat{f}^{\hspace*{1.5mm}i}(\tau))$.
\end{enumerate}

\end{prop}

\medskip

\begin{defn}\label{defn 9.3}
{\rm
Let $H_r$ be an exponential stratum. A reduced $f$-path $\tau \subset G_r$
containing edges from $H_r$ is called {\it $A$-perfect} if

\begin{enumerate}
\item[(i)] all $r$-cancelation points in $\tau$ are non-deletable,\\ the corresponding $r$-cancelation areas are edge paths,

\medskip

\item[(ii)] the $A$-decomposition of $\tau$ begins with an $A$-area,
i.e. it has the form\\ $\tau\equiv A_1b_1\dots A_kb_k$, $k\geqslant 1$,

\medskip

\item[(iii)] $[\tau f(\tau)]\equiv \tau\cdot [f(\tau)]$ and the turn at the point between $\tau$ and $[f(\tau)]$ is legal.
\end{enumerate}

\noindent
A vertex in $D_f$ is called {\it $A$-perfect} if the corresponding $f$-path in $\Gamma$ is
$A$-perfect.}
\end{defn}

Note that the first edge of such $\tau$ lies in $H_r$.
The following proposition can be proved straightforward and we leave it for the reader.

\begin{prop}\label{prop 9.4}
Let $H_r$ be an exponential stratum and let $\tau$ be an $A$-perfect path
in $G_r$ with the $A$-decomposition $\tau\equiv A_1b_1\dots A_kb_k$.

For $1\leqslant j\leqslant k$, we set $\tau_{0,j}\equiv [A_jb_j\dots A_kb_kf(A_1b_1\dots A_{j-1}b_{j-1})]$ and for $i\geqslant 1$ we set
$\tau_{i,j}\equiv[f^i(\tau_{0,j})]$.
Then the following statements are satisfied:

\begin{enumerate}
\item[(1)] For any $1\leqslant j\leqslant k$ and $i\geqslant 0$ the path $\tau_{i,j}$ is $A$-perfect.

\medskip

\item[(2)] For any $1\leqslant j\leqslant k$ and $i\geqslant 0$ the vertex $\tau_{i,j}$ of $D_f$
lies in the $\tau$-subgraph. Moreover $\tau_{i,j}\equiv\widehat{f}^{\hspace*{1.5mm}m_{i,j}}(\tau)$
for some computable $m_{i,j}$ satisfying $m_{0,1}=0$, $m_{i,j}<m_{i,j+1}$, and $m_{i,k}<m_{i+1,1}$.

\medskip

\item[(3)] All $A$-perfect vertices of the $\tau$-subgraph are $\tau_{i,j}$, $1\leqslant j\leqslant k$, $i\geqslant 0$.

\medskip

\item[(4)]
For every vertex $\sigma$ in the $\tau$-subgraph, at least one of the paths $\sigma$, $\widehat{f}(\sigma),\dots ,\widehat{f}^{\hspace*{1.5mm}l(\sigma)}(\sigma)$ coincides with $\tau_{i,j}$ for some $i,j$.
\medskip

\end{enumerate}
\end{prop}


\section{$\mu$-subgraphs in the case where $\mu$ and $[\mu f(\mu)]$ are $r$-legal}

\begin{prop}\label{prop 10.1}
Let $H_r$ be an exponential stratum. Let $\mu$ be a reduced $r$-legal $f$-path
with endpoints in $H_r$ and with the first edge in $G_{r-1}$, and the last edge in~$H_r$.
Then the $\mu$-subgraph of $D_f$ contains a vertex which is an $r$-perfect path in $G_r$. This vertex can be efficiently found.
\end{prop}

{\it Proof.} We write $\mu\equiv \sigma\tau$, where $\sigma$ is a nontrivial path in $G_{r-1}$ and $\tau$ is
a nontrivial path in $G_r$ with the first and the last edges from $H_r$.
We prove that $\mu':=\widehat{f}^{\hspace*{1.5mm}l(\sigma)}(\mu)$ is $r$-perfect.


Clearly, $\mu'\equiv \tau[f(\sigma)]$, $\mu'$ is $r$-legal and lies in the $\mu$-subgraph.
It remains to prove that the turn between $\mu'$ and $[f(\mu')]$ in $\mu' \cdot [f(\mu')]$ is legal.

\vspace*{0mm}
\unitlength 1mm
\linethickness{0.4pt}
\ifx\plotpoint\undefined\newsavebox{\plotpoint}\fi 
\begin{picture}(60,39.5)(0,20)
\put(15,49.5){\circle*{1.00}}
\put(15,49.5){\line(1,0){18}}
\put(24,51.5){\makebox(0,0)[lc]{\small $\sigma$}}
\put(33,49.5){\thicklines\line(1,0){23}}
\put(44.5,51.5){\makebox(0,0)[lc]{\small$\tau$}}
\put(56,49.5){\circle*{1.00}}
\put(35.5,47.5){\makebox(0,0)[ct]{$\underbrace{\hspace{41mm}}$}}
\put(35.5,41.00){\makebox(0,0)[lc]{$\mu$}}
\put(56,49.5){\line(1,0){22}}
\put(67,51.75){\makebox(0,0)[cc]{\small $[f(\sigma)]$}}
\put(78,49.5){\thicklines\line(1,0){30}}
\put(93,51.75){\makebox(0,0)[cc]{\small $[f(\tau)]$}}
\put(108,49.5){\circle*{1.00}}
\put(82.5,47.5){\makebox(0,0)[ct]{$\underbrace{\hspace{52mm}}$}}
\put(82.5,41.00){\makebox(0,0)[lc]{$[f(\mu)]$}}
\put(108,49.5){\line(1,0){18}}
\put(117,51.75){\makebox(0,0)[cc]{\small $[f^2(\sigma)]$}}
\put(55.5,37){\makebox(0,0)[ct]{$\underbrace{\hspace{45mm}}$}}
\put(55.5,32.5){\makebox(0,0)[cc]{$\mu'$}}
\put(102,37){\makebox(0,0)[ct]{$\underbrace{\hspace{48mm}}$}}
\put(102,32.5){\makebox(0,0)[cc]{$[f(\mu')]$}}
\end{picture}

\vspace*{-10mm}
\begin{center}
Figure 9.
\end{center}

\bigskip

By (RTT-ii), $[f(\sigma)]$ is a nontrivial path in $G_{r-1}$. So,
the last edge of $\mu'$ lies in $G_{r-1}$. On the other hand,
the first edge of $[f(\mu')]\equiv [f(\tau)][f^2(\sigma)]$ lies in $H_r$.
Therefore, the turn under consideration is mixed and hence legal.\hfill$\Box$

\begin{prop}\label{prop 10.2}
Let $H_r$ be an exponential stratum and let $\mu\subset G_r$ be a reduced $f$-path
such that $\mu$ is $r$-legal and $[\mu f(\mu)]$ is $r$-legal or trivial.
Then the $\mu$-subgraph of $D_f$ contains a vertex $\mu'$
for which at least one of the following statements is satisfied:

\medskip

{\rm (1)} the corresponding path $\mu'$ lies in $G_{r-1}$;

\medskip

{\rm (2)} the corresponding path $\mu'$ is $r$-perfect;

\medskip

{\rm (3)} the $\mu'$-subgraph is finite.

\medskip
Moreover, there exists an efficient algorithm which constructs such a vertex $\mu'$ and indicates
one of the statements {\rm (1)-(3)} which is satisfied for $\mu'$.

\end{prop}


{\it Proof.} We proceed by induction on the number of $r$-edges in $\mu$.
We assume that $\mu$ contains an $r$-edge, otherwise (1) is satisfied for $\mu':=\mu$.
Write $\mu\equiv b_1\cdot b_2\cdot b_3$,
where the (possibly trivial) edge paths $b_1,b_3$ lie in $G_{r-1}$ and the path $b_2$ begins and ends with $r$-edges.
We claim that the path $\widetilde{\mu}:= b_2\cdot [b_3f(b_1)]$ satisfies the following properties:

\begin{enumerate}
\item[i)] The path $\widetilde{\mu}$ is a reduced $f$-path which begins with an $r$-edge;\\
it has the same number of $r$-edges as $\mu$.

\item[ii)] The vertex $\widetilde{\mu}$ of $D_f$ lies in the $\mu$-subgraph.

\item[iii)] The path $\widetilde{\mu}$ is $r$-legal and the path $[\widetilde{\mu}f(\widetilde{\mu})]$ is $r$-legal or trivial.

\end{enumerate}

Indeed, i) is obvious; ii) is valid since
$\widetilde{\mu}\equiv \widehat{f}^{\hspace*{1.5mm}l(b_1)}(\mu)$.


The claim iii) follows from the formulas $$\widetilde{\mu}\equiv[\bar{b}_1\mu f(b_1)],\hspace*{5mm} [\widetilde{\mu}f(\widetilde{\mu})]\equiv[\bar{b}_1[\mu f(\mu)]f^2(b_1)].$$

Indeed, the $r$-legality of $\widetilde{\mu}$ follows from the assumption that $\mu$ is $r$-legal
and from the fact that the product of an $r$-legal path with a path from $G_{r-1}$
followed by tightening is $r$-legal or trivial.
Note that $\widetilde{\mu}$ is nontrivial since $\mu$ contains an $r$-edge.
Analogously, the path $[\widetilde{\mu}f(\widetilde{\mu})]$ is $r$-legal or trivial, since, by assumption, $[\mu f(\mu)]$ is $r$-legal or trivial.



Thus, replacing $\mu$ by $\widetilde{\mu}$, we may additionally assume that
$\mu$ begins with an edge $E$ from $H_r$.
Since $\mu$ is $r$-legal, $f(E)$ is an initial subpath of $[f(\mu)]$.
Let $\tau$ be the initial subpath of $\mu$ such that
$\mu\equiv \tau \cdot \overline{I(\bar{\mu},[f(\mu)])}$.  We consider the following four cases
(see Figure~10).

\medskip

{\tt Case (1)} The path $\tau$ is nontrivial (hence, it contains $E$) and
$f(E)$ is an initial subpath of $I(\bar{\mu},[f(\mu)])$.


{\tt Case (2)} The path $\tau$ is nontrivial (hence, it contains $E$) and $I(\bar{\mu},[f(\mu)])$ is
a proper initial subpath of $f(E)$.


{\tt Case (3)} The path $\tau$ is trivial (hence, $\mu\equiv\overline{I(\bar{\mu},[f(\mu)])}$\,\,) and
$f(E)$ is a proper initial subpath of $I(\bar{\mu},[f(\mu)])$.


{\tt Case (4)} The path $\tau$ is trivial (hence, $\mu\equiv\overline{I(\bar{\mu},[f(\mu)])}$\,\,) and
$I(\bar{\mu},[f(\mu)])$ is an initial subpath of $f(E)$.

\newpage

\vspace*{4mm}
\hspace*{10mm}
{
\unitlength=1.00mm
\special{em:linewidth 0.4pt}
\linethickness{0.4pt}
\begin{picture}(105.00,46.00)
\put(5.00,25.00){\thicklines\line(1,0){8.00}}
\put(13.00,25.00){\line(1,0){10.00}}
\put(23.00,25.00){\line(0,1){20.00}}
\put(26.00,45.00){\thicklines\line(0,-1){8.00}}
\put(26.00,37.00){\line(0,-1){12.00}}
\put(26.00,25.00){\line(1,0){15.00}}
\put(27.50,37.00){\line(0,-1){2.50}}
\put(27.50,34.00){\line(0,-1){1.5}}
\put(27.50,32.00){\line(0,-1){1.5}}
\put(27.50,30.00){\line(0,-1){1.5}}
\put(27.50,28.00){\line(0,-1){1.5}}
\put(28,26.50){\line(1,0){1.5}}
\put(30,26.50){\line(1,0){1.5}}
\put(32,26.50){\line(1,0){1.5}}
\put(34,26.50){\line(1,0){1.5}}
\put(21.50,37.00){\line(0,-1){10.50}}
\put(21.50,26.50){\line(-1,0){8.50}}
\put(13.00,25.00){\circle*{1.00}}
\put(13.00,26.50){\circle*{1.00}}
\put(21.50,37.00){\circle*{1.00}}
\put(23.00,37.00){\circle*{1.00}}
\put(26.00,37.00){\circle*{1.00}}
\put(5.00,25.00){\circle*{1.00}}
\put(27.00,41.00){\makebox(0,0)[lc]{$f(E)$}}
\put(30.00,29.00){\makebox(0,0)[lb]{$[f(\mu_1)]$}}
\put(19.00,29.00){\makebox(0,0)[rb]{$\mu_1$}}
\put(9.00,24.00){\makebox(0,0)[ct]{$E$}}
\put(14.50,20.00){\makebox(0,0)[ct]{$\underbrace{\hspace{1.9cm}}$}}
\put(33.00,20.00){\makebox(0,0)[ct]{$\underbrace{\hspace{1.6cm}}$}}
\put(14.50,17.00){\makebox(0,0)[ct]{$\mu$}}
\put(33.00,17.00){\makebox(0,0)[ct]{$[f(\mu)]$}}
\put(24.00,4.00){\makebox(0,0)[cc]{(1)}}
\put(80.00,4.00){\makebox(0,0)[cc]{(2)}}
\put(56.00,25.00){\thicklines\line(1,0){5.00}}
\put(61.00,25.00){\line(1,0){13.00}}
\put(74.00,25.00){\line(0,1){20.00}}
\put(56.00,25.00){\circle*{1.00}}
\put(61.00,25.00){\circle*{1.00}}
\put(77.00,45.00){\thicklines\line(0,-1){20.00}}
\put(77.00,25.00){\thicklines\line(1,0){5.00}}
\put(82.00,25.00){\circle*{1.00}}
\put(82.00,25.00){\line(1,0){23.00}}
\put(61.00,23.50){\line(1,0){21.00}}
\put(61.00,23.50){\circle*{1.00}}
\put(82.00,23.50){\circle*{1.00}}
\put(58.00,26.00){\makebox(0,0)[cb]{$E$}}
\put(79.00,36.00){\makebox(0,0)[lc]{$f(E)$}}
\put(65.50,18.00){\makebox(0,0)[ct]{$\underbrace{\hspace{1.9cm}}$}}
\put(91.00,18.00){\makebox(0,0)[ct]{$\underbrace{\hspace{2.8cm}}$}}
\put(65.50,15.00){\makebox(0,0)[ct]{$\mu$}}
\put(91.00,15.00){\makebox(0,0)[ct]{$[f(\mu)]$}}
\put(71.00,22.00){\makebox(0,0)[ct]{$\mu_1$}}
\put(21.50,32.00){\vector(0,1){0.00}}
\put(27.50,30.00){\vector(0,-1){0.00}}
\put(72.00,23.50){\vector(1,0){0.00}}
\put(26.00,45.00){\circle*{1.00}}
\put(27.50,37.00){\circle*{1.00}}
\put(77.00,45.00){\circle*{1.00}}
\end{picture}
}

\vspace*{10mm}
\hspace*{22mm}
{
\unitlength=1.00mm
\special{em:linewidth 0.4pt}
\linethickness{0.4pt}
\begin{picture}(83.00,56.00)
\put(15.00,20.00){\line(0,1){24.00}}
\put(15.00,44.00){\thicklines\line(0,1){12.00}}
\put(15.00,56.00){\circle*{1.00}}

\put(18.00,44.00){\circle*{1.00}}
\put(18.00,44.00){\thicklines\line(0,-1){12.00}}
\put(18.00,38.00){\vector(0,1){0}}
\put(18.00,32.00){\circle*{1.00}}
\put(15.00,20.00){\line(1,0){20.00}}

\put(16.50,29.00){\vector(0,-1){0}}
\put(16.50,44.00){\line(0,-1){10.00}}
\put(16.50,33.5){\line(0,-1){1.5}}
\put(16.50,31.5){\line(0,-1){1.5}}
\put(16.50,29.5){\line(0,-1){1.5}}
\put(16.50,27.5){\line(0,-1){1.5}}
\put(16.50,25.5){\line(0,-1){1.5}}
\put(16.50,23.5){\line(0,-1){1.5}}
\put(17,22.00){\line(1,0){1.5}}
\put(19,22.00){\line(1,0){1.5}}
\put(21,22.00){\line(1,0){1.5}}
\put(23,22.00){\line(1,0){1.5}}
\put(25,22.00){\line(1,0){1.5}}
\put(11.00,20.00){\thicklines\line(0,1){12.00}}
\put(11.00,32.00){\line(0,1){24.00}}
\put(11.00,32.00){\circle*{1.00}}
\put(11.00,20.00){\circle*{1.00}}
\put(15.00,44.00){\circle*{1.00}}
\put(10.00,26.00){\makebox(0,0)[rc]{$E$}}
\put(18.00,51.00){\makebox(0,0)[lc]{$f(E)$}}

\put(19.00,38.00){\makebox(0,0)[lc]{$\mu_1$}}
\put(17.50,23.00){\makebox(0,0)[lb]{$[f(\mu_1)]$}}
\put(25.00,16.00){\makebox(0,0)[ct]{$\underbrace{\hspace{2cm}}$}}
\put(24.00,13.00){\makebox(0,0)[ct]{$[f(\mu)]$}}
\put(12.00,37.00){\makebox(0,0)[rc]{$\mu\ \begin{cases}\\ \\ \\ \\ \\ \\ \end{cases}$}}
\put(19.00,3.00){\makebox(0,0)[cc]{(3)}}
\put(57.00,20.00){\thicklines\line(0,1){12.00}}
\put(57.00,32.00){\line(0,1){24.00}}
\put(57.00,32.00){\circle*{1.00}}
\put(57.00,20.00){\circle*{1.00}}
\put(56.00,26.00){\makebox(0,0)[rc]{$E$}}
\put(64.00,38.00){\makebox(0,0)[lc]{$f(E)$}}
\put(72.00,16.00){\makebox(0,0)[ct]{$\underbrace{\hspace{2.3cm}}$}}
\put(70.00,13.00){\makebox(0,0)[ct]{$[f(\mu)]$}}
\put(58.00,37.00){\makebox(0,0)[rc]{$\mu\ \begin{cases}\\ \\ \\ \\ \\ \\ \end{cases}$}}
\put(65.00,3.00){\makebox(0,0)[cc]{(4)}}
\put(56.00,42.00){\makebox(0,0)[rc]{$\sigma$}}
\put(61.00,56.00){\circle*{1.00}}
\put(61.00,20.00){\thicklines\line(0,1){36.00}}
\put(61.00,20.00){\thicklines\line(1,0){12.00}}
\put(73.00,20.00){\line(1,0){10.00}}
\put(73.00,20.00){\circle*{1.00}}
\put(67.00,19.00){\makebox(0,0)[ct]{$\delta$}}
\end{picture}
}

\vspace*{-5mm}
\begin{center}
Figure 10.

\begin{minipage}[t]{120mm}
In (1) and (2), the path $\mu$ is depictured by the large left corner, and in (3), (4),
it is depictured by a vertical line. The path $[f(\mu)]$ is depictured by the large right corner.
\end{minipage}
\end{center}

\medskip

{\tt Consider Cases (1) and (3)}.

Let $\mu_1:\equiv [\overline{E}\mu f(E)]$. Clearly, $\mu_1$ lies in the
$\mu$-subgraph. Moreover, the paths $\mu_1$ and $[\mu_1f(\mu_1)]$ are $r$-legal or trivial.
Indeed:

in both cases $\mu_1$ is a subpath of the $r$-legal path $\mu$;

in Case (1), $[\mu_1f(\mu_1)]$ is a subpath of one of the $r$-legal paths $\mu$ or $[\mu f(\mu)]$;

in Case (3), $[\mu_1f(\mu_1)]$ is a subpath of one of the $r$-legal paths $\mu$ or $[f(\mu)]$.


We may assume that $\mu_1$ is nontrivial, otherwise the statement (3) is valid for $\mu':=\mu_1$.
Since the number of $r$-edges in $\mu_1$ is less than in $\mu$, we can proceed by induction and
find the desired $\mu'$ in the $\mu_1$-subgraph.

\medskip

{\tt Consider Case (2).} 
Let $E_1$ be the second edge of $[\mu f(\mu)]$.
In Figure~11, we distinguish the cases where $\tau$ contains at least two edges (picture (1))
and where $\tau$ contains exactly one edge (picture (2)).

We set $\mu_1:\equiv [\overline{E}\mu f(E)]$. The path $\mu_1$ is $r$-legal as a subpath of $[\mu f(\mu)]$,
it begins with $E_1$ and ends with the last edge of $f(E)$, and it lies in the $\mu$-subgraph.

Suppose that $E_1$ is an $r$-edge. Then $[f(\mu_1)]$ begins with the first edge of $f(E_1)$.
Therefore the turn between $\mu_1$ and $[f(\mu_1)]$ coincides with the turn between $f(E)$ and $f(E_1)$. This
turn is legal, since the turn between $E$ and $E_1$ is an $r$-turn in the $r$-legal path $[\mu f(\mu)]$.
Hence, $\mu_1$ is $r$-perfect and we have the statement~(2) for $\mu':\equiv\mu_1$.

Now suppose that $E_1$ lies in $G_{r-1}$.
By Proposition~\ref{prop 10.1} applied to $\mu_1$, we can construct $\mu'$ satisfying the statement~(2).

\medskip
\medskip

\hspace*{-9mm}
{
\unitlength=1.00mm
\special{em:linewidth 0.4pt}
\linethickness{0.4pt}
\begin{picture}(143.00,45.50)
\put(111.00,25.00){\thicklines\line(1,0){14.00}}
\put(125.00,25.00){\line(0,1){20.00}}
\
\put(128.00,45.00){\thicklines\line(0,-1){20.00}}
\put(128.00,25.00){\thicklines\line(1,0){8.00}}
\
\put(136.00,25.00){\line(1,0){7.00}}
\put(136.00,25.00){\circle*{1.00}}
\put(128.00,23.50){\line(1,0){8.00}}
\put(136.00,23.50){\circle*{1.00}}
\put(128.00,23.50){\circle*{1.00}}
\put(119.00,26.00){\makebox(0,0)[cb]{$E$}}
\put(130.00,33.00){\makebox(0,0)[lc]{$f(E)$}}
\put(132.00,22.00){\makebox(0,0)[ct]{$\mu_1$}}
\put(118.50,18.00){\makebox(0,0)[ct]{$\underbrace{\hspace{1.4cm}}$}}
\put(135.00,18.00){\makebox(0,0)[ct]{$\underbrace{\hspace{1.6cm}}$}}
\put(118.50,15.00){\makebox(0,0)[ct]{$\mu$}}
\put(135.00,15.00){\makebox(0,0)[ct]{$[f(\mu)]$}}
\put(126.00,4.00){\makebox(0,0)[cc]{(2)}}
\
\
\put(1.00,25.00){\thicklines\line(1,0){5.00}}
\put(6.00,25.00){\line(1,0){13.00}}
\put(19.00,25.00){\line(0,1){20.00}}
\put(1.00,25.00){\circle*{1.00}}
\put(6.00,25.00){\circle*{1.00}}
\
\put(22.00,45.00){\thicklines\line(0,-1){20.00}}
\put(22.00,25.00){\thicklines\line(1,0){5.00}}
\
\put(27.00,25.00){\circle*{1.00}}
\put(41.00,25.00){\circle*{1.00}}
\put(27.00,25.00){\line(1,0){23.00}}
\put(6.00,23.50){\line(1,0){21.00}}
\put(6.00,23.50){\circle*{1.00}}
\put(27.00,23.50){\circle*{1.00}}
\put(3.00,26.00){\makebox(0,0)[cb]{$E$}}
\put(9.00,26.00){\makebox(0,0)[cb]{$E_1$}}
\put(34.00,26.00){\makebox(0,0)[cb]{$f(E_1)$}}
\put(24.00,36.00){\makebox(0,0)[lc]{$f(E)$}}
\put(10.50,18.00){\makebox(0,0)[ct]{$\underbrace{\hspace{1.9cm}}$}}
\put(36.00,18.00){\makebox(0,0)[ct]{$\underbrace{\hspace{2.8cm}}$}}
\put(10.50,15.00){\makebox(0,0)[ct]{$\mu$}}
\put(36.00,15.00){\makebox(0,0)[ct]{$[f(\mu)]$}}
\put(16.00,22.00){\makebox(0,0)[ct]{$\mu_1$}}
\put(12.00,25.00){\circle*{1.00}}
\put(22.00,45.00){\circle*{1.00}}
\put(17.00,23.50){\vector(1,0){0.00}}
\put(128.00,45.00){\circle*{1.00}}
\put(111.00,25.00){\circle*{1.00}}
\put(125.00,25.00){\circle*{1.00}}
\put(133.00,23.50){\vector(1,0){0.00}}
\put(20.00,4.00){\makebox(0,0)[cc]{(1)}}

\put(132.00,25.00){\circle*{1.00}}
\put(130.50,25.50){\makebox(0,0)[cb]{$E_1$}}
\end{picture}
}

\begin{center} Figure 11.
\end{center}

\medskip

{\tt Consider Case (4).} Let $\mu\equiv E\cdot \sigma$, then $f(E)\equiv \bar{\sigma}\cdot \overline{E}\cdot \delta$
for some path $\delta$. This implies $\widehat{f}(\mu)\equiv \overline{E}\cdot \delta$ and $\widehat{f}^{\hspace*{1.5mm}2}(\mu)\equiv \mu$. Hence,
the $\mu$-subgraph is finite and for $\mu':\equiv \mu$ we have the statement (3).\hfill $\Box$


\section{$\mu$-subgraphs in the case where $\mu$ is $r$-legal, but $[\mu f(\mu)]$ is not}

Recall that $C_{\star}$ is the constant from Lemma~\ref{lem 3.12} and $R_{\star}$ is the constant from Proposition~\ref{prop 8.5}.5).

\begin{prop}\label{prop 10.3} Let $H_r$ be an exponential stratum.
Let $\mu\subset G_r$ be a reduced $f$-path which is $r$-legal.
Suppose that $[\mu f(\mu)]$ contains a non-deletable $r$-cancelation point
and the $r$-cancelation area in $[\mu f(\mu)]$ is an edge path. Then the $\mu$-subgraph contains a vertex $\mu'$ for which at least one of the following statements is satisfied:

\medskip

{\rm (1)} $l(\mu')\leqslant 3C_{\star}+2R_{\star}+||f||$;

\medskip

{\rm (2)} the path $\mu'$ in $\Gamma$ is $A$-perfect.

\medskip

Moreover, there exists $k\leqslant l(\mu)$ such that $\mu':= \widehat{f}^{\hspace*{1.5mm}k}(\mu)$ satisfies {\rm (1)} or {\rm (2)}.
In particular, such $\mu'$ can be efficiently found.


\end{prop}

{\it Proof.}
For $i=0,1$, let $z^i$ be the unique non-deletable $r$-cancelation point in the path  $[f^i(\mu)f^{i+1}(\mu)]$
and let $A(z^i)$ be the $r$-cancelation area in this path. We have $A(z^1)\equiv [f\bigl(A(z^0)\bigr)]$. Let $a$ and $b$ be the initial and the terminal vertices of $A(z^0)$. Then $f(a)$ and $f(b)$ are the initial and the terminal vertices of $A(z^1)$.

\medskip

Let $\mu\equiv E_1\dots E_n$, and  let $\mu_m\equiv E_1\dots E_m$ be the initial subpath of $\mu$ until the vertex $z^0$. Thus $I(\overline{\mu},[f(\mu)])\equiv \bar{E}_n\dots \bar{E}_{m+1}$.
Denote $\mu_i:=E_1\dots E_i$ for $i=1,\dots ,n$.

\medskip

{\tt Case 1.} Suppose that there exists $i\in \{1,\dots,m\}$ such that $E_i\dots E_n$ cancels out in $E_i\dots E_n\cdot [f(E_1\dots E_{i-1})]$ and that $i$ is minimal with this property.

Then, using induction, one can show that
$$\widehat{f}^{\hspace*{1.5mm} i-1}(\mu)\equiv [E_i\dots E_n\cdot [f(E_1\dots E_{i-1})]].$$
By definition of $i$, we have $[f(E_1\dots E_{i-1})]\equiv \bar{E}_n\dots \bar{E}_i\ell$ for some edge path $\ell$.
Since $i\leqslant m$, we can write $[f(\mu_{i-1})]\equiv \bar{E}_n\dots \bar{E}_{m+1} \bar{E}_m\dots \bar{E}_i\ell$. Hence, the maximal common initial subpath of $[f(\mu_{i-1})]$ and $[f(\mu)]$ is $\bar{E}_n\dots \bar{E}_{m+1}$.
By Lemma~\ref{Cooper_1}, we have $l(\bar{E}_m\dots \bar{E}_i\ell)\leqslant C_{\star}$, and hence $l(\ell)\leqslant C_{\star}$.
Since $\widehat{f}^{\hspace*{1.5mm} i-1}(\mu)\equiv~\!\ell$, we have the statement (1) for $\mu':=\widehat{f}^{\hspace*{1.5mm} i-1}(\mu)$.

\medskip

{\tt Case 2.} Suppose that cancelations in $E_i\dots E_n\cdot [f(E_1\dots E_{i-1})]$ don't meet $E_i$
for each $i=1,\dots,m$. Then, for these $i$, one can deduce by induction that
$$\widehat{f}^{\hspace*{1.5mm} i}(\mu)\equiv [E_{i+1}\dots E_n\cdot [f(E_1\dots E_i)]].\eqno{(15.1)}$$

{\it Notation.} For two vertices $u,v$ in $[f(\mu)]$, we write $u\prec v$ if $u$ and $v$ are the initial and the terminal
vertices of some nontrivial subpath of $[f(\mu)]$.

\medskip

{\it Case 2.1.} Let $z^0 \preccurlyeq f(z^0)$.

We write $\mu\equiv p\cdot q$, where $p\equiv E_1E_2\dots E_s$ is the initial subpath of $\mu$ until the vertex $a$, and $q$ is the terminal subpath of $\mu$ starting from the vertex $a$.
Since $s<m$, the equation (15.1) implies  $\widehat{f}^{\hspace*{1.5mm} s}(\mu)\equiv [qf(p)]$.
We show that $\mu':=\widehat{f}^{\hspace*{1.5mm} s}(\mu)$ satisfies the proposition.
Consider three variants for the position of $f(a)$ in $[f(\mu)]$.

\medskip

{\it Variant 1.} Suppose that $b\prec f(a)$ (see Figure 12).

\unitlength 1mm 
\linethickness{0.4pt}
\ifx\plotpoint\undefined\newsavebox{\plotpoint}\fi 
\begin{picture}(90,40)(-5,0)
\put(2,10){\line(0,1){20}}
\put(2,10){\line(1,0){40}}
\put(42,10){\line(0,1){20}}
\put(43,6.5){\makebox(0,0)[cc]{\small $z^0$}}
\put(33,10){\makebox(0,0)[cc]{\small $($}}
\put(33,6){\makebox(0,0)[cc]{\small $a$}}
\put(22,2){\makebox(0,0)[ct]{\large $\underbrace{\hspace{40mm}}_{\mu}$}}
\put(3,11.3){\line(1,0){2}}
\put(6,11.3){\line(1,0){2}}
\put(9,11.3){\line(1,0){2}}
\put(12,11.3){\line(1,0){2}}
\put(15,11.3){\line(1,0){2}}
\put(18,11.3){\line(1,0){2}}
\put(21,11.3){\line(1,0){2}}
\put(24,11.3){\line(1,0){2}}
\put(27,11.3){\line(1,0){2}}
\put(30,11.3){\line(1,0){2}}
\put(34.2,11.3){\line(1,0){0.8}}
\put(36.2,11.3){\line(1,0){0.8}}
\put(38.2,11.3){\line(1,0){0.8}}
\put(40.2,11.3){\line(1,0){0.8}}
\put(41,11.3){\line(0,1){0.8}}
\put(41,13.3){\line(0,1){0.8}}
\put(41,15.3){\line(0,1){0.8}}
\put(41,17.3){\line(0,1){0.8}}
\put(41,19.3){\line(0,1){0.8}}
\put(41,21.3){\line(0,1){0.8}}
\put(41,23.3){\line(0,1){0.8}}
\put(41,25.3){\line(0,1){0.8}}
\put(41,27.3){\line(0,1){0.8}}
\put(41,29.3){\line(0,1){0.8}}
\put(3,11.3){\line(0,1){2.3}}
\put(3,14.6){\line(0,1){2.3}}
\put(3,17.9){\line(0,1){2.3}}
\put(3,21.2){\line(0,1){2.3}}
\put(3,24.5){\line(0,1){2.3}}
\put(3,27.8){\line(0,1){2.3}}
\put(2,33){\makebox(0,0)[cc]{\small $\alpha(\mu)$}}
\put(40,33){\makebox(0,0)[cc]{\small $\omega(\mu)$}}
\put(5,23){\makebox(0,0)[cc]{\small $p$}}
\put(39,23){\makebox(0,0)[cc]{\small $q$}}
\put(44,10){\line(0,1){20}}
\put(44,10){\line(1,0){40}}
\put(84,10){\line(0,1){20}}
\put(85,6.5){\makebox(0,0)[cc]{\small $z^1$}}
\put(92,19.5){\makebox(0,0)[cc]{\small $f(z^0)$}}
\put(54,10){\makebox(0,0)[cc]{\small $)$}}
\put(54,6){\makebox(0,0)[cc]{\small $b$}}
\put(74,10){\makebox(0,0)[cc]{\small $($}}
\put(74,6){\makebox(0,0)[cc]{\small $f(a)$}}
\put(64,2){\makebox(0,0)[ct]{\large $\underbrace{\hspace{40mm}}_{[f(\mu)]}$}}
\put(45,11.3){\line(1,0){2}}
\put(48,11.3){\line(1,0){2}}
\put(51,11.3){\line(1,0){2}}
\put(54,11.3){\line(1,0){2}}
\put(57,11.3){\line(1,0){2}}
\put(60,11.3){\line(1,0){2}}
\put(63,11.3){\line(1,0){2}}
\put(66,11.3){\line(1,0){2}}
\put(69,11.3){\line(1,0){2}}
\put(72,11.3){\line(1,0){2}}
\put(45,11.3){\line(0,1){2.3}}
\put(45,14.6){\line(0,1){2.3}}
\put(45,17.9){\line(0,1){2.3}}
\put(45,21.2){\line(0,1){2.3}}
\put(45,24.5){\line(0,1){2.3}}
\put(45,27.8){\line(0,1){2.3}}
\put(51.5,23){\makebox(0,0)[cc]{\small $[f(p)]$}}
\put(86,10){\line(0,1){20}}
\put(86,10){\line(1,0){40}}
\put(126,10){\line(0,1){20}}
\put(96,10){\makebox(0,0)[cc]{\small $)$}}
\put(96,6){\makebox(0,0)[cc]{\small $f(b)$}}
\put(106,2){\makebox(0,0)[ct]{\large $\underbrace{\hspace{40mm}}_{[f^2(\mu)]}$}}
\put(41.5,9){\line(1,0){3}}
\put(41.5,11){\line(1,0){3}}
\put(41.5,9){\line(0,1){2}}
\put(44.5,9){\line(0,1){2}}
\put(83.5,16){\line(1,0){3}}
\put(83.5,18){\line(1,0){3}}
\put(83.5,16){\line(0,1){2}}
\put(86.5,16){\line(0,1){2}}
{\linethickness{1pt}
\put(32.5,10){\line(1,0){41}}
}
\end{picture}

\vspace*{10mm}
\begin{center}
Figure 12.
\end{center}

\medskip

Then $\mu'$ is the subpath of $[\mu f(\mu)]$ from $a$ to $f(a)$.
We prove that $\mu'$ is $A$-perfect. Then the statement (2) will be fulfilled.
First note that $\mu'$ has the $A$-decomposit\-ion $\mu'\equiv A(z^0)\cdot \ell$, where $\ell$ is the subpath of $[f(\mu)]$ from $b$ to $f(a)$.
We have $[f(\mu')]\equiv A(z^1)\cdot [f(\ell)]$. By Proposition~\ref{prop 8.5}.3),
the first edge of the $r$-cancelation area $A(z^1)$ lies in $H_r$. Then, in $\mu'\cdot [f(\mu')]$, the turn at the point between $\mu'$ and $[f(\mu')]$ is either mixed in $(G_{r-1},G_r)$
or an $r$-turn. In the first case, this turn is legal. In the second case it is $r$-legal, since it belongs to the $r$-legal path $[f(\mu)]$. Thus, $\mu'$ is $A$-perfect.

\medskip

For two vertices $u,v$ in a reduced edge path $\sigma$, let $l_{\sigma}(u,v)$ be the number of edges in the subpath of $\sigma$ from $u$ to $v$.

\medskip

{\it Variant 2.} Suppose that $z^0\preccurlyeq f(a)\preccurlyeq b$.


Then $l(\mu')=l_{[\mu f(\mu)]}(a,f(a))\leqslant l_{[\mu f(\mu)]}(a,b)\leqslant 2R_{\star}$, hence $\mu'$ satisfies the statement (1).

\medskip

{\it Variant 3.} Suppose that $\alpha([f(\mu)])\preccurlyeq f(a)\prec z^0$.

Then
$$\begin{array}{ll}
l(\mu') & \leqslant l_{\mu}(a,z^0)+l_{[f(\mu)]}(f(a),z^0)\vspace*{1mm}\\
& \leqslant l_{\mu}(a,z^0)+l_{[f(\mu)]}(f(a),f(z^0))\vspace*{1mm}\\
& \leqslant l_{\mu}(a,z^0)+ l_{[f(\mu)]}(f(a),z^1) +l_{[f(\mu)]}(z^1,f(z^0))\leqslant 2R_{\star}+C_{\star}.
\end{array}
$$
Hence, $\mu'$ satisfies the statement (1).

\medskip

{\it Case 2.2.} Suppose that $f(z^0)\prec z^0$.

We show how to find $\mu'$ satisfying (1).
By Lemma~\ref{Cooper_1}, for each initial subpath $\mu_j$ of $\mu$, we have
$[f(\mu_j)]\equiv\tau_j\ell_j$, where $\tau_j$ is the maximal common initial subpath of $[f(\mu_j)]$ and $[f(\mu)]$, and $l(\ell_j)\leqslant C_{\star}$.
Recall that $\mu_m\equiv E_1\dots E_m$ is the initial subpath of $\mu$ which terminates at $z^0$. Since $E_m$ is an $r$-edge and $\mu$ is $r$-legal,
$[f(\mu_m)]\equiv \tau_m$ is the initial subpath of $[f(\mu)]$ which terminates at $f(z^0)$.
Observe:

\medskip

(a) $\omega(\mu_m)\succ \omega(\tau_m)$ (since $z^0\succ f(z^0)$),

\vspace*{1mm}

(b) $\omega(\mu_n)\prec \omega(\tau_n)$ (since $\omega(\mu_n)=\alpha([f(\mu)])\prec \omega([f(\mu)])=\omega(\tau_n)$),

\vspace*{1mm}

(c) the $l$-distance in $[f(\mu)]$ between $\omega(\tau_j)$ and $\omega(\tau_{j+1})$ is at most $2C_{\star}+||f||$.

\medskip

Indeed, (c) follows from $\tau_{j+1}\ell_{j+1}=f(\mu_{j+1})=f(\mu_jE_{j+1})=\tau_j\ell_j f(E_{j+1})$.

Claims (a)-(c) imply that
there exists $m\leqslant k<n$ such that,
for all $m\leqslant j\leqslant k$, we have $\omega(\mu_j)\succ \omega(\tau_j)$, but $\omega(\mu_{k+1})\preccurlyeq \omega(\tau_{k+1})$. From this and (c), we have:

\medskip

(d) the $l$-distance in $[f(\mu)]$ between $\omega(\mu_k)$ and $\omega(\tau_k)$ is at most $2C_{\star}+||f||+1$.

\vspace*{1mm}

(e) For $m\leqslant j\leqslant k$,
we have $\widehat{f}^{\hspace*{1.5mm}j}(\mu)\equiv [E_{j+1}\dots E_n\cdot \tau_j\ell_j]$.

\medskip

We prove (e) by induction. For $j=m$ this is valid by (15.1).
Suppose (e) is valid for some $m\leqslant j <k$ and we prove it for $j+1$.
Since $\omega(\mu_j)\succ \omega(\tau_j)$, the first edge of $\widehat{f}^{\hspace*{1.5mm}j}(\mu)$ is $E_{j+1}$.
Then $\widehat{f}^{\hspace*{1.5mm}{j+1}}(\mu)\equiv [\bar{E}_{j+1}\widehat{f}^{\hspace*{1.5mm}j}(\mu)
f(E_{j+1})]\equiv [\bar{E}_{j+1}[\bar{\mu}_j\mu f(\mu_j)]f(E_{j+1})]\equiv [\bar{\mu}_{j+1}\mu f(\mu_{j+1})]
\equiv [\bar{\mu}_{j+1}\mu\tau_{j+1}\ell_{j+1}]\equiv  [E_{j+2}\dots E_n \tau_{j+1}\ell_{j+1}]$, and we are done.

\medskip


From (e) we have $\widehat{f}^{\hspace*{1.5mm}k}(\mu)\equiv [E_{k+1}\dots E_n\cdot \tau_k\ell_k]$,
and from (d) we conclude that $l(\widehat{f}^{\hspace*{1.5mm}k}(\mu))=l([E_{k+1}\dots E_n\cdot
\tau_k])+l(\ell_k)\leqslant 3C_{\star}+||f||+1\leqslant 3C_{\star}+R_{\star}+||f||$.
Thus, the statement (1) is satisfied for $\mu':=\widehat{f}^{\hspace*{1.5mm}k}(\mu)$. \hfill $\Box$


\medskip

\section{\!\! $\mu$-subgraphs in the case where $\mu$ is not $r$-legal, but\!~$r$-superstable}

Let $H_r$ be an exponential stratum.
Recall that a reduced path $\tau \subset G_r$ is called $r$-stable if all $r$-cancelation points in $\tau$
are non-deletable.
By Theorem~\ref{prop 7.2}, for any reduced edge path $\tau\subset G_r$,
one can compute a natural number $i_0$ such that $[f^{i_0}(\tau)]$ is $r$-stable.
The number of $r$-cancelation points in $[f^{i_0}(\tau)]$ is denoted by $N_r(\tau)$.
Note that $N_r(\tau)$ is the minimum of the numbers of $r$-cancelation points in the paths $[f^i(\tau)]$, $i\geqslant 0$. The $r$-superstable paths were introduced in Definition~\ref{defn 12.1}.

\begin{prop}\label{heute}  Let $H_r$ be an exponential stratum. Suppose that $\mu$ is an $r$-superstable reduced $f$-path in $G_r$ with $N_r(\mu)\geqslant 1$. Then the $\mu$-subgraph contains a vertex $\mu_1$ for which at least one of the following statements is satisfied:

\medskip

{\rm (1)}  $l(\mu_1)\leqslant C_{\star}+2R_{\star}$;

\medskip

{\rm (2)} $N_r(\mu_1)<N_r(\mu)$;

\medskip

{\rm (3)} the path $\mu_1$ in $\Gamma$ is $A$-perfect.

\medskip

Moreover, there exists $s\leqslant l(\mu)$ such that $\mu_1:= \widehat{f}^{\hspace*{1.5mm}s}(\mu)$ satisfies {\rm (1), (2)},  or {\rm (3)}.
In particular, such $\mu_1$ can be efficiently found.

\end{prop}

\medskip

{\it Proof.}
Let $y_1,\dots,y_k$ $(k\geqslant 1)$ be all $r$-cancelation points in $\mu$ (they are non-deletable, since $\mu$
is $r$-stable) and let
$y_1^1,\dots,y_k^1$ be their $1$-successors in $[f(\mu)]$.
The terminal point of $I_0:=I(\bar{\mu},[f(\mu)])$ is denoted by $z^0$.
For two vertices $v,u\in \mu$, we denote by $l_{\mu}(v,u)$ the length of the subpath of $\mu$ from $v$ to $u$.

\medskip

{\it Case 1.} Suppose that $y_k$ lies in  $I_0$.

The proof is illustrated by Figure~13, where we distinguish two cases:

\hspace*{5mm}(a) $y_k$ lies in the interior of $I_0$. In this case $y_k=y_1^1$.

\hspace*{5mm}(b) $y_k$ coincides with the terminal vertex of $I_0$.

\medskip

The point $y_1$ divides $\mu$ into two subpaths, say $p$ and $q$, so $\mu\equiv pq$.

{\it Claim.} For every decomposition $p\equiv p_1p_2$, where $p_1$ is a proper initial edge-subpath of $p$,
cancelations in the product $p_2q\cdot [f(p_1)]$ don't touch the first edge of~$p_2$.

\medskip

To prove this claim, it suffices to note that $[f(p_1)]$ is $r$-legal, but $p_2q$ is not.
The first follows from the fact that $p$, and hence $p_1$, are $r$-legal; the second follows from the fact that
$p_2q$ contains the $r$-cancelation point $y_k$.


Now we set $\mu_1:=\widehat{f}^{\hspace*{1.5mm}l(p)}(\mu)$.
Clearly, $\mu_1$ lies in the $\mu$-subgraph.
We assert that $N_r(\mu_1)< N_r(\mu)$. To prove this, we describe the path $\mu_1$ precisely.
First, by the above claim we have $\mu_1\equiv [qf(p)]$.
Second, we observe that $\mu_1\equiv \mu_{y_1y_k}\cdot \ell$,
where $\mu_{y_1y_k}$ is the initial subpath of $q$ from $y_1$ to $y_k$ and $\ell$ is the final subpath of $[f(p)]$ from $y_k$ to $f(y_1)$.
Since $[f(p)]$ and hence $\ell$ are $r$-legal, the $r$-cancelation points of $\mu_1$
are contained in the set $\{y_1,\dots,y_k\}\setminus \{y_1\}$. Hence, $N_r(\mu_1)<N_r(\mu)$.

\unitlength 1mm 
\linethickness{0.4pt}
\ifx\plotpoint\undefined\newsavebox{\plotpoint}\fi 
\begin{picture}(120,40)(-16,0)
\put(20,19){\line(1,0){3}}
\put(20,21){\line(1,0){3}}
\put(20,19){\line(0,1){2}}
\put(23,19){\line(0,1){2}}
\put(2,10){\line(1,0){18}}
\put(20,10){\line(0,1){22}}
\put(22,10){\line(0,1){22}}
\put(22,10){\line(1,0){16.5}}
{\linethickness{0.8pt}
\put(10,10){\line(1,0){10}}
\put(20,10){\line(0,1){10}}
\put(22,20){\line(1,0){10}}
}
\put(10,10){\circle*{1.00}}
\put(10.5,7){\makebox(0,0)[cc]{\small $y_1$}}
\put(20,20){\circle*{1.00}}
\put(22,20){\circle*{1.00}}
\put(17,22){\makebox(0,0)[cc]{\small $y_k$}}
\put(25.5,22.5){\makebox(0,0)[cc]{\small $y_1^1$}}
\put(38,20){\makebox(0,0)[cc]{\small $f(y_1)$}}
\put(32,20){\circle*{0.8}}
\put(22,20){\line(1,0){10}}
\put(11.25,4){\makebox(0,0)[ct]{\large $\underbrace{\hspace{17mm}}_{\mu}$}}
\put(30.25,4){\makebox(0,0)[ct]{\large $\underbrace{\hspace{16.5mm}}_{[f(\mu)]}$}}
\put(2,10.8){\line(1,0){2}}
\put(4.6,10.8){\line(1,0){2}}
\put(7.2,10.8){\line(1,0){2}}
\put(6,13){\makebox(0,0)[cc]{\small $p$}}
\put(22.8,31.8){\line(0,-1){4}}
\put(22.8,27){\line(0,-1){4}}
\put(22.8,22.2){\line(0,-1){2}}
\put(27,21){\line(1,0){2}}
\put(29.6,21){\line(1,0){2.3}}
\put(29,29){\makebox(0,0)[cc]{\small $[f(p)]$}}
\put(62,10){\line(1,0){18}}
\put(80,10){\line(0,1){22}}
\put(82,10){\line(0,1){22}}
\put(82,10){\line(1,0){16.5}}
\put(70,10){\circle*{1.00}}
\put(70.5,7){\makebox(0,0)[cc]{\small $y_1$}}
\put(80,10){\circle*{1.00}}
\put(82.6,31.8){\line(0,-1){4}}
\put(82.6,27){\line(0,-1){4}}
\put(82.6,22.2){\line(0,-1){4}}
\put(82.6,17.4){\line(0,-1){4}}
\put(82.6,12.8){\line(0,-1){2}}
\put(82.6,10.8){\line(1,0){2.3}}
\put(85.7,10.8){\line(1,0){2.3}}
\put(88.8,10.8){\line(1,0){2.3}}
\put(91.1,10.8){\line(0,1){2}}
\put(91.1,13.2){\line(0,1){2}}
\put(89,26){\makebox(0,0)[cc]{\small $[f(p)]$}}
\put(62,10.8){\line(1,0){2}}
\put(64.6,10.8){\line(1,0){2}}
\put(67.2,10.8){\line(1,0){2}}
\put(66,13){\makebox(0,0)[cc]{\small $p$}}
\put(80,7){\makebox(0,0)[cc]{\small $y_k$}}
\put(92.5,7){\makebox(0,0)[cc]{\small $y_1^1$}}
\put(98,18){\makebox(0,0)[cc]{\small $f(y_1)$}}
\put(92,10){\circle*{1}}
\put(92,15){\circle*{0.8}}
\put(92,10){\line(0,1){5}}
\put(71.25,4){\makebox(0,0)[ct]{\large $\underbrace{\hspace{17mm}}_{\mu}$}}
\put(90.25,4){\makebox(0,0)[ct]{\large $\underbrace{\hspace{16.5mm}}_{[f(\mu)]}$}}
{\linethickness{0.8pt}
\put(70,10){\line(1,0){10}}
\put(82,10){\line(1,0){10}}
\put(92,10){\line(0,1){5}}
}
\put(80,9){\line(1,0){3}}
\put(80,11){\line(1,0){3}}
\put(80,9){\line(0,1){2}}
\put(83,9){\line(0,1){2}}
\put(21,-10){\makebox(0,0)[ct]{(a)}}
\put(82,-10){\makebox(0,0)[ct]{(b)}}
\end{picture}

\vspace*{12mm}
\begin{center}

Figure 13.\\
\begin{minipage}[t]{112mm}
The path $\mu_1$ from $y_1$ to $f(y_1)$ is depictured by the bold line. The point $y_1$ on the left picture can lie in $I_0$.
The point $y_1^1$ on the right picture can coincide with $y_k$. The case $k=1$ is possible.

\end{minipage}
\end{center}

\medskip

{\it Case 2.} Suppose that $y_k$ does not lie in $I_0$, but $y_1^1$ lies in $I_0$. In this case
$y_1^1=z^0$, where $z^0$ is the terminal vertex of $I_0$.


\medskip

{\it Case 2.1.} Suppose that $l_{\mu}(y_1,y_1^1)>R_{\star}+C_{\star}$. The proof in this case is illustrated by Figure~14.

\medskip

Let $u$ be the terminal vertex of the $r$-cancelation area $A(y_1)$ in $\mu$.
The point $u$ divides $\mu$ into two subpaths, say $p$ and $q$, so $\mu\equiv pq$.

\medskip

By definition of $R_{\star}$, we have $l_{\mu}(y_1,u)\leqslant R_{\star}$, hence
$l_{\mu}(u,y_1^1)>C_{\star}$ and $u$ lies in $\mu$ between $y_1$ and $y_1^1$.

\medskip

{\it Claim.} For every decomposition $p\equiv p_1p_2$, where $p_1$ is a proper initial edge-subpath of $p$,
cancelations in the product $p_2q\cdot [f(p_1)]$ don't touch $p_2$.

\medskip

Indeed, since $p_1$ is an initial subpath of $\mu$, by Lemma~\ref{Cooper_1}, we have $[f(p_1)]=ab$
for some initial subpath $a$ of $[f(\mu)]$ and some path $b$ of length at most $C_{\star}$.
This and the fact that $l_{\mu}(u,y_1^1)>C_{\star}$ imply the claim.

Let $\mu_1:=\widehat{f}^{\hspace*{1.5mm}l(p)}(\mu)$. Clearly, $\mu_1$ lies in the $\mu$-subgraph.
We assert  that $N_r(\mu_1)< N_r(\mu)$. To prove this, we describe the path $\mu_1$ precisely.
First, by the above claim we have $\mu_1\equiv [qf(p)]$. Then $\mu_1\equiv \mu_{uy_1^1}\cdot \ell$, where $\ell$ is the subpath of $[f(\mu)]$ from $y_1^1=z^0$ to $f(u)$. Since $\ell$ is $r$-legal, the $r$-cancelation points of $\mu_1$ are contained in the set $\{y_1,\dots,y_k,y_1^1\}\setminus \{y_1\}$.
However, $y_1^1$ is not an $r$-cancelation point of $\mu_1$, see Case (4) in Remark~\ref{vier_cases}. Hence, $N_r(\mu_1)< N_r(\mu)$.

\unitlength 1mm 
\linethickness{0.4pt}
\ifx\plotpoint\undefined\newsavebox{\plotpoint}\fi 
\begin{picture}(50,40)(-40,-5)
\put(2,10){\line(1,0){28}}
\put(30,10){\line(0,1){18}}
\put(2,11.3){\line(1,0){2}}
\put(5,11.3){\line(1,0){2}}
\put(8,11.3){\line(1,0){2}}
\put(11,11.3){\line(1,0){2}}
\put(14,11.3){\line(1,0){2}}
\put(4,13){\makebox(0,0)[cc]{\small $p$}}
\put(39,25){\makebox(0,0)[cc]{\small $[f(p)]$}}
\put(32,10){\line(0,1){18}}
\put(33,28){\line(0,-1){2}}
\put(33,25){\line(0,-1){2}}
\put(33,22){\line(0,-1){2}}
\put(33,19){\line(0,-1){2}}
\put(33,16){\line(0,-1){2}}
\put(33,13){\line(0,-1){2}}
\put(33,11){\line(1,0){2}}
\put(36,11){\line(1,0){2}}
\put(39,11){\line(1,0){2}}
{\linethickness{0.8pt}
\put(32,10){\line(1,0){10}}
\put(16,10){\line(1,0){14}}
}
\put(32,10){\line(1,0){20.5}}
\put(16,10){\circle*{1.00}}
\put(32,10){\circle*{1.00}}
\put(16.5,7){\makebox(0,0)[cc]{\small $u$}}
%
%
\put(32.5,7){\makebox(0,0)[cc]{\small $y_1^1$}}
\put(42,10){\circle*{1.00}}
\put(42.5,7){\makebox(0,0)[cc]{\small $f(u)$}}
\put(8,10){\circle*{1.00}}
\put(8.5,7){\makebox(0,0)[cc]{\small $y_1$}}
\put(32,10){\line(2,1){5}}
\put(37,12.5){\circle*{0.8}}
\put(40,16){\makebox(0,0)[cc]{\small $f(y_1)$}}
\put(16,4){\makebox(0,0)[ct]{\large $\underbrace{\hspace{27mm}}_{\mu}$}}
\put(42.5,4){\makebox(0,0)[ct]{\large $\underbrace{\hspace{20mm}}_{[f(\mu)]}$}}
\put(29,9){\line(1,0){3}}
\put(29,11){\line(1,0){3}}
\put(29,9){\line(0,1){2}}
\put(32,9){\line(0,1){2}}
\end{picture}

\vspace*{2mm}
\begin{center}
Figure 14.

The path $\mu_1$ is depictured by the bold line from $u$ to $f(u)$.
\end{center}

\medskip

{\it Case 2.2.} Suppose that $l_{\mu}(y_1,y_1^1)\leqslant R_{\star}+C_{\star}$.

Let $p$ be the initial subpath of $\mu$ which terminates at $y_1$ (see Figure 15); then $\mu\equiv pq$ for some $q$. We set $\mu_1:=\widehat{f}^{\hspace*{1.5mm}l(p)}(\mu)$. As in Case 1, one can prove that $\mu_1\equiv [qf(p)]$, and hence
$\mu_1\equiv [\mu_{y_1y_1^1}\cdot \ell]$,
where $\mu_{y_1y_1^1}$ is the initial subpath of $q$ from $y_1$ to $y_1^1$
and $\ell$ is the final subpath of $[f(p)]$ from $y_1^1$ to $f(y_1)$.
Using Lemma~\ref{Cooper_1}, we deduce that $l(\ell)\leqslant C_{\star}$, hence
$$l(\mu_1)\leqslant l_{\mu}(y_1,y_1^1)+l(\ell)\leqslant R_{\star}+2C_{\star},$$
and we have the statement (1).

\unitlength 1mm 
\linethickness{0.4pt}
\ifx\plotpoint\undefined\newsavebox{\plotpoint}\fi 
\begin{picture}(50,40)(-40,-5)
\put(2,11.3){\line(1,0){2}}
\put(5,11.3){\line(1,0){2}}
\put(8,11.3){\line(1,0){2}}
\put(11,11.3){\line(1,0){2}}
\put(14,11.3){\line(1,0){2}}
\put(4,13){\makebox(0,0)[cc]{\small $p$}}
\put(39,25){\makebox(0,0)[cc]{\small $[f(p)]$}}
\put(2,10){\line(1,0){28}}
\put(30,10){\line(0,1){18}}
\put(32,10){\line(0,1){18}}
\put(32,10){\line(1,0){10}}
\put(32,10){\line(1,0){20.5}}
\put(16,10){\circle*{1.00}}
\put(32,10){\circle*{1.00}}
\put(16.5,7){\makebox(0,0)[cc]{\small $y_1$}}
\put(33,7){\makebox(0,0)[cc]{\small $y_1^1$}}
\put(32,10){\line(2,1){5}}
\put(37,12.5){\circle*{0.8}}
\put(40,16){\makebox(0,0)[cc]{\small $f(y_1)$}}
\put(16,4){\makebox(0,0)[ct]{\large $\underbrace{\hspace{27mm}}_{\mu}$}}
\put(42.5,4){\makebox(0,0)[ct]{\large $\underbrace{\hspace{20mm}}_{[f(\mu)]}$}}
{\linethickness{0.8pt}
\put(16,10){\line(1,0){14}}
}
\put(31.9,10){\line(2,1){5}}
\put(32.1,10){\line(2,1){5}}
\put(30,9){\line(1,0){3}}
\put(30,11){\line(1,0){3}}
\put(30,9){\line(0,1){2}}
\put(33,9){\line(0,1){2}}

\put(33,28){\line(0,-1){2}}
\put(33,25){\line(0,-1){2}}
\put(33,22){\line(0,-1){2}}
\put(33,19){\line(0,-1){2}}
\put(33,16){\line(0,-1){2}}
\put(33,13){\line(0,-1){1}}
\put(33,12){\line(1,0){0.5}}
\put(34,12.3){\line(1,0){0.5}}
\put(35,12.6){\line(1,0){0.5}}
\put(36,12.9){\line(1,0){0.5}}

\end{picture}

\vspace*{2mm}
\begin{center}
Figure 15.

The path $\mu_1$ is depictured by the bold line from $y_1$ to $f(y_1)$.
\end{center}

\bigskip

{\it Case 3.} Suppose that both, $y_k$ and $y_1^1$ don't lie in $I_0$.

\medskip

Then $y_1,\dots,y_k,y_1^1,\dots ,y_k^1$ don't lie in $I_0$.

We may assume that $z^0$ is an $r$-cancelation
point in $[\mu f(\mu)]$.

Indeed, suppose the opposite.
We set $\mu_1:=\widehat{f}^{\hspace*{1.5mm}l(p)}(\mu)$, where $p$ is the initial subpath of $\mu$ with the endpoint $y_1$ (see Figure~16). Then $\mu_1$ lies in the $\mu$-subgraph and the $r$-cancelation points of $\mu_1$ are contained in the set $\{y_2,\dots,y_k,z^0\}$. But $z^0$ is not an $r$-cancelation point in $\mu_1$, since we have supposed that $z^0$ is not an $r$-cancelation point in $[\mu f(\mu)]$.  So, $N_r(\mu_1)=N_r(\mu)-1$ and we were done.

\vspace*{10mm}
\hspace*{00mm}
{
\unitlength=1mm
\special{em:linewidth 0.4pt}
\linethickness{0.6pt}
\begin{picture}(92.00,45.00)(-20,00)
\put(0,31.3){\line(1,0){2}}
\put(3.2,31.3){\line(1,0){2}}
\put(6.4,31.3){\line(1,0){2}}
\put(9.6,31.3){\line(1,0){2}}
\put(12.8,31.3){\line(1,0){2}}
\put(4,33){\makebox(0,0)[cc]{\small $p$}}
\put(00.00,30.00){\line(1,0){40.00}}
\put(28.00,30.00){\circle*{1.00}}
\put(29,26.5){\makebox(0,0)[ct]{\small $y_k$}}
\put(23.00,30.00){\circle*{1.00}}
\put(15.00,30.00){\circle*{1.00}}
\put(16,26.5){\makebox(0,0)[ct]{\small $y_1$}}
\put(14.5,30.00){\thicklines\line(1,0){39.5}}
\put(41.00,30.00){\circle*{1.00}}
\put(41.00,26.5){\makebox(0,0)[ct]{\small $z^0$}}
\put(40.00,30.00){\line(0,1){18.00}}
\put(42.00,48.00){\line(0,-1){18.00}}
\put(49.00,46){\makebox(0,0)[cc]{\small $[f(p)]$}}
\put(42.9,48){\line(0,-1){2}}
\put(42.9,45){\line(0,-1){2}}
\put(42.9,42){\line(0,-1){2}}
\put(42.9,39){\line(0,-1){2}}
\put(42.9,36){\line(0,-1){2}}
\put(42.9,33){\line(0,-1){2}}
\put(42.9,31){\line(1,0){1.5}}
\put(45.4,31){\line(1,0){1.5}}
\put(48.4,31){\line(1,0){1.5}}
\put(51.4,31){\line(1,0){1.5}}
\put(52.9,31){\line(0,1){1.6}}
\put(52.9,33.6){\line(0,1){1.6}}
\put(42.00,30.00){\line(1,0){40.00}}
\put(70.00,30.00){\circle*{1.00}}
\put(71.5,26.5){\makebox(0,0)[ct]{\small $y_k^1$}}
\put(70.00,30.00){\line(0,1){10.00}}
\put(70.00,40.00){\circle*{0.8}}
\put(71.5,46.00){\makebox(0,0)[ct]{\small $f(y_k)$}}
\put(65.00,30.00){\circle*{1.00}}
\put(65.00,30.00){\line(0,1){6.00}}
\put(65.00,36.00){\circle*{0.8}}
\put(54.00,30.00){\circle*{1.00}}
\put(55.5,26.5){\makebox(0,0)[ct]{\small $y_1^1$}}
\put(54.00,30.00){\thicklines\line(0,1){5.00}}
\put(54.00,35.00){\circle*{0.8}}
\put(54,41){\makebox(0,0)[ct]{\small $f(y_1)$}}
\put(20.00,20.00){\makebox(0,0)[ct]{\large $\underbrace{\hspace{40mm}}_{\mu}$}}
\put(62.00,20.00){\makebox(0,0)[ct]{\large $\underbrace{\hspace{40mm}}_{[f(\mu)]}$}}
\end{picture}
}

\vspace*{-10mm}
\begin{center}
Figure 16.

The path $\mu_1$ is depictured by the bold line from $y_1$ to $f(y_1)$.
\end{center}

\medskip

Thus, we may assume that $y_1,\dots,y_k,z^0,y_1^1,\dots ,y_k^1$ are all $r$-cancelation points in $[\mu f(\mu)]$.
Since $\mu$ is $r$-superstable, they are non-deletable there.

\medskip

Let $A(y_1^1)$ be the $r$-cancelation area in $[f(\mu)]$ containing $y_1^1$.
Let $A'(z^0)$ and $A'(y_1^1)$ be the $r$-cancelation areas in $[\mu f(\mu)]$ containing $z^0$ and $y_1^1$, respectively.

{\it Claim.} We have $A(y_1^1)=A'(y_1^1)$. In particular, $\alpha(A(y_1^1))$ lies to the right from $z^0$ in $[\mu f(\mu)]$.

\medskip

Indeed, since the interiors of $A'(z^0)$ and $A'(y_1^1)$ don't intersect, $A'(y_1^1)$ lies to the right from $z^0$ in $[\mu f(\mu)]$. Hence
$A'(y_1^1)$ lies in $[f(\mu)]$. Therefore $A(y_1^1)=A'(y_1^1)$.

\medskip

The further proof is illustrated by Figure~17.

Let $u$ be the initial vertex of the $r$-cancelation area $A(y_1)$ in $\mu$. The point $u$ divides $\mu$ into two subpaths, say $p$ and $q$, so $\mu\equiv pq$. We set $\mu_1:=\widehat{f}^{\hspace*{1.5mm}l(p)}(\mu)$. Clearly, $\mu_1$ lies in the $\mu$-subgraph. As in Case 1, we can prove that $\mu_1\equiv [qf(p)]$  and hence $\mu_1$ is the subpath of $[\mu f(\mu)]$
from $\alpha(A(y_1))$ to $\alpha(A(y_1^1))$.

\vspace*{10mm}
\hspace*{00mm}
{
\unitlength=1mm
\special{em:linewidth 0.4pt}
\linethickness{0.6pt}
\begin{picture}(92.00,45.00)(-20,00)
\put(0,31.3){\line(1,0){2}}
\put(3,31.3){\line(1,0){2}}
\put(6,31.3){\line(1,0){2}}
\put(9,31.3){\line(1,0){2}}
\put(12,31.3){\line(1,0){1}}
\put(4,33){\makebox(0,0)[cc]{\small $p$}}
\put(00.00,30.00){\line(1,0){40.00}}
\put(30.00,30.00){\makebox(0,0)[cc]{$)$}}
\put(26.00,30.00){\makebox(0,0)[cc]{$($}}
\put(28.00,30.00){\circle*{1.00}}
\put(29,26.5){\makebox(0,0)[ct]{\small $y_k$}}
\put(25.5,30.00){\makebox(0,0)[cc]{$)$}}
\put(21.00,30.00){\makebox(0,0)[cc]{$($}}
\put(23.00,30.00){\circle*{1.00}}
\put(17.5,30.00){\makebox(0,0)[cc]{$)$}}
\put(13.00,30.00){\makebox(0,0)[cc]{$($}}
\put(15.00,30.00){\circle*{1.00}}
\put(16,26.5){\makebox(0,0)[ct]{\small $y_1$}}
\put(12,28){\makebox(0,0)[ct]{\small $u$}}
\put(12.5,30.00){\thicklines\line(1,0){38.5}}
\put(36.5,30.00){\makebox(0,0)[cc]{$($}}
\put(45.5,30.00){\makebox(0,0)[cc]{$)$}}
\put(41.00,30.00){\circle*{1.00}}
\put(41.00,26.5){\makebox(0,0)[ct]{\small $z^0$}}
\put(40.00,30.00){\line(0,1){18.00}}
\put(42.00,48.00){\line(0,-1){18.00}}
\put(49.00,46){\makebox(0,0)[cc]{\small $[f(p)]$}}
\put(43,48){\line(0,-1){2}}
\put(43,45){\line(0,-1){2}}
\put(43,42){\line(0,-1){2}}
\put(43,39){\line(0,-1){2}}
\put(43,36){\line(0,-1){2}}
\put(43,33){\line(0,-1){2}}
\put(43,31){\line(1,0){2}}
\put(46,31){\line(1,0){2}}
\put(49,31){\line(1,0){2}}
\put(42.00,30.00){\line(1,0){40.00}}
\put(72.5,30.00){\makebox(0,0)[cc]{$)$}}
\put(68.00,30.00){\makebox(0,0)[cc]{$($}}
\put(70.00,30.00){\circle*{1.00}}
\put(71.5,26.5){\makebox(0,0)[ct]{\small $y_k^1$}}
\put(70.00,30.00){\line(0,1){10.00}}
\put(70.00,40.00){\circle*{0.8}}
\put(71.5,46.00){\makebox(0,0)[ct]{\small $f(y_k)$}}
\put(67.5,30.00){\makebox(0,0)[cc]{$)$}}
\put(63.00,30.00){\makebox(0,0)[cc]{$($}}
\put(65.00,30.00){\circle*{1.00}}
\put(65.00,30.00){\line(0,1){6.00}}
\put(65.00,36.00){\circle*{0.8}}
\put(56.5,30.00){\makebox(0,0)[cc]{$)$}}
\put(51.5,30.00){\makebox(0,0)[cc]{$($}}
\put(54.00,30.00){\circle*{1.00}}
\put(55.5,26.5){\makebox(0,0)[ct]{\small $y_1^1$}}
\put(54.00,30.00){\line(0,1){5.00}}
\put(54.00,35.00){\circle*{0.8}}
\put(54,41){\makebox(0,0)[ct]{\small $f(y_1)$}}
\put(20.00,20.00){\makebox(0,0)[ct]{\large $\underbrace{\hspace{40mm}}_{\mu}$}}
\put(62.00,20.00){\makebox(0,0)[ct]{\large $\underbrace{\hspace{40mm}}_{[f(\mu)]}$}}
\end{picture}
}

\vspace*{-10mm}
\begin{center}
Figure 17.

The path $\mu_1$ is depictured by the bold line from $u$ to the endpoint of $[f(p)]$.
\end{center}

\medskip

We prove that $\mu_1$ is $A$-perfect.
Obviously, the $A$-decomposition of $\mu_1$ starts with the $r$-cancelation area $A(y_1)$.
It remains to prove that the turn in $\mu_1\cdot [f(\mu_1)]$ between $\mu_1$ and $[f(\mu_1)]$ is legal.
The first edge of $[f(\mu_1)]$ is the first edge of $[f(A(y_1))]=A(y_1^1)$, so this edge lies in $H_r$.
If the last edge of $\mu_1$ lies in $G_{r-1}$, then the turn is mixed and hence legal.
If the last edge of $\mu_1$ lies in $H_r$, then this turn is an $r$-turn. But this turn
lies in the $r$-legal subpath of $[\mu f(\mu)]$ from $z^0$ to $y_1^1$, hence it is $r$-legal.
Thus, $\mu_1$ is $A$-perfect and lies in the $\mu$-subgraph; the proof is completed. \hfill $\Box$

\section{Cones over vertices of $D_f$}

The main purpose of this section is to define the $\mu${\it -cone}, ${\rm \bf Cone}(\mu)$, for a vertex $\mu\in D_f$
(see Definition~\ref{12_c} and Figure~21), and to prove Corollary~\ref{corollary12.2}.
This corollary implies that if
the associated function $B_{\mu}$ is constant, then the $\mu$-cone is {\it ($f$ or $f^2$)-invariant},
i.e. for each $\sigma\in {\rm \bf Cone}(\mu)$, at least one of the vertices $[f(\sigma)]$
or $[f^2(\sigma)]$ lies in ${\rm \bf Cone}(\mu)$.

\begin{notation}\label{not12.0}
{\rm
1) For any vertex $\tau\in D_f$, there is a unique path in $D_f$ from $\tau$ to $[f(\tau)]$ with the label $\tau$.
We denote this path by $p_{\tau}$.

\medskip

2)
If $p$ is a reduced path in $D_f$ from $\tau$ to $\sigma$ with the label $E_1E_2\dots E_k$,
then there is a reduced path in $D_f$ from $[f(\tau)]$ to $[f(\sigma)]$ with the label  $[f(E_1E_2\dots E_k)]$. We denote this path by $f_{\bullet}(p)$.
}
\end{notation}

\begin{prop}\label{prop12.4} If $\tau$ is an alive vertex in $D_f$, then the first edge of $p_{\tau}$ coincides
with the first edge of the $\tau$-subgraph.
\end{prop}

{\it Proof.} Let $\tau\equiv E_1E_2\dots E_k$. Then the path $p_{\tau}\subset D_f$ starts at $\tau$
and the labels of consecutive edges of $p_{\tau}$ are $E_1,\dots ,E_k$.
The $\tau$-subgraph also starts at $\tau$, and the first edge of the $\tau$-subgraph is
labeled by $E_1$. Therefore the first edges coincide.~$\Box$

\vspace*{17mm}
{\hspace*{-45mm}
\unitlength 1mm 
\linethickness{0.4pt}
\ifx\plotpoint\undefined\newsavebox{\plotpoint}\fi 
\begin{picture}(71.25,31)(0,30)

\put(81.25,46){\makebox(0,0)[cc]{$\tau$}}
\put(81.25,49.25){\circle*{1.00}}
\put(81.25,49.25){\line(1,0){10}}
\put(85,49.25){\makebox(0,0)[cc]{$\vartriangleright$}}
\put(91.25,49.25){\circle*{1.00}}
\put(91.25,49.25){\line(1,0){10}}
\put(95,49.25){\makebox(0,0)[cc]{$\vartriangleright$}}
\put(101.25,49.25){\circle*{1.00}}
\put(101.25,49.25){\line(1,0){10}}
\put(105,49.25){\makebox(0,0)[cc]{$\vartriangleright$}}
\put(111.25,49.25){\circle*{1.00}}
\put(111.25,49.25){\line(1,0){10}}
\put(115,49.25){\makebox(0,0)[cc]{$\vartriangleright$}}
\put(121.25,49.25){\circle*{1.00}}
\put(121.25,49.25){\line(1,0){10}}
\put(125,49.25){\makebox(0,0)[cc]{$\vartriangleright$}}
\put(131.25,49.25){\circle*{1.00}}
\put(131.25,49.25){\line(1,0){10}}
\put(135,49.25){\makebox(0,0)[cc]{$\vartriangleright$}}
\put(141.25,49.25){\circle*{1.00}}
\put(146,49.25){\makebox(0,0)[cc]{\small $\dots$}}

\put(101.25,49.25){\line(1,1){21}}
\put(108.25,56.25){\circle*{1.00}}
\put(115.25,63.25){\circle*{1.00}}
\put(122.25,70.25){\circle*{1.00}}

\put(124,72){\makebox(0,0)[lc]{\small $[f(\tau)]$}}

\put(81.25,50.75){\circle*{0.5}}
\put(83,50.75){\circle*{0.5}}
\put(84.75,50.75){\circle*{0.5}}
\put(86.5,50.75){\circle*{0.5}}
\put(88.25,50.75){\circle*{0.5}}
\put(90,50.75){\circle*{0.5}}
\put(91.75,50.75){\circle*{0.5}}
\put(93.5,50.75){\circle*{0.5}}
\put(95.25,50.75){\circle*{0.5}}
\put(97,50.75){\circle*{0.5}}
\put(98.75,50.75){\circle*{0.5}}
\put(100.5,50.75){\circle*{0.5}}
\put(101.7,51.95){\circle*{0.5}}
\put(102.9,53.15){\circle*{0.5}}
\put(104.1,54.35){\circle*{0.5}}

\put(105.3,55.65){\circle*{0.5}}
\put(106.5,56.85){\circle*{0.5}}
\put(107.7,58.05){\circle*{0.5}}

\put(108.9,59.25){\circle*{0.5}}
\put(110.1,60.45){\circle*{0.5}}
\put(111.3,61.65){\circle*{0.5}}

\put(112.5,62.85){\circle*{0.5}}
\put(113.7,64.05){\circle*{0.5}}
\put(114.9,65.25){\circle*{0.5}}

\put(116.1,66.45){\circle*{0.5}}
\put(117.3,67.65){\circle*{0.5}}
\put(118.5,68.85){\circle*{0.5}}

\put(119.7,69.95){\circle*{0.5}}

\put(105,61){\makebox(0,0)[lc]{\small $p_{\tau}$}}

\end{picture}
}

\vspace*{-16mm}
\begin{center}
Figure 18.
\end{center}


\medskip

{\bf Remark.} The second edge of the $\tau$-subgraph and the second edge of the path $p_{\tau}$ can be different.
Indeed, the second vertex of the $\tau$-subgraph is $[\bar{E}_1\tau f(E_1)]=[E_2E_3\dots E_kf(E_1)]$.
The label of the second edge of the $\tau$-subgraph is the first edge of the $f$-path $[\bar{E}_1\tau f(E_1)]$.
This label is not necessarily $E_2$; however the label of the second edge of $p_{\tau}$ is $E_2$.

\medskip

\begin{defn}\label{12_a}
{\rm A path $p$ in $D_f$ is called {\it directed} if either $p$ is a vertex in $D_f$ or
$p\equiv E_1E_2\dots E_n$, where $E_1,\dots ,E_n$ are edges in~$D_f$
and the preferable direction at $\alpha(E_i)$ is the direction of $E_i$ for each $i=1,2,\dots ,n$.

The set of all repelling vertices of $D_f$ is denoted by ${\text{\rm \bf Rep}}(D_f)$.
}
\end{defn}


\begin{lem}\label{lem12.3}
Let $\mu$ and $\sigma$ be two vertices in $D_f$. Suppose that $p$ is a reduced path in $D_f$ from $\sigma$ to a vertex $\tau$ in the $\mu$-subgraph and that the inverse to the last edge of $p$ does not lie in the $\mu$-subgraph. Then there exists a vertex $\sigma_1\in \{\sigma\}\cup\, {\text{\rm \bf Rep}}(D_f)$ such that $\sigma_1$ lies in $p$ and the subpath of $p$ from $\sigma_1$ to $\tau$ is directed.
\end{lem}


{\it Proof.} We may assume that $p$ is nontrivial. Let $p\equiv E_1E_2\dots E_n$. If $E_n$ is a repelling edge, we  take $\sigma_1:=\omega(E_n)$.
If not, then the preferable direction at $\alpha(E_n)$ is carried by $E_n$ (see Figure 19).
If $E_{n-1}$ is a repelling edge, we take $\sigma_1:=\omega(E_{n-1})$.
If not, then the preferable direction at $\alpha(E_{n-1})$ is carried by $E_{n-1}$.
Continuing by induction, we complete the proof. \hfill $\Box$

\medskip

\vspace*{30mm}
{\hspace*{-45mm}
\unitlength 1mm 
\linethickness{0.4pt}
\ifx\plotpoint\undefined\newsavebox{\plotpoint}\fi 
\begin{picture}(71.25,31)(0,30)

\put(81.25,46){\makebox(0,0)[cc]{$\mu$}}
\put(101.25,46){\makebox(0,0)[cc]{$\tau$}}
\put(81.25,49.25){\circle*{1.00}}
\put(81.25,49.25){\line(1,0){10}}
\put(85,49.25){\makebox(0,0)[cc]{$\vartriangleright$}}
\put(91.25,49.25){\circle*{1.00}}
\put(91.25,49.25){\line(1,0){10}}
\put(95,49.25){\makebox(0,0)[cc]{$\vartriangleright$}}
\put(101.25,49.25){\circle*{1.00}}
\put(101.25,49.25){\line(1,0){10}}
\put(105,49.25){\makebox(0,0)[cc]{$\vartriangleright$}}
\put(111.25,49.25){\circle*{1.00}}
\put(111.25,49.25){\line(1,0){10}}
\put(115,49.25){\makebox(0,0)[cc]{$\vartriangleright$}}
\put(121.25,49.25){\circle*{1.00}}
\put(121.25,49.25){\line(1,0){10}}
\put(125,49.25){\makebox(0,0)[cc]{$\vartriangleright$}}
\put(131.25,49.25){\circle*{1.00}}
\put(131.25,49.25){\line(1,0){10}}
\put(135,49.25){\makebox(0,0)[cc]{$\vartriangleright$}}
\put(141.25,49.25){\circle*{1.00}}
\put(146,49.25){\makebox(0,0)[cc]{\small $\dots$}}

\put(101.25,59.25){\circle*{1.00}}
\put(101.25,49.25){\line(0,1){10}}
\put(101.3,55.5){\makebox(0,0)[cc]{$\triangledown$}}

\put(101.25,69.25){\circle*{1.00}}
\put(101.25,59.25){\line(0,1){10}}
\put(101.3,65.5){\makebox(0,0)[cc]{$\triangledown$}}

\put(101.25,79.25){\circle*{1.00}}
\put(101.25,69.25){\line(0,1){10}}

\put(101.25,89.25){\circle*{1.00}}
\put(101.25,79.25){\line(0,1){10}}

\put(101.25,79.25){\line(1,0){5}}
\put(104.3,79.25){\makebox(0,0)[cc]{$\vartriangleright$}}

\put(103,89){\makebox(0,0)[lc]{\small $\sigma$}}
\put(103,69){\makebox(0,0)[lc]{\small $\sigma_1$}}


\end{picture}
}

\vspace*{-12mm}
\begin{center}
Figure 19.
\end{center}

\begin{defn}\label{12_b}
{\rm
Let $\mu$ and $\sigma$ be two vertices in $D_f$ and let $k\geqslant 1$ and $\ell\geqslant 0$ be integer numbers.
We say that $\sigma$ is {\it $(k,\ell)$-close} to $\mu$
if there exist a vertex $\tau$ in~$D_f$, a directed path of length $k$ from $\mu$ to $\tau$,
and a directed path of length $\ell$ from $\sigma$ to $\tau$ such that the last edges of these paths do not coincide.


Such $\tau$ is called a {\it $\sigma$-entrance} into the $\mu$-subgraph.
}
\end{defn}


{\bf Remark.} The vertex $\sigma$ can be $(k,\ell)$-close to $\mu$ for several pairs $(k,\ell)$. Also there can be several
$\sigma$-entrances into the $\mu$-subgraph for given $\sigma$ and $\mu$.

\vspace*{-17mm}
{\hspace*{-37mm}
\unitlength 1mm 
\linethickness{0.4pt}
\ifx\plotpoint\undefined\newsavebox{\plotpoint}\fi 
\begin{picture}(71.25,61)(0,30)

\put(81.25,46){\makebox(0,0)[cc]{$\mu$}}
\put(101.25,46){\makebox(0,0)[cc]{$\tau$}}
\put(111.25,46){\makebox(0,0)[cc]{$\sigma$}}
\put(81.25,49.25){\circle*{1.00}}
\put(81.25,49.25){\line(1,0){10}}
\put(85,49.25){\makebox(0,0)[cc]{$\vartriangleright$}}
\put(91.25,49.25){\circle*{1.00}}
\put(91.25,49.25){\line(1,0){10}}
\put(95,49.25){\makebox(0,0)[cc]{$\vartriangleright$}}
\put(101.25,49.25){\circle*{1.00}}
\put(101.25,49.25){\line(1,0){10}}
\put(105,49.25){\makebox(0,0)[cc]{$\vartriangleright$}}
\put(111.25,49.25){\circle*{1.00}}
\put(111.25,49.25){\line(1,0){10}}
\put(115,49.25){\makebox(0,0)[cc]{$\vartriangleright$}}
\put(121.25,49.25){\circle*{1.00}}
\put(121.25,49.25){\line(0,1){10}}
\put(121.3,52.25){\makebox(0,0)[cc]{$\vartriangle$}}
\put(121.25,59.25){\circle*{1.00}}
\put(121.25,59.25){\line(0,1){10}}
\put(121.3,62.25){\makebox(0,0)[cc]{$\vartriangle$}}
\put(121.25,69.25){\circle*{1.00}}
\put(101.25,69.25){\circle*{1.00}}
\put(101.25,69.25){\line(1,0){10}}
\put(108,69.25){\makebox(0,0)[cc]{$\vartriangleleft$}}
\put(111.25,69.25){\circle*{1.00}}
\put(111.25,69.25){\line(1,0){10}}
\put(118,69.25){\makebox(0,0)[cc]{$\vartriangleleft$}}
\put(91.25,49.25){\circle*{1.00}}
\put(91.25,49.25){\line(1,0){10}}
\put(95,49.25){\makebox(0,0)[cc]{$\vartriangleright$}}
\put(101.25,59.25){\circle*{1.00}}
\put(101.25,49.25){\line(0,1){10}}
\put(101.35,56.25){\makebox(0,0)[cc]{$\triangledown$}}
\put(101.25,69.25){\circle*{1.00}}
\put(101.25,59.25){\line(0,1){10}}
\put(101.35,66.25){\makebox(0,0)[cc]{$\triangledown$}}
\end{picture}
}

\vspace*{-14mm}
\begin{center}
Figure 20.
\end{center}

For example, in Figure 20, $\sigma$ is simultaneously $(3+8k,0)$-close and $(2,7+8k)$-close to $\mu$ for each $k=0,1,\dots$. Moreover, $\sigma$ and $\tau$ are $\sigma$-entrances into the $\mu$-subgraph. If the $\mu$-subgraph
is a ray and $\sigma\in D_f^0$ is given, then the set of $\sigma$-entrances into the $\mu$-subgraph is either empty or consists of a single element. We stress that $k\geqslant 1$ in the above definition.


\medskip

\begin{defn}\label{12_c}
{\rm Given two vertices $\mu, \sigma\in D_f^0$, we define ${\rm \bf Entr}(\sigma,\mu)$ as the set of all
$\sigma$-entrances into the $\mu$-subgraph. The set $${\rm \bf Cone}(\mu)=\{\sigma\in D_f^0\,|\,{\rm \bf Entr}(\sigma,\mu)\neq \emptyset \}$$
is called the {\it $\mu$-cone}. For each $k\geqslant 1$, the set
$${\rm \bf Cone}^{(k)}(\mu)=\{\sigma\in {\rm \bf Cone}(\mu)\,|\, \sigma\hspace*{2mm} {\text {\rm is}}\hspace*{2mm} (k,\ell){\text{\rm -close to}}\hspace*{2mm} \mu\hspace*{2mm}{\text{\rm  for some}}
\hspace*{2mm} \ell\geqslant 0\}$$
is called the {\it $k$-branch} of the $\mu$-cone.
}
\end{defn}

Observe that if $\sigma\in {\rm \bf Cone}(\mu)$ and $\sigma$ is not dead, then $\widehat{f}(\sigma)\in {\rm \bf Cone}(\mu)$. The $\mu$-cone is the union of its $\mu$-branches.

\bigskip

\vspace*{30mm}
{\hspace*{-45mm}
\unitlength 1mm 
\linethickness{0.4pt}
\ifx\plotpoint\undefined\newsavebox{\plotpoint}\fi 
\begin{picture}(71.25,31)(0,30)

\put(81.25,46){\makebox(0,0)[cc]{$\mu$}}
\put(81.25,49.25){\circle*{1.00}}
\put(81.25,49.25){\line(1,0){10}}
\put(85,49.25){\makebox(0,0)[cc]{$\vartriangleright$}}
\put(91.25,49.25){\circle*{1.00}}
\put(91.25,49.25){\line(1,0){10}}
\put(95,49.25){\makebox(0,0)[cc]{$\vartriangleright$}}
\put(101.25,49.25){\circle*{1.00}}
\put(101.25,49.25){\line(1,0){10}}
\put(105,49.25){\makebox(0,0)[cc]{$\vartriangleright$}}
\put(111.25,49.25){\circle*{1.00}}
\put(111.25,49.25){\line(1,0){10}}
\put(115,49.25){\makebox(0,0)[cc]{$\vartriangleright$}}
\put(121.25,49.25){\circle*{1.00}}
\put(121.25,49.25){\line(1,0){10}}
\put(125,49.25){\makebox(0,0)[cc]{$\vartriangleright$}}
\put(131.25,49.25){\circle*{1.00}}
\put(131.25,49.25){\line(1,0){10}}
\put(135,49.25){\makebox(0,0)[cc]{$\vartriangleright$}}
\put(141.25,49.25){\circle*{1.00}}
\put(146,49.25){\makebox(0,0)[cc]{\small $\dots$}}

\put(101.25,59.25){\circle*{1.00}}
\put(101.25,49.25){\line(0,1){10}}
\put(101.3,55.5){\makebox(0,0)[cc]{$\triangledown$}}

\put(101.25,69.25){\circle*{1.00}}
\put(101.25,59.25){\line(0,1){10}}
\put(101.3,65.5){\makebox(0,0)[cc]{$\triangledown$}}

\put(101.25,79.25){\circle*{1.00}}
\put(101.25,69.25){\line(0,1){10}}
\put(101.3,75.5){\makebox(0,0)[cc]{$\triangledown$}}

\put(101.25,89.25){\circle*{1.00}}
\put(101.25,79.25){\line(0,1){10}}
\put(101.3,85.5){\makebox(0,0)[cc]{$\triangledown$}}

\put(111.25,79.25){\circle*{1.00}}
\put(101.25,79.25){\line(1,0){10}}
\put(107.5,79.25){\makebox(0,0)[cc]{$\vartriangleleft$}}

\put(91.25,79.25){\circle*{1.00}}
\put(91.25,79.25){\line(1,0){10}}
\put(95,79.25){\makebox(0,0)[cc]{$\vartriangleright$}}

\put(121.25,59.25){\circle*{1.00}}
\put(121.25,49.25){\line(0,1){10}}
\put(121.3,55.5){\makebox(0,0)[cc]{$\triangledown$}}

\put(121.25,69.25){\circle*{1.00}}
\put(121.25,59.25){\line(0,1){10}}
\put(121.3,65.5){\makebox(0,0)[cc]{$\triangledown$}}

\put(121.25,79.25){\circle*{1.00}}
\put(121.25,69.25){\line(0,1){10}}
\put(121.3,75.5){\makebox(0,0)[cc]{$\triangledown$}}

\put(121.25,89.25){\circle*{1.00}}
\put(121.25,79.25){\line(0,1){10}}
\put(121.3,85.5){\makebox(0,0)[cc]{$\triangledown$}}

\put(131.25,59.25){\circle*{1.00}}
\put(121.25,59.25){\line(1,0){10}}
\put(127.5,59.25){\makebox(0,0)[cc]{$\vartriangleleft$}}

\put(131.25,69.25){\circle*{1.00}}
\put(131.25,59.25){\line(0,1){10}}
\put(131.3,65.5){\makebox(0,0)[cc]{$\triangledown$}}

\put(131.25,79.25){\circle*{1.00}}
\put(131.25,69.25){\line(0,1){10}}
\put(131.3,75.5){\makebox(0,0)[cc]{$\triangledown$}}

\put(141.25,69.25){\circle*{1.00}}
\put(131.25,69.25){\line(1,0){10}}
\put(137.5,69.25){\makebox(0,0)[cc]{$\vartriangleleft$}}

\end{picture}
}

\vspace*{-12mm}
\begin{center}
Figure 21.
\end{center}


\begin{center}
\begin{minipage}{11cm}
All vertices on this figure, except of $\mu$, lie in ${\rm \bf Cone}(\mu)$.
\end{minipage}
\end{center}

\medskip

\begin{notation}
{\rm Let $\mu,\tau$ be vertices in $D_f$. We write $\mu\preccurlyeq \tau$ if there exists $i\geqslant 0$ with $\widehat{f}^{\hspace*{1.5mm i}}(\mu)=\tau$.
We also write $\mu\prec\tau$ if there exists $i>0$ with $\widehat{f}^{\hspace*{1.5mm i}}(\mu)=\tau$. Clearly, the relation $\preccurlyeq$ is reflexive and transitive, but
not necessarily antisymmetric (see Figure~20).}
\end{notation}

Recall that each vertex in $D_f$ is an $f$-path in $\Gamma$, hence it has an $l$-length.
The following proposition implies that ${\rm \bf Cone}(\mu)$ is ($f$ or $f^2$)-invariant
if all vertices of the $\mu$-subgraph have $l$-lengths
larger than $C_{\star}$ and if the $\mu$-cone does not contain repelling vertices.

\medskip

\begin{prop}\label{prop 11.5}  
Let $\mu$ be a vertex in $D_f$, $\sigma\in {\rm \bf Cone}(\mu)$ and $\tau\in {\rm \bf Entr}(\sigma,\mu)$.\break
If $l(\tau)> C_{\star}$, then there exist
$\sigma_1\in {\rm \bf Cone}(\mu)$ and $\tau_1\in {\rm \bf Entr}(\sigma_1,\mu)$ such that

\begin{enumerate}
\item[1)] $\sigma_1\in \{[f(\sigma)],[f^2(\sigma)]\}\cup {\text{\rm \bf Rep}}(D_f)$;
\item[2)] $\tau\preccurlyeq\tau_1$;
\item[3)] there are directed paths  from $\tau$ to $\tau_1$
and from $\sigma_1$ to $\tau_1$ of lengths at most $$\max\{ l(\sigma)+{\text{\rm dist}}(\sigma,\tau),\,\\
l(\tau)+l([f(\sigma)])+C_{\star}\},$$
where ${\text{\rm dist}}(\sigma,\tau)$ is the number of edges in the shortest directed path from $\sigma$ to $\tau$.
\end{enumerate}
\end{prop}


{\it Proof.} We consider the path $p_{\sigma}$ from $\sigma$ to $[f(\sigma)]$ with the label $\sigma$.
Note that $f_{\bullet}(p_{\sigma})$ is the path from $[f(\sigma)]$ to $[f^2(\sigma)]$ with the label $[f(\sigma)]$.
Also we consider the shortest directed path $p_{\sigma,\tau}$ from $\sigma$ to $\tau$.

{\it Case 1.} Suppose that the maximal common initial subpath of $p_{\sigma,\tau}$ and $p_{\sigma}$
is a
\hspace*{18.5mm} proper subpath of $p_{\sigma,\tau}$ (see Figure~22).


We write $p_{\sigma,\tau}\equiv p_1p_2$ and $p_{\sigma}\equiv p_1p_3$, where $p_1$ is the maximal common initial subpath of $p_{\sigma,\tau}$ and $p_{\sigma}$. By Lemma~\ref{lem12.3}, there exists a vertex $\sigma_1\in \{[f(\sigma)]\}\cup {\text{\rm \bf Rep}}(D_f)$ in $\bar{p}_3$, such that the subpath of $\bar{p}_3$ from $\sigma_1$ to the end of $\bar{p}_3$ is directed. Let $\bar{p}_4$ be this subpath.
Then $\bar{p}_4p_2$ is the directed path
from $\sigma_1$ to $\tau$, and we can set $\tau_1:=\tau$. In this case the statements 1) and 2) are evident. The statement 3) is valid since ${\text{\rm dist}}(\sigma_1,\tau_1)\leqslant l(\bar{p}_4p_2)\leqslant l(p_{\sigma})+{\text{\rm dist}}(\sigma,\tau)=l(\sigma)+{\text{\rm dist}}(\sigma,\tau)$.

\vspace*{-4mm}

{\hspace*{-42mm}
\unitlength 1mm 
\linethickness{0.4pt}
\ifx\plotpoint\undefined\newsavebox{\plotpoint}\fi 
\begin{picture}(71.25,61)(0,30)

\put(81.25,46){\makebox(0,0)[cc]{$\mu$}}
\put(102.25,46){\makebox(0,0)[cc]{$\tau=\tau_1$}}
\put(99,81){\makebox(0,0)[cc]{$\sigma$}}
\put(81.25,49.25){\circle*{1.00}}
\put(81.25,49.25){\line(1,0){10}}
\put(85,49.25){\makebox(0,0)[cc]{$\vartriangleright$}}
\put(91.25,49.25){\circle*{1.00}}
\put(91.25,49.25){\line(1,0){10}}
\put(95,49.25){\makebox(0,0)[cc]{$\vartriangleright$}}
\put(101.25,49.25){\circle*{1.00}}
\put(101.25,49.25){\line(1,0){10}}
\put(105,49.25){\makebox(0,0)[cc]{$\vartriangleright$}}
\put(111.25,49.25){\circle*{1.00}}
\put(111.25,49.25){\line(1,0){10}}
\put(115,49.25){\makebox(0,0)[cc]{$\vartriangleright$}}
\put(121.25,49.25){\circle*{1.00}}
\put(121.25,49.25){\line(1,0){10}}
\put(125,49.25){\makebox(0,0)[cc]{$\vartriangleright$}}
\put(131.25,49.25){\circle*{1.00}}
\put(131.25,49.25){\line(1,0){10}}
\put(135,49.25){\makebox(0,0)[cc]{$\vartriangleright$}}
\put(141.25,49.25){\circle*{1.00}}
\put(146,49.25){\makebox(0,0)[cc]{\small $\dots$}}
\put(101.25,59.25){\line(1,1){21}}
\put(102.25,61.25){\circle*{0.5}}
\put(103.25,62.25){\circle*{0.5}}
\put(104.25,63.25){\circle*{0.5}}
\put(105.25,64.25){\circle*{0.5}}
\put(106.25,65.25){\circle*{0.5}}
\put(107.25,66.25){\circle*{0.5}}
\put(108.25,67.25){\circle*{0.5}}
\put(109.25,68.25){\circle*{0.5}}
\put(110.25,69.25){\circle*{0.5}}
\put(111.25,70.25){\circle*{0.5}}
\put(112.25,71.25){\circle*{0.5}}
\put(113.25,72.25){\circle*{0.5}}
\put(114.25,73.25){\circle*{0.5}}
\put(115.25,74.25){\circle*{0.5}}
\put(116.25,75.25){\circle*{0.5}}
\put(117.25,76.25){\circle*{0.5}}
\put(118.25,77.25){\circle*{0.5}}
\put(119.25,78.25){\circle*{0.5}}
\put(120.25,79.25){\circle*{0.5}}
\put(121.25,80.25){\circle*{0.5}}
%
%
\put(102.25,61.25){\circle*{0.5}}
\put(102.25,62.25){\circle*{0.5}}
\put(102.25,63.25){\circle*{0.5}}
\put(102.25,64.25){\circle*{0.5}}
\put(102.25,65.25){\circle*{0.5}}
\put(102.25,66.25){\circle*{0.5}}
\put(102.25,67.25){\circle*{0.5}}
\put(102.25,68.25){\circle*{0.5}}
\put(102.25,69.25){\circle*{0.5}}
\put(102.25,70.25){\circle*{0.5}}
\put(102.25,71.25){\circle*{0.5}}
\put(102.25,72.25){\circle*{0.5}}
\put(102.25,73.25){\circle*{0.5}}
\put(102.25,74.25){\circle*{0.5}}
\put(102.25,75.25){\circle*{0.5}}
\put(102.25,76.25){\circle*{0.5}}
\put(102.25,77.25){\circle*{0.5}}
\put(102.25,78.25){\circle*{0.5}}
\put(102.25,79.25){\circle*{0.5}}
\put(102.25,59.25){\circle*{0.5}}
\put(103.25,60.25){\circle*{0.5}}
\put(104.25,61.25){\circle*{0.5}}
\put(105.25,62.25){\circle*{0.5}}
\put(106.25,63.25){\circle*{0.5}}
\put(107.25,64.25){\circle*{0.5}}
\put(108.25,65.25){\circle*{0.5}}
\put(109.25,66.25){\circle*{0.5}}
\put(110.25,67.25){\circle*{0.5}}
\put(111.25,68.25){\circle*{0.5}}
\put(112.25,69.25){\circle*{0.5}}
\put(102.25,58.25){\circle*{0.5}}
\put(102.25,57.25){\circle*{0.5}}
\put(102.25,56.25){\circle*{0.5}}
\put(102.25,55.25){\circle*{0.5}}
\put(102.25,54.25){\circle*{0.5}}
\put(102.25,53.25){\circle*{0.5}}
\put(102.25,52.25){\circle*{0.5}}
\put(102.25,51.25){\circle*{0.5}}
\put(102.25,50.25){\circle*{0.5}}
\put(102.25,49.25){\circle*{0.5}}

\put(122.5,80.55){\circle*{1.00}}

\put(124,82){\makebox(0,0)[lc]{\small $[f(\sigma)]$}}
\put(113.5,69){\makebox(0,0)[lc]{\small $\sigma_1$}}
\put(112.75,76.75){\makebox(0,0)[lc]{\small $p_3$}}
\put(103.75,71.75){\makebox(0,0)[lc]{\small $p_1$}}
\put(103.75,54.75){\makebox(0,0)[lc]{\small $p_2$}}
\put(108,62.75){\makebox(0,0)[lc]{\small $p_4$}}
\put(111.5,69.5){\circle*{1.00}}

\put(101.25,59.25){\circle*{1.00}}
\put(101.25,49.25){\line(0,1){10}}
\put(101.3,55.5){\makebox(0,0)[cc]{$\triangledown$}}

\put(101.25,69.25){\circle*{1.00}}
\put(101.25,59.25){\line(0,1){10}}
\put(101.3,65.5){\makebox(0,0)[cc]{$\triangledown$}}

\put(101.25,79.25){\circle*{1.00}}
\put(101.25,69.25){\line(0,1){10}}
\put(101.3,75.5){\makebox(0,0)[cc]{$\triangledown$}}

\end{picture}
}

\vspace*{-12mm}
\begin{center}
Figure 22.
\end{center}

\medskip

{\it Case 2.} Suppose that $p_{\sigma,\tau}$ is an initial segment of $p_{\sigma}$.


Then the vertex $\tau$ lies on the reduced path $p_{\sigma}$.
By Lemma~\ref{lem 3.12},
$[f(\tau)]$ lies in the $C_{\star}$-neighborhood of some vertex $\nu$ on the reduced path $f_{\bullet}(p_{\sigma})$
from $[f(\sigma)]$ to $[f^2(\sigma)]$ (see Figure~23).
Let $p_1$ be a path from $[f(\tau)]$ to $\nu$ of length at most~$C_{\star}$.
Recall that $p_{\tau}$ is the reduced path from $\tau$ to $[f(\tau)]$ with the label $\tau$. Hence, $l(p_{\tau})=l(\tau)>C_{\star}$ by assumption.
By Proposition~\ref{prop12.4}, the first edge of $p_{\tau}$ carries the preferable direction at $\tau$.

Let $q_{\tau,\nu}$ be the reduced form of the path ${p}_{\tau}p_1$.
Then $q_{\tau,\nu}$ is a path from $\tau$ to $\nu$, it has length at most $l(\tau)+C_{\star}$, and
the first edge of $q_{\tau,\nu}$ coincides with the first edge of $p_{\tau}$, i.e. with the directed edge emanating from $\tau$.

The vertex $\nu$ divides the reduced path $f_{\bullet}(p_{\sigma})$ into two subpaths, say $f_{\bullet}(p_{\sigma})\equiv q_1q_2$. Then one of the paths $q_{\tau,\nu}{\bar q}_1$ or $q_{\tau,\nu}q_2$
is reduced. W.l.o.g., we may assume that $q_{\tau,\nu}q_2$ is reduced. This path terminates at $[f^2(\sigma)]$
and its first edge is the directed edge emanating from $\tau$. We write $q_{\tau,\nu}q_2\equiv qp$, where $q$ is the maximal initial subpath of $q_{\tau,\nu}q_2$ that lies in the $\tau$-subgraph.  By Lemma~\ref{lem12.3}, there
is a vertex $\sigma_1\in \{[f^2(\sigma)]\}\cup  {\text{\rm \bf Rep}}(D_f)$ such that $\sigma_1$ lies in $\bar{p}$ and the subpath of $\bar{p}$ from $\sigma_1$ to $\tau_1:=\omega(q)$ is directed.
Clearly, $\sigma_1\in {\rm \bf Cone}(\mu)$ and $\tau_1\in {\rm \bf Entr}(\sigma_1,\mu)$, and $\tau\prec \tau_1$. The statement 3) follows from the estimate
$$
\begin{array}{ll}
l(q_{\tau,\nu}q_2) & \leqslant l(p_{\tau})+l(p_1)+l(q_2)\vspace*{1mm}\\
& \leqslant l(p_{\tau})+l(p_1)+l(q_1q_2)\vspace*{1mm}\\
& \leqslant l(\tau)+ C_{\star}+l([f(\sigma)]).
\end{array}
$$
We use here that the label of $p_{\tau}$ is $\tau$ and the label of $f_{\bullet}(p_{\sigma})$ is $[f(\sigma)]$. \hfill $\Box$

\medskip

\vspace*{-7mm}
{\hspace*{-45mm}
\unitlength 1mm 
\linethickness{0.4pt}
\ifx\plotpoint\undefined\newsavebox{\plotpoint}\fi 
\begin{picture}(71.25,81)(0,30)

\put(81.25,46){\makebox(0,0)[cc]{$\mu$}}
\put(101.25,46){\makebox(0,0)[cc]{$\tau$}}
\put(99,81){\makebox(0,0)[cc]{$\sigma$}}
\put(131.75,76){\makebox(0,0)[cc]{$\nu$}}
\put(116.5,92){\makebox(0,0)[cc]{$[f(\sigma)]$}}
\put(149.5,98.5){\makebox(0,0)[cc]{$[f^2(\sigma)]$}}
\put(131,67){\makebox(0,0)[cc]{$p_1$}}
\put(137,60){\makebox(0,0)[cc]{$p_{\tau}$}}
\put(148,84){\makebox(0,0)[cc]{$q_2$}}
\put(120,78){\makebox(0,0)[cc]{$q_1$}}
\put(81.25,49.25){\circle*{1.00}}
\put(81.25,49.25){\line(1,0){10}}
\put(85,49.25){\makebox(0,0)[cc]{$\vartriangleright$}}
\put(91.25,49.25){\circle*{1.00}}
\put(91.25,49.25){\line(1,0){10}}
\put(95,49.25){\makebox(0,0)[cc]{$\vartriangleright$}}
\put(101.25,49.25){\circle*{1.00}}
\put(101.25,49.25){\line(1,0){10}}
\put(105,49.25){\makebox(0,0)[cc]{$\vartriangleright$}}
\put(111.25,49.25){\circle*{1.00}}
\put(111.25,49.25){\line(1,0){10}}
\put(115,49.25){\makebox(0,0)[cc]{$\vartriangleright$}}
\put(121.25,49.25){\circle*{1.00}}
\put(121.25,49.25){\line(1,0){10}}
\put(125,49.25){\makebox(0,0)[cc]{$\vartriangleright$}}
\put(131.25,49.25){\circle*{1.00}}
\put(131.25,49.25){\line(1,0){10}}
\put(135,49.25){\makebox(0,0)[cc]{$\vartriangleright$}}
\put(141.25,49.25){\circle*{1.00}}
\put(146,49.25){\makebox(0,0)[cc]{\small $\dots$}}

\put(121.25,49.25){\line(1,1){21}}
\put(142.25,70.25){\circle*{1.00}}

\put(144,72){\makebox(0,0)[lc]{\small $[f(\tau)]$}}

\put(101.25,59.25){\circle*{1.00}}
\put(101.25,49.25){\line(0,1){10}}
\put(101.3,55.5){\makebox(0,0)[cc]{$\triangledown$}}

\put(101.25,69.25){\circle*{1.00}}
\put(101.25,59.25){\line(0,1){10}}
\put(101.3,65.5){\makebox(0,0)[cc]{$\triangledown$}}

\put(101.25,79.25){\circle*{1.00}}
\put(101.25,69.25){\line(0,1){10}}
\put(101.3,75.5){\makebox(0,0)[cc]{$\triangledown$}}

\linethickness{0.5pt}
\qbezier(116.5,88)(133,55.75)(149.5,94.5)
\qbezier(131.5,73)(132.125,61.25)(142.25,70.5)
\put(131.5,73.25){\circle*{1.00}}
\put(116.5,88){\circle*{1.00}}
\put(149.5,94.5){\circle*{1.00}}
\put(101.25,48.25){\circle*{0.5}}
\put(102.7,48.25){\circle*{0.5}}
\put(104.15,48.25){\circle*{0.5}}
\put(105.6,48.25){\circle*{0.5}}
\put(107.05,48.25){\circle*{0.5}}
\put(108.5,48.25){\circle*{0.5}}
\put(109.95,48.25){\circle*{0.5}}
\put(111.4,48.25){\circle*{0.5}}
\put(112.85,48.25){\circle*{0.5}}
\put(114.3,48.25){\circle*{0.5}}
\put(115.75,48.25){\circle*{0.5}}
\put(117.2,48.25){\circle*{0.5}}
\put(118.65,48.25){\circle*{0.5}}
\put(120.1,48.25){\circle*{0.5}}
\put(121.55,48.25){\circle*{0.5}}
\put(122.8,49.5){\circle*{0.5}}
\put(124.05,50.75){\circle*{0.5}}
\put(125.3,52){\circle*{0.5}}
\put(126.55,53.25){\circle*{0.5}}
\put(127.8,54.5){\circle*{0.5}}
\put(129.05,55.75){\circle*{0.5}}
\put(130.3,57){\circle*{0.5}}
\put(131.55,58.25){\circle*{0.5}}
\put(132.8,59.5){\circle*{0.5}}
\put(134.05,60.75){\circle*{0.5}}
\put(135.3,62){\circle*{0.5}}
\put(136.55,63.25){\circle*{0.5}}
\put(137.8,64.5){\circle*{0.5}}
\put(139.05,65.75){\circle*{0.5}}
\put(140.3,67){\circle*{0.5}}
\put(141.55,68.25){\circle*{0.5}}
\put(142.8,69.5){\circle*{0.5}}

\end{picture}
}

\vspace*{-15mm}
\begin{center}
Figure 23.
\end{center}

\vspace*{2mm}

Proposition~\ref{prop 11.5} can be applied only if the assumption $l(\tau)>C_{\star}$ is satisfied.
In Section~18, we need to control the stronger inequality $l(\tau)> 3C_{\star}+2R_{\star}+||f||$.
Therefore we introduce a function $B_{\mu}$ which counts the number of
cases where this inequality fails.

\begin{defn} {\rm Let $\mu$ be a vertex in $D_f$.
\begin{enumerate}
\item[1)] A $k$-branch of the $\mu$-cone is called {\it thin} if it contains a vertex $\sigma\in D_f$ which is repelling
or satisfies
$l(\sigma)\leqslant 3C_{\star}+2R_{\star}+||f||$.

\item[2)] For each $n\in \mathbb{N}$, let $B_{\mu}(n)$ be the number of thin $k$-branches in the $\mu$-cone, where $1\leqslant k< n$.
\end{enumerate}
}
\end{defn}

\begin{defn}\label{zug} Let $\mathcal{B}$ be the total number of repelling vertices in $D_f$ and of $f$-paths in $\Gamma$ of $l$-length at most $3C_{\star}+2R_{\star}+||f||$.
\end{defn}


\begin{rmk}\label{remark_12}
{\rm The function $B_{\mu}$ is nondecreasing and the constant $\mathcal{B}$ is computable by Proposition~\ref{prop 4.1}.
If the $\mu$-subgraph is a ray, $B_{\mu}$ is bounded from above by the constant $\mathcal{B}$.
Therefore if there exists $n\in \mathbb{N}$ with $B_{\mu}(n)>\mathcal{B}$, then the $\mu$-subgraph is finite.
}
\end{rmk}

\begin{cor}\label{corollary12.2} Let $\mu$ be a vertex of $D_f$ and let $\sigma\in {\rm \bf Cone}(\mu)$. Suppose that we know
the numbers $k,\ell$ such that $\sigma$ is $(k,\ell)$-close to $\mu$. Then at least one of the following holds:
\begin{enumerate}
\item[i)] There exists a computable number $k_1>k$ such that $B_{\mu}(k_1)>B_{\mu}(k)$.
\item[ii)] There exist computable numbers $k',\ell'$ such that $k'\geqslant k$ and at least one of the vertices $[f(\sigma)]$, $[f^2(\sigma)]$ is $(k',\ell')$-close to $\mu$.
\end{enumerate}
Moreover, we can efficiently determine which of these cases occurs.
\end{cor}


{\it Proof.} Let $\tau:=\widehat{f}^{\hspace*{1.5mm} k}(\mu)$. If $l(\tau)\leqslant C_{\star}$, then we have i) for $k_1=k+1$.
If $l(\tau)> C_{\star}$, then there exist $\sigma_1\in {\rm \bf Cone}(\mu)$ and $\tau_1\in {\rm \bf Entr}(\sigma_1,\mu)$
satisfying the claims 1)-3) of Proposition~\ref{prop 11.5}.
By 1), we have $\sigma_1\in \{[f(\sigma)],[f^2(\sigma)]\}\cup\, {\text{\rm \bf Rep}}(D_f)$.
By 2) and 3), we can compute $k'\geqslant k$ and $\ell'$ such that $\sigma_1$ is $(k',\ell')$-close to $\mu$.
If $\sigma_1\in \{[f(\sigma)],[f^2(\sigma)]\}$, we have ii).
If $\sigma_1\in {\text{\rm \bf Rep}}(D_f)$, then the $k'$-branch of the $\mu$-cone is thin and we have i) with $k_1:=k'+1$.
\hfill $\Box$

\section{$r$-perfect and $A$-perfect vertices in $\mu$-subgraphs}

Until the end of this section we will assume that $H_r$ is a fixed exponential stratum.
For any vertex $\mu$ in $D_f$, we define the {\it $r$-cone} over $\mu$ as follows: $${\rm \bf Cone}_r(\mu):=\{{\sigma\in \rm \bf Cone}(\mu)\,|\, {\rm the}\hspace*{2mm} f {\text{\rm -path}}\hspace*{2mm} \sigma\hspace*{2mm}  {\rm lies \,\, in}\hspace*{2mm} G_r\}.$$
Observe that if $\sigma\in {\rm \bf Cone}_r(\mu)$ and $\sigma$ is not dead, then $\widehat{f}(\sigma)\in {\rm \bf Cone}_r(\mu)$.
\begin{defn}
{\rm Let $\mu$ be a vertex in $D_f$.
With each pair $(\sigma,k)\in {\rm \bf Cone}_r(\mu)\times \mathbb{N}$, where $\sigma$ is
$(k,\ell)$-close to $\mu$ for some $\ell\geqslant 0$, we associate the {\it complexity} of $(\sigma,k)$ with respect to $\mu$:
 $$Compl_{\mu}((\sigma,k))=\bigl(\mathcal{B}-B_{\mu}(k),\, N_{r}(\sigma)\bigr).$$
\noindent
Here $N_r(\sigma)$ is the number from Definition~\ref{defn 7.1}.

\medskip

We use the symbol $\preccurlyeq_{lex}$ for the usual lexicographical order on $\mathbb{Z}\times \mathbb{Z}$.
}
\end{defn}

\begin{rmk}\label{rmk_today}
{\rm
1) The first component of the complexity is bounded from above by the computable constant $\mathcal{B}$. By Remark~\ref{remark_12}, this component is nonnegative if the $\mu$-subgraph is a ray.
It can be negative if the $\mu$-subgraph is finite.
We don't claim that this component is computable.


2) The second component is computable by Theorem~\ref{prop 7.2}.
It is nonnegative and can be roughly bounded from above by a function
depending only on $k,\ell$, and $l(\mu)$:
$$N_r(\sigma)\leqslant l(\sigma)\leqslant l(\mu)+(k+\ell)(||f||+1).$$

Indeed, the first inequality is obvious. The second follows from
$\widehat{f}^k(\hspace*{1.5mm}\mu)=\widehat{f}^{\hspace*{1.5mm}\ell}(\sigma)$ by induction with the help of the inequality
$|\l(\widehat{f}(\tau))- l(\tau)|\leqslant ||f||+1$ which holds for each non-dead vertex $\tau\in D_f$.

\medskip

3) $N_r(\sigma)=N_r([f^i(\sigma)])$ for all $i\geqslant 1$. This follows directly from the definition of~$N_r$.
}
\end{rmk}

\begin{lem}\label{super} Let $\mu$ be a vertex in $D_f$ such that the $f$-path $\mu$ lies in $G_r$, and let $\sigma\in {\rm \bf Cone}_r(\mu)$. Suppose that we know
the numbers $k,\ell$ such that $\sigma$ is $(k,\ell)$-close to $\mu$. Then at least one of the following holds:
\begin{enumerate}
\item[(1)] the $\mu$-subgraph is finite;
\item[(2)] there exists an $r$-superstable $\sigma'\in {\rm \bf Cone}_r(\mu)$ such that $\sigma'$ is $(k',\ell')$-close to $\mu$
for some $k'\geqslant k$ and $\ell'\geqslant 0$, and $$Compl_{\mu}((\sigma',k'))\preccurlyeq_{lex}Compl_{\mu}((\sigma,k)).$$
\end{enumerate}
Moreover, there is an efficient algorithm deciding which of these cases occurs;
in Case {\rm (1)} it computes the vertices of the $\mu$-subgraph, and in Case {\rm (2)} it computes $\sigma',k',\ell'$.

\end{lem}

\medskip

{\it Proof.}
We apply Corollary~\ref{corollary12.2}.
In Case i), we can compute $k_1>k$ with $B_{\mu}(k_1)>B_{\mu}(k)$.
We may assume that $\widehat{f}^{\hspace*{1.5mm}k_1}(\mu)$ exists, otherwise we have (1). Set $\sigma_1:=\widehat{f}^{\hspace*{1.5mm}k_1}(\mu)$.
Clearly, $\sigma_1$ is $(k_1,0)$-close to $\mu$, hence $\sigma_1\in {\rm \bf Cone}_r(\mu)$.
Moreover, $$Compl_{\mu}((\sigma_1,k_1))\prec_{lex}Compl_{\mu}((\sigma,k)).$$
If $\sigma_1$ is $r$-superstable, we are done. If not, we restart with $\sigma:=\sigma_1$ and $(k,\ell):=(k_1,0)$. By Remark~\ref{remark_12}, if we restart more than $\mathcal{B}$ times and each time have Case~i), then the $\mu$-subgraph is finite.

Thus, we may assume that, from the beginning and in the further process, we have Case ii). We start with $\sigma_0:=\sigma$ and $(k_0,\ell_0):=(k,\ell)$, and compute the number $S$ such that $[f^S(\sigma)]$ is $r$-superstable, see Lemma~\ref{lem 12.2}. By Case ii), for $i=0,\dots,S$, we compute consequently $\sigma_{i+1}\in \{[f(\sigma_i)],[f^2(\sigma_{i})]\}$ and
numbers $k_{i+1},\ell_{i+1}$ such that $k_{i+1}\geqslant k_i$ and $\sigma_{i+1}$ is $(k_{i+1},\ell_{i+1})$-close to $\mu$.
We have $\sigma_{i+1}\in {\rm \bf Cone}_r(\mu)$.
Observe that $\sigma_S\in \{[f^S(\sigma)],\dots ,[f^{2S}(\sigma)]\}$. Hence, $\sigma_S$ is $r$-superstable.
Moreover, $N_r(\sigma)=N_r(\sigma_S)$ by definition of $N_r$.
Since $k_S\geqslant k_0=k$, we have
$$Compl_{\mu}((\sigma_S,k_S))\preccurlyeq_{lex}Compl_{\mu}((\sigma,k)),$$ and we obtain (2) with $\sigma':=\sigma_S$
and $(k',\ell'):=(k_S,\ell_S)$.
\hfill $\Box$

\medskip

\begin{prop}\label{prop 12.3} Let $H_r$ be an exponential stratum. There exists an efficient algorithm which,
given a reduced $f$-path $\mu\subset G_r$, finds a vertex $\mu'$ in the $\mu$-subgraph with one of the following properties:

\medskip

{\rm 1)} the $f$-path $\mu'$ lies in $G_{r-1}$;

{\rm 2)} $\mu'$ is $r$-perfect;

{\rm 3)} $\mu'$ is $A$-perfect;

{\rm 4)} the $\mu'$-subgraph is finite.

\medskip

Moreover, the algorithm indicates which of these cases occurs.
In Case~{\rm 4)} it efficiently computes all vertices of the $\mu$-subgraph.

\end{prop}

{\it Proof.} We may assume that $\mu$ is not a dead vertex in $D_f$. Moreover, we may assume that all vertices in ${\rm \bf Cone}_r(\mu)$ which we construct in the process below are non-dead. Otherwise we have Case 4).

Suppose that we can find a vertex $\mu''\in {\rm \bf Cone}_r(\mu)$ which satisfies one of the properties 1)-4) and that we can find numbers $k'',\ell''$ such that
$\mu''$ is $(k'',\ell'')$-close to $\mu$.
Then we will be able to find the desired $\mu'$ just by following along the $\mu''$-subgraph
until the $\mu$-subgraph and then along the $\mu$-subgraph
(see the statements (1) and~(2) in Propositions~\ref{prop 9.2} and~\ref{prop 9.4}). We show how to find such $\mu''$.

\medskip

\begin{itemize}

\item[{\it Step 1.}] We go to Step~2 with $\sigma:=\widehat{f}({\mu})$. Clearly, $\sigma\in {\rm \bf Cone}_r(\mu)$.


\medskip

\item[{\it Step 2.}] Let $\sigma\in {\rm \bf Cone}_r(\mu)$ and we know the numbers $(k,\ell)$
such that $\sigma$ is $(k,\ell)$-close to $\mu$.
By Lemma~\ref{super} we can either prove that the $\mu$-subgraph is finite, or
find an $r$-superstable vertex $\mu_0\in {\rm \bf Cone}_r(\mu)$
%
%
and numbers $k_0\geqslant k$ and $\ell_0$ such that $\mu_0$
is $(k_0,\ell_0)$-close to~$\mu$ and $$Compl_{\mu}((\mu_0,k_0))\preccurlyeq_{lex}Compl_{\mu}((\sigma,k)).$$
We may assume that the second possibility occur. If $N_r(\mu_0)=0$, we go to Step~3. If $N_r(\mu_0)\geqslant 1$, we go to Step~4.



\medskip

\item[{\it Step 3.}] Suppose that $N_r(\mu_0)=0$. Recall that $\mu_0$ is $r$-superstable.

\noindent If $N_r([\mu_0f(\mu_0)])=0$, we can find the desired $\mu''$ by Proposition~\ref{prop 10.2}.

\noindent If $N_r([\mu_0f(\mu_0)])=1$, we can apply Proposition~\ref{prop 10.3} to $\mu_0$ and find $\mu_1$ in the $\mu_0$-subgraph such that one of the following cases  holds:

\medskip

{\rm (3a)} $l(\mu_1)\leqslant 3C_{\star}+2R_{\star}+||f||$;


{\rm (3b)} the path $\mu_1$ in $\Gamma$ is $A$-perfect.

\medskip

\noindent
By Step 2 we know that $\mu_0$ is $(k_0,\ell_0)$-close to $\mu$, hence $\mu_1$ is $(k_1,\ell_1)$-close to $\mu$ for some computable $k_1\geqslant k_0$ and $\ell_1$.

In Case (3a) we have $B_{\mu}(k_1+1)>B_{\mu}(k_1)\geqslant B_{\mu}(k_0)$.
Hence, $$Compl_{\mu}((\sigma,k_1+1))\prec_{lex} Compl_{\mu}((\mu_0,k_0))$$

for $\sigma:=\widehat{f}^{\hspace*{1.5mm}k_1+1}(\mu)$.
We go to Step 2 with this $\sigma$ and $(k,\ell):=(k_1+1,0)$.

In Case (3b) we are done for $\mu'':=\mu_1$.

\medskip

\item[{\it Step 4.}] Suppose that $N_r(\mu_0)\geqslant 1$.

\noindent
By Proposition~\ref{heute},
we can find  a vertex $\mu_1$ in the $\mu_0$-subgraph such that one of the following holds:

\medskip

{\rm (4a)}  $l(\mu_1)\leqslant 3C_{\star}+2R_{\star}+||f||$;

{\rm (4b)} $N_r(\mu_1)<N_r(\mu_0)$;

{\rm (4c)} the path $\mu_1$ in $\Gamma$ is $A$-perfect.







\medskip

\noindent In Case (4c) we are done for $\mu'':=\mu_1$. In Case~(4a) we proceed as in Case~(3a).

\noindent
Consider Case (4b).
By Step 2 we know that $\mu_0$ is $(k_0,\ell_0)$-close to $\mu$, hence $\mu_1$ is $(k_1,\ell_1)$-close to $\mu$ for some computable $k_1\geqslant k_0$ and $\ell_1$.
Since the function $B_{\mu}$ is nondecreasing and  $N_r(\mu_1)<N_r(\mu_0)$,  we have
$$Compl_{\mu}((\mu_1,k_1)\prec_{lex} Compl_{\mu}((\mu_0,k_0)),$$
and we go to Step~2 with $\sigma:=\mu_1$ and $(k,\ell):=(k_1,\ell_1)$.

Note that the new pair $(\sigma,k)$ in Steps~3 and~4 has smaller complexity than the old one.
If the first component of the complexity falls more than $\mathcal{B}$ times,
then the $\mu$-subgraph is finite by Remark~\ref{remark_12}. If the first component remains unchanged,
the second component can fall only finitely many times since it is nonnegative. We can compute or estimate the second component by Remark~\ref{rmk_today}.
Thus, the process stops in a finite (and computable) number of steps.
\hfill $\Box$
\end{itemize}


\section{$E$-perfect vertices in $\mu$-subgraphs}

\medskip

Let $H_r$ be a polynomial stratum.
There exists a permutation $\sigma$ on the set of $r$-edges
and, for each $r$-edge $E$, there exists an edge path $c_E$ (which is trivial or is an edge path in $G_{r-1}$) such that $f(E)=c_E\cdot \sigma(E)\cdot \overline{c}_{\,\overline{E}}$.
Then, for each $i\geqslant 0$ and for each $r$-edge $E$, one can compute a path $c_{i,E}$ (which is trivial or is an edge path in $G_{r-1}$) such that
$f^i(E)\equiv c_{i,E}\cdot \sigma^i(E) \cdot \overline{c}_{i,\overline{E}}$.
For any edge path $\mu\subset G_r$, let $\mathcal{N}(\mu)$ be
the number of $r$-edges in $\mu$.
Clearly, if $\mu$ is a reduced nontrivial $f$-path in $G_r$, then $\mathcal{N}(\widehat{f}(\mu))\leqslant \mathcal{N}(\mu)$.

\begin{defn}\label{defn 13.1} {\rm Let $H_r$ be a polynomial stratum.
A vertex $\mu\in D_f$ is called {\it $E$-perfect} if $\mu\equiv E_1b_1E_2\dots E_kb_k$,
where $k\geqslant 1$, $E_1,\dots,E_k$ are $r$-edges, $b_1,\dots ,b_k$ are paths which lie in $G_{r-1}$ or trivial, and $\mathcal{N}(\mu')=\mathcal{N}(\mu)$ for every vertex $\mu'$ in the $\mu$-subgraph.
}
\end{defn}

\begin{prop}\label{prop 13.2} Let $H_r$ be a polynomial stratum.
Let $\mu\equiv E_1b_1\dots E_kb_k$ be a reduced $f$-path in $G_r$,
where $k\geqslant 1$, $E_1,\dots,E_k$ are $r$-edges, and $b_1,\dots ,b_k$ are paths which lie in $G_{r-1}$ or trivial.
For $1\leqslant j\leqslant k$ and  $i\geqslant 1$, we set
$$\begin{array}{ll}
\mu_{0,j}\equiv & [E_jb_j\dots E_kb_kf(E_1b_1\dots E_{j-1}b_{j-1})],\vspace*{2mm}\\
\mu_{i,j}\equiv & [\overline{c}_{i,E_j}\,f^i(\mu_{0,j})f(c_{i,E_j})].
\end{array}$$

The elements of the sequence $\mu_{0,1},\dots \mu_{0,k},\mu_{1,1},\dots ,\mu_{1,k},\dots $ will be denoted by
$\mu_1,\mu_2,\dots $. (Clearly, all these elements are reduced $f$-paths in $G_{r}$.)
Then the following statements are satisfied.

\begin{enumerate}
\item[(1)] $\mu$ is $E$-perfect if and only if $\mathcal{N}(\widehat{f}(\mu))=\mathcal{N}(\mu)$.

\item[(2)]
One can efficiently find a vertex in the $\mu$-subgraph which is $E$-perfect
or lies in $G_{r-1}$ (considered as an $f$-path), or is dead.

\item[(3)] If $\mu$ is $E$-perfect, then $\mu_1,\mu_2,\dots $ are all $E$-perfect vertices in the $\mu$-subgraph.
\end{enumerate}
\end{prop}


{\it Proof.} (1) If $k=1$, then $\mu$ is $E$-perfect. Suppose that $k\geqslant 2$. Then (1) follows by induction from the next claim.

{\it Claim.} The condition (a) below implies the condition (b).

\begin{enumerate}

\item[(a)] $\mathcal{N}(\mu)=\mathcal{N}(\widehat{f}(\mu))=k$.


\item[(b)] $\mathcal{N}(\mu')=\mathcal{N}(\widehat{f}(\mu'))=k$,
where $\mu':=\widehat{f}^{\hspace*{1.5mm} 1+l(b_1)}(\mu)$.
\end{enumerate}

\medskip

{\it Proof.} We have
$\widehat{f}(\mu)\equiv [b_1E_2b_2\dots E_k\cdot b_kc_{1,E_1}\cdot \sigma(E_1)\cdot \overline{c}_{1,\overline{E}_1}].$
If (a) is valid, then $b_1$ is an initial subpath of  $\widehat{f}(\mu)$ and we have
$\mu'\equiv [E_2b_2\dots E_k\cdot b_kc_{1,E_1}\cdot \sigma(E_1)\cdot \overline{c}_{1,\overline{E}_1}f(b_1)],$
hence $\mathcal{N}(\mu')=\mathcal{N}(\widehat{f}(\mu))=k$. Then $$\widehat{f}(\mu')=[b_2E_3\dots E_k\cdot b_kc_{1,E_1}\cdot \sigma(E_1)\cdot \overline{c}_{1,\overline{E}_1}f(b_1)c_{1,E_2}\cdot \sigma(E_2)\cdot \overline{c}_{1,\overline{E}_2}].$$
Suppose that (b) is not valid, i.e. $\mathcal{N}(\widehat{f}(\mu'))<k$. Then $[\sigma(E_1)\cdot \overline{c}_{1,\overline{E}_1}f(b_1)c_{1,E_2}\cdot \sigma(E_2)]$ is trivial.
This is possible only if $E_2=\overline{E}_1$ (hence $b_1$ is a loop) and $[\overline{c}_{1,\overline{E}_1}f(b_1)c_{1,\overline{E}_1}]$
is trivial. The latter is equivalent that $b_1$ is trivial. But then $[E_1b_1E_2]$ is trivial and $\mu$ is not reduced,
a contradiction.~\hfill $\Box$

(2) follows from (1), and (3) can be proved by direct computations.\hfill $\Box$

\begin{prop}\label{prop 13.3} Let $H_r$ be a polynomial stratum.
For every two $E$-perfect vertices $\mu,\tau$ in $D_f$ one can efficiently decide,
whether $\tau$ lies in the $\mu$-subgraph.

\end{prop}

{\it Proof.}
By Proposition~\ref{prop 13.2}.(3), $\tau$ lies in the $\mu$-subgraph if and only if $\tau\equiv \mu_{i,j}$
for some $i\geqslant 0$ and $1\leqslant j\leqslant  k$, where $k=\mathcal{N}(\mu)$.

Let $m$ be the number of edges in $H_r$ including the inverses. Since the filtration for $f$ is maximal, we have $\sigma^m=id$.
Then, for each $r$-edge $E$ we have $$f^m(E)\equiv c_{m,E}\cdot E\cdot \overline{c}_{m,\overline{E}}.\eqno{(19.1)}$$
Since $f$, restricted to any edge, is a piecewise-linear map, we can find a subdivision $E=E'E''$ such that $f^m(E')\equiv c_{m,E}E'$
and $f^m(E'')\equiv E''\overline{c}_{m,\overline{E}}$.
This implies $$c_{m,E}= f^m(E')\overline{E'}.\eqno{(19.2)}$$

\medskip

{\it Claim.} For any integer $a,b,s,t\geqslant 0$ and for each $r$-edge $E$ the following is satisfied:

1) $c_{a+b,E}=f^b(c_{a,E})c_{b,\sigma^a(E)}$.

\vspace*{1mm}

2) $c_{ms,E}=f^{ms}(E')\overline{E'}$.

\vspace*{1mm}

3) $c_{ms+t,E}=f^{ms+t}(E')f^t(\overline{E'})c_{t,E}$.

\medskip

{\it Proof.} 1) follows from the definition of $c_{i,E}$. From 1) and using $\sigma^m=id$, we have
$$c_{ms,E}=f^{m(s-1)}(c_{m,E})\dots f^m(c_{m,E})c_{m,E}.$$ This and the equation (19.2) imply 2).
3) follows from 1) and 2).\hfill $\Box$

Using Claim 3), we deduce
$$\begin{array}{ll}
\mu_{ms+t,j} & \equiv [\overline{c}_{ms+t,E_j}f^{ms+t}(\mu_{0,j})f(c_{{ms+t},E_j})]\vspace*{2mm}\\
& \equiv [\overline{c}_{t,E_j}f^t(E_j')\cdot f^{ms+t}(\overline{E_j'}\mu_{0,j}f(E_j'))\cdot f^{t+1}(\overline{E_j'})
f(c_{t,E_j})].
\end{array}\eqno{(19.3)}
$$
Thus, $\tau\equiv\mu_{i,j}$ for some $i\geqslant 0$ and $1\leqslant j\leqslant k$ if and only if
there exist $s\geqslant 0$, $0\leqslant t<m$, and $1\leqslant j\leqslant k$, such that
$$
[f^t(\overline{E_j'})c_{t,E_j}\tau f(\overline{c}_{t,E_j})f^{t+1}(E_j')]=
[f^{ms}(f^t(\overline{E_j'}\mu_{0,j}f(E_j')))].
$$
For fixed $t,j$, and using Corollary~\ref{cor 3.8} for $f^m$, we can efficiently decide whether there exists $s\geqslant 0$ satisfying the above equation. Hence, we can efficiently decide whether there exist $i,j$ with $\tau\equiv \mu_{i,j}$.  \hfill $\Box$










\section{Finiteness and Membership problems for $\mu$-subgraphs}

We continue to work with the PL-relative train track $f:\Gamma \rightarrow \Gamma$ satisfying\break (RTT-iv).
For such $f$ we prove Propositions~\ref{prop 14.2} and~\ref{prop 14.5} which solve the Finiteness and the Membership problems from Section~7.
This will complete the proof of the main Theorem~\ref{thm 1.1}.

\medskip

\begin{lem}~\label{lem 14.1} Let $H_r$ be an exponential stratum in $\Gamma$ with the Perron-Frobenius eigenvalue~$\lambda_r$
and let $\mu\subset G_r$ be an $r$-perfect path.
Let $\mu=\mu_0,\mu_1,\dots$ be consecutive vertices of the $\mu$-subgraph in~$D_f$.
Then, for all $i\geqslant 0$, we have $$L_r(\mu_{i+1})\geqslant L_r(\mu_i)>0$$ and
$$L_r(\mu_{m_i})=\lambda_r^iL_r(\mu),$$ where $0=m_0<m_1<m_2<\dots $ are computable numbers
from Proposition~\ref{prop 9.2}. In particular, the $\mu$-subgraph is infinite.
\end{lem}


{\it Proof.} By  Proposition~\ref{prop 9.2}.\!~(3), $\mu_i$ contains edges from $H_r$ for each $i\geqslant 0$.
Hence $L_r(\mu_i)>0$.
The other formulas follow from the statements~(3) and~(2) of Proposition~\ref{prop 9.2}. These formulas and $\lambda_r>1$ imply that the $\mu$-subgraph is infinite.~\hfill $\Box$

\medskip

\begin{prop}~\label{prop 14.2} Given a vertex $\mu$ in $D_f$, one can efficiently decide, whether the $\mu$-subgraph
is finite or not. Moreover, one can efficiently compute the vertices of the
$\mu$-subgraph if it is finite.
\end{prop}

{\it Proof.}  Let $r$ be the minimal number such that the $f$-path $\mu$ lies in $G_r$.

First suppose that $H_r$ is an exponential stratum. By Proposition~\ref{prop 12.3},
we can efficiently find a vertex $\mu'$ in the $\mu$-subgraph with one of the following properties:

\medskip

{\rm 1)} the $f$-path $\mu'$ lies in $G_{r-1}$;

{\rm 2)} $\mu'$ is $r$-perfect;

{\rm 3)} $\mu'$ is $A$-perfect;

{\rm 4)} the $\mu'$-subgraph is finite.

\medskip

Moreover, the algorithm indicates which of these cases occurs, and in Case~4) it efficiently computes the vertices of the $\mu$-subgraph. Note that in Cases~1)-3) the $\mu'$-subgraph, and hence the $\mu$-subgraph, can be finite or infinite. So, we analyze these cases.

\medskip

In Case~1) we apply induction.
In Case~2) the $\mu'$-subgraph (and hence the $\mu$-subgraph) is infinite by Lemma~\ref{lem 14.1}.


Consider Case~3). By Proposition~\ref{prop 9.4}.\!~(2), there exist natural numbers
$m_{1,1}<m_{2,1}<m_{3,1}<\dots$ such that
$\widehat{f}^{\hspace*{1.5mm}m_{i,1}}(\mu')\equiv [f^i(\mu')]$, $i> 0$.
Hence, the $\mu'$-subgraph is finite if and only if there exist $0< i<j$ such
that $[f^i(\mu')]=[f^j(\mu')]$.
This problem is efficiently decidable by Corollary~\ref{cor 3.9}. In Case 4) we are done.

\medskip

Now suppose that $H_r$ is a polynomial stratum. By Proposition~\ref{prop 13.2}.\!~(2),
we can efficiently find a vertex $\mu'$ in the $\mu$-subgraph with one of the following properties:

\medskip

{\rm 1)} the $f$-path $\mu'$ lies in $G_{r-1}$ or is trivial;

{\rm 2)} $\mu'$ is $E$-perfect.

\medskip

In Case~1) we apply induction. Consider Case~2). Let $m$ be the number of edges in $H_r$ including the inverses.
Let $\mu'\equiv E_1b_1\dots E_kb_k$, where $k\geqslant 1$, $E_1,\dots, E_k$ are\break $r$-edges, and $b_1,\dots ,b_k$ are paths which lie in $G_{r-1}$ or trivial.
By Proposition~\ref{prop 13.2}.\!~(3), the $\mu'$-subgraph contains the vertices
$\mu'_{ms,1}$, $s\geqslant 0$. Hence, the $\mu'$-subgraph is finite if and only if there exist $0\leqslant s_1<s_2$ such
that $\mu'_{ms_{_1},1}\equiv \mu'_{ms_{_2},1}$. (Recall that $\mu'_{0,1}\equiv \mu'$.)
By the formula~(19.3), this is equivalent to
$$
[f^{ms_{_1}}(\overline{E_1'}\mu' f(E_1'))]=[f^{ms_{_2}}(\overline{E_1'}\mu' f(E_1'))].
$$
The problem of existence of such $s_{_1},s_{_2}$ is efficiently decidable by Corollary~\ref{cor 3.9}.

If $H_r$ is a zero stratum, then following along the $\mu$-subgraph at most $l(\mu)$ steps, we can find a vertex $\mu'$ in the $\mu$-subgraph such that the $f$-path $\mu'$ lies in $G_{r-1}$ or is trivial. Then we apply induction.

It follows from the above proof that if the $\mu$-subgraph is finite, then we can efficiently compute its vertices.\hfill $\Box$

\medskip

\begin{prop}~\label{prop 14.3} For every two vertices $\mu$, $\tau$ in $D_f$, where $\mu$ is $r$-perfect, one can efficiently decide whether $\tau$ lies in the $\mu$-subgraph.
\end{prop}

{\it Proof.} We may assume that the $f$-path $\tau$ lies in $G_r$ (otherwise $\tau$ does not lie in the $\mu$-subgraph).
Let $\mu=\mu_0,\mu_1,\dots,$ be consecutive vertices of the $\mu$-subgraph.
Using Lemma~\ref{lem 14.1}, we can compute the minimal $i$ such that $L_r(\mu_i)> L_r(\tau)$. Then $\tau$ lies in the $\mu$-subgraph if and only if $\tau$ coincides with one of the vertices $\mu_0,\mu_1,\dots,\mu_{i-1}$.~\hfill $\Box$

\medskip

\begin{prop}~\label{prop 14.4} For every two vertices $\mu$, $\tau$ in $D_f$, where $\mu$ is $A$-perfect,
one can efficiently decide whether $\tau$ lies in the $\mu$-subgraph.
\end{prop}


{\it Proof.} Due to Proposition~\ref{prop 14.2}, we may assume that the $\mu$-subgraph is infinite.
Let $\mu\equiv A_1b_1\dots A_kb_k$ be the $A$-decomposition of $\mu$.
We use the following notation from Proposition~\ref{prop 9.4}:
$$
\begin{array}{ll}
\mu_{0,j}\equiv & [A_jb_j\dots A_kb_kf(A_1b_1\dots A_{j-1}b_{j-1})],
\vspace*{2mm}\\
\mu_{i,j}\equiv & [f^i(\mu_{0,j})],
\end{array}\eqno{(20.1)}
$$
where $1\leqslant j\leqslant k$ and $i\geqslant 1$. By  Proposition~\ref{prop 9.4}.\!~(4),
for every vertex $\sigma$ in the $\mu$-subgraph, at least one of the paths $\sigma$, $\widehat{f}(\sigma),\dots ,\widehat{f}^{\hspace*{1.5mm}l(\sigma)}(\sigma)$ coincides with $\mu_{i,j}$ for some $i,j$.

Thus, we first decide, whether one of the paths  $\tau$, $\widehat{f}(\tau),\dots ,\widehat{f}^{\hspace*{1.5mm}l(\tau)}(\tau)$ coincides with $\mu_{i,j}$ for some $i,j$.
In view of (20.1), this can be done with the help of Corollary~\ref{cor 3.8}.
If the answer is negative, then $\tau$ does not lie in the $\mu$-subgraph.
If it is positive, then we can find $t,i,j$ such that $\widehat{f}^{\hspace*{1.5mm}t}(\tau)\equiv \mu_{i,j}$.
Recall that by Proposition~\ref{prop 9.4}.\!~(2),
$\mu_{i,j}\equiv \widehat{f}^{\hspace*{1.5mm}{m_{i,j}}}(\mu)$ for computable $m_{i,j}$.
Then $\tau$ lies in the $\mu$-subgraph if and only if $m_{i,j}\geqslant t$ and $\tau\equiv \widehat{f}^{\hspace*{1.5mm}{m_{i,j}-t}}(\mu)$. \hfill $\Box$



\medskip

\begin{prop}~\label{prop 14.5} Given two vertices $\mu$, $\tau$ in $D_f$, one can efficiently decide
whether $\tau$ lies in the $\mu$-subgraph.
\end{prop}


{\it Proof.} By Proposition~\ref{prop 14.2}, we can efficiently decide whether the $\mu$-subgraph and the $\tau$-subgraph are finite or not.
If the $\mu$-subgraph is finite, we can compute all its vertices and verify, whether $\tau$ is one of them.
Suppose that the $\mu$-subgraph is infinite. Then, if the $\tau$-subgraph is finite, the vertex $\tau$ cannot
lie in the $\mu$-subgraph.
So, we may assume that the $\tau$-subgraph is also infinite.
Let $r$ be the minimal number such that the $f$-path $\mu$ lies in $G_r$.
We will use induction on $r$.

First suppose that $H_r$ is an exponential stratum. Then, by Proposition~\ref{prop 12.3}, we can efficiently
find a vertex~$\mu'$ in the $\mu$-subgraph with one of the following properties:

\medskip

{\rm 1)} the $f$-path $\mu'$ lies in $G_{r-1}$;

{\rm 2)} $\mu'$ is $r$-perfect;

{\rm 3)} $\mu'$ is $A$-perfect.

\medskip

First we check, whether $\tau$ belongs to the segment of the $\mu$-subgraph from $\mu$ to~$\mu'$. If yes, we are done. If not, we replace $\mu$ by $\mu'$ and consider the above cases.
In Case 1) we proceed by induction, in Case~2) by Proposition~\ref{prop 14.3}, and
in Case~3) by Proposition~\ref{prop 14.4}.

Now suppose that $H_r$ is a polynomial stratum. Then, by Proposition~\ref{prop 13.2}.\!~(2), we
can efficiently find a vertex~$\mu'$ in the $\mu$-subgraph with one of the following properties:

\medskip

{\rm 1)} the $f$-path $\mu'$ lies in $G_{r-1}$ or is trivial;

{\rm 2)} $\mu'$ is $E$-perfect.

\medskip

In Case~1) we proceed by induction.  Suppose we have Case~2). We may assume that the $f$-path $\tau$ lies in $G_r$, otherwise
$\tau$ does not lie in the $\mu$-subgraph.
By Proposition~\ref{prop 13.2}, there exists $k\leqslant l(\tau)$ such that either $\widehat{f}^{\hspace*{1.5mm}k}(\tau)$ lies in $G_{r-1}$ or is dead, or $\widehat{f}^{\hspace*{1.5mm}k}(\tau)$ is $E$-perfect. If $\widehat{f}^{\hspace*{1.5mm}k}(\tau)$ lies in $G_{r-1}$ or is dead, then $\tau$ does not lie in the $\mu$-subgraph. Suppose that $\widehat{f}^{\hspace*{1.5mm}k}(\tau)$ is $E$-perfect. By Proposition~\ref{prop 13.3}, we can decide, whether $\widehat{f}^{\hspace*{1.5mm}k}(\tau)$ lies in the $\mu'$-subgraph,
and hence in the $\mu$-subgraph (these sub\-graphs differ by a finite segment). If $\widehat{f}^{\hspace*{1.5mm}k}(\tau)$ does not lie in the $\mu$-subgraph, then $\tau$ does not lie in the $\mu$-subgraph. If $\widehat{f}^{\hspace*{1.5mm}k}(\tau)$ lies in the $\mu$-subgraph, say $\widehat{f}^{\hspace*{1.5mm}k}(\tau)=\widehat{f}^{\hspace*{1.5mm}t}(\mu)$, then $\tau$ lies in the $\mu$-subgraph if and only if $t\geqslant k$ and $\tau=\widehat{f}^{\hspace*{1.5mm}t-k}(\mu)$.



Finally, if $H_r$ is a zero stratum, we follow along the $\mu$-subgraph at most $l(\mu)$ steps until we arrive at a vertex $\mu'\in D_f$
which, considered as an $f$-path, lies in $G_{r-1}$. Then we apply induction.
\hfill $\Box$

\end{document}